\documentclass[a4paper,10pt]{article}

\usepackage{amsmath}
\usepackage{amssymb}
\usepackage{amsthm}
\usepackage{amsfonts}
\usepackage{xcolor}
\usepackage{graphicx}
\usepackage{mathtools}
\usepackage{enumerate}
\usepackage{multirow}
\usepackage{rotating}
\usepackage{setspace}
\usepackage{lineno}
\usepackage{url}

\theoremstyle{theorem}
\newtheorem{thm}{Theorem}

\theoremstyle{definition}

\usepackage{geometry}

\begin{document}
	
\title{\large{Optimal Control Studies on Age Structural Modeling of COVID-19 in Presence of Saturated Medical Treatment of Holling Type III }}

\vspace{0.1in}
\author{{\bf\large Bishal Chhetri$^{1,a}$,  D. K. K. Vamsi$^{*,a}$,  Carani B Sanjeevi$^{b, c}$}\hspace{2mm} \\
	{\it\small $^{a}$Department of Mathematics and Computer Science, Sri Sathya Sai Institute of Higher Learning, Prasanthi Nilayam}, \\
	{\it \small Puttaparthi, Anantapur District - 515134, Andhra Pradesh, India}\\
	{\it\small $^{b}$  Vice-Chancellor, Sri Sathya Sai Institute of Higher Learning -  SSSIHL, India}\\
	{\it\small $^{c}$ Department of Medicine, Karolinska Institute, Stockholm, Sweden }\\
	{\it\small bishalchhetri@sssihl.edu.in, dkkvamsi@sssihl.edu.in$^{*}$,}\\
	{\it\small sanjeevi.carani@sssihl.edu.in, sanjeevi.carani@ki.se}\\
	{\small $^{1}$  First Authors},
	{ \small $^{*}$ Corresponding Author}
	\vspace{1mm}
}

		\date{}
		\maketitle

\begin{abstract} \vspace{.25cm}

 COVID pandemic has catalysed the development of novel coronavirus vaccines across pharmaceutical companies and research organizations. In a situation where the vaccine is still unavailable and the disease is spreading exponentially, an effective control measures seems to be the need of an hour to at least contain the size of the disease.
In this context, age related transmissibility of COVID-19 has become a public health concern. This disease has caused the most severe health issues for adults over the age of 60 with particularly fatal results for those 80 years and old. In this situation an age structured modelling could be helpful in understanding the spread of the disease across different age groups. In this study initially, we propose an age structured model and calculate the equilibrium points and basic reproduction number. Later we propose an optimal control problem to understand the roles of treatment in controlling the epidemic. From the Stability analysis we see that the infection free equilibrium remains asymptotically stable whenever $R_0 < 1$ and as $R_0$ crosses unity we have the infected equilibrium to be stable. From the sensitivity analysis parameters $u_{11}$, $b_1$, $\beta_1$, $d_1$ and $\mu$ were found to be sensitive.  Findings from  the Optimal Control studies suggests that the infection among the adult population(age $\geq 30)$ is least considering the second control $u_{12}$ whereas, when both the controls $u_{11}$ and $u_{12}$ are considered together the infectives is minimum in case of young populations(age $ \leq 30$). The cumulative infected population reduced the maximum when the second control was considered followed by considering both the controls together. The control $u_{12}$  was effective for mild epidemic $(R_0 \in(1, 2))$ whereas control $u_{11}$ was found to be highly effective when epidemic was severe $(R_0 \in(2, 7))$ for the population of age group $(\leq 30)$. Whereas for age group $(\geq 30)$ the control $u_{12}$ was highly effective for the entire  range of basic reproduction number. The effect of saturation level in treatment is also explored numerically.

\end{abstract}

\section*{Introduction}

Coronavirus disease 2019 (COVID-19) is a contagious respiratory and vascular disease caused by severe acute respiratory syndrome coronavirus 2 (SARS-CoV-2). The Coronavirus caused disease COVID-19 has been declared a pandemic by WHO. As on 20 November 2020, around 55 million have been affected and around 1.3 million have lost their worldwide \cite{chh}. The social and economic life of people is drastically affected by this COVID-19 Pandemic.

This pandemic has catalysed the development of novel coronavirus vaccines across pharmaceutical companies and research organizations. Many potential vaccines for COVID-19 are being studied and some are under clinical trials. Drugs such as remdesivir, favipiravir, ivermectin, lopinavir/ritonavir, mRNA-1273, phase I trial (NCT04280224) and AVT technology are being used as therapeutic agents by different countries for treating Covid-19  \cite{10,8,9,tu2020review}.
 Therefore, in situation where the vaccine is still unavailable and the disease is spreading exponentially, an effective control measures has become the need of an hour to at least contain the size of the disease.
 
 In this context age related transmissibility of  COVID-19 has become a public health concern. It has caused the most severe health issues for adults over the age of 60  with particularly fatal results for those 80 years and older. This is due to the number of underlying health conditions present in older populations \cite{ch}. A study with 425 patients indicated that the median age was 59  years (range, 15 to 89) and reported different case distributions in four age groups groups: 0 – 14, 15 – 44 , 45 73 – 64 , and > 65 years \cite{bai2020presumed}. However, there is no  enough epidemiological evidence to classify the age groups in  transmission. Age  structure modelling studies for diseases such as influenza and Dengue can be found in (\cite{pongsumpun2003transmission}, \cite{kouokam2013disease}). \\

Motivated by the above, in this study, we have proposed a non-linear age structured compartmental model in which the  population
is divided into two age groups.  The first group with age $\leq$ 30(we call this group by young population) and the other group with age$\geq$ 30(we call this group by adult population) and in each of these age groups we have three compartments, namely susceptible denoted by $S_1$(young) and $S_2$(adult), infected denoted by $I_1$(young) and $I_2$(adult) and recovered denoted by $R_1$(young) and $R_2$(adult). Later, we frame an optimal control problem to understand the role of treatments by considering  them as controls in the two age groups considered. 
 
 We have used Holling type III recovery rate function of infected adult individuals wherein the treatment provided initially is less, owing to less infected individuals, and as the epidemic
progresses, the treatment increases accordingly and non-linearly saturates due to limitations on medical facilities.

The section-wise split-up of the paper is as follows: In the next section we talk about the formulation of age structured model followed by the Positivity and Boundedness of the model. In section 3 the equilibrium points and the basic reproduction number is calculated following which the stability of equilibrium is analysed. In section 5 we perform the sensitivity analysis of the parameters of the model. We formulate an optimal control problem and do the simulations to see the role of the controls in reducing  the infection  in section 6. Finally, we present the discussions and conclusions.
\newpage
\section{Model Formulation}

We have divided the entire population into two age groups, one with age $\leq 30 $(we call this group as young population) and the other ones with age $\geq 30$(we call this group as adult population). In each of these groups we  have three compartments namely susceptible($S_1, S_2$), Infected($I_1, I_2$) and the Recovered($R_1, R_2$). At any time an individual can be in one of the states namely
susceptible, infected, or removed. Susceptibles are the ones who are not infected but can get infected when it comes in contact with he infected ones. Infected are those who spread the infection to others and ones an infected individual recovers from the disease then he develops immunity against the disease and  is said to be recovered. The parameters $\beta_1, \beta_2, \beta_3, \beta_4$ are the infection rates corresponding to two age groups considered. We assume that there is a constant birth rate of the young susceptible population denoted by parameter $b_1$ and also the recovered can again become susceptible again with rate $\delta_1$ and $\delta_2$. Transition rate or maturation rate describing transition from age  group to another is used and is denoted by parameters $m$. We have assumed that the natural mortality rate is common across all the compartments and is denoted by parameter $\mu$, whereas the disease induced death rate differs among age groups denoted by $d_1$ and $d_2$. As the recovery from a disease differs among different age groups we denote the recovery rate for the first age group by $u_{11}$ and for the second group by $u_{12}$. This
model leads to the system of ordinary differential equations
as follows:
\begin{eqnarray}
   	\frac{dS_{1}}{dt}& =&  b_{1}+ \delta_{1}R_{1} \ - \beta_{1} S_{1}I_{1} - \beta_{2} S_{1}I_{2} -\mu S_{1} -m S_{1} \label{sec2equ1} \\
   	\frac{dI_{1}}{dt} &=&  \beta_{1} S_{1}I_{1} + \beta_{2} S_{1}I_{2} \ - d_{1}I_{1}   \ - \mu I_{1}-\mu_{11}I_{1}   \label{sec2equ2}\\ 
   	\frac{dR_{1}}{dt} &=&   \mu_{11}I_{1} \ -  \mu R_{1}-\delta_{1}R_{1} -m R_{1} \label{sec2equ3}\\
   	\frac{dS_{2}}{dt}& =& m S_{1} + \delta_{2}R_{2} \ - \beta_{3} S_{2}I_{1} - \beta_{4} S_{2}I_{2} -\mu S_{2}  \label{sec2equ4} \\
   	\frac{dI_{2}}{dt} &=&  \beta_{3} S_{2}I_{1} + \beta_{4} S_{2}I_{2} \ - d_{2}I_{2}   \ - \mu I_{2}-\frac{\mu_{12}I_{2}^2}{1+\alpha I_{2}^2}   \label{sec2equ5}\\ 
   	\frac{dR_{2}}{dt} &=&  m R_{1}+\frac{\mu_{12}I_{2}^2}{1+\alpha I_{2}^2}  \ -  \mu R_{2}-\delta_{2}R_{2}  \label{sec2equ6}
   \end{eqnarray} 
   
    \begin{table}[ht!]
	\caption{Parameters and their Meanings} 
     	
     	\centering 
     	\begin{tabular}{|l|l|} 
     		\hline\hline
     		
     		\textbf{Parameters} &  \textbf{Biological Meaning} \\  
     		\hline\hline 
     		$S_1$ & Suceptible young population  \\
     			\hline\hline 
     		$S_2$ & Suceptible adult population  \\

     			\hline\hline 
     		$I_1$ & Infected young population  \\
     			\hline\hline 
     		$I_2$ & Infected adult population  \\
     		\hline\hline
     		$R_1$ & Recovered young population  \\
     			\hline\hline 
     		$R_2$ & Recovered adult population  \\
     		\hline\hline
     		$b_{1}$ & Constant birth rate of young population  \\
     		\hline\hline
     		$\delta_{1}$ & rate of recovered  young becoming suceptible  \\
     		\hline\hline
     		$\beta_{1}, \beta_{2}$ & Rate at which Suceptible young population are infected  \\
     		& because of infected young ones and infected adults \\
     		
     		\hline\hline
     		$\mu$ & natural death rate \\
     	
     		\hline\hline
     		$d_{1}, d_{2}$ & disease induced death rate of young and adult popultaion \\
     		 
     		\hline\hline
     		
     		$\mu_{11}$ & recovery rate of young ones because of treatment \\
     	\hline\hline
     		
     		$\delta_{2}$ & rate of recovered adult becoming suceptible again\\
     		   
     		\hline\hline

     		$\beta_{3}, \beta_{4}$ & Rate at which Suceptible adult population are infected  \\
     		& because of infected young ones and infected adults  \\
     		\hline\hline
     	$ \mu_{12},   \alpha$ & $\mu_{12}$ is the treatment rate of infected adult individuals \\
     	& and $ \alpha$ is the non-negative
constant related with saturation in
treatment.\\
     		\hline\hline
     		
     $m$ & maturation rate 	\\
     \hline\hline
     	\end{tabular}
     \end{table} \vspace{.25cm}

  \newpage
\section{Positivity and Boundedness}
     
      We now show that if the initial conditions of the system (\ref{sec2equ1})-(\ref{sec2equ6}) are positive, then the solution remain positive for any future time. Using the  equations (\ref{sec2equ1})-(\ref{sec2equ6}), we get,
\begin{align*}
\frac{dS_1}{dt} \bigg|_{S_1=0} &= \bigg(b_{1}S_{2}+\delta_{1}R_{1} \bigg)\geq 0 ,  &  
\frac{dI_1}{dt} \bigg|_{I_1=0} &= \beta_{2} S_1 I_2  \geq 0,
\\ \\
\frac{dR_1}{dt} \bigg|_{R_1=0} &= \mu_{11}I_1 \geq 0 , &
\frac{dS_2}{dt} \bigg|_{S_2=0} &= \bigg(m S_{1}+\delta_{2}R_{2} \bigg)\geq 0 
\\ \\
\frac{dI_2}{dt} \bigg|_{I_2=0} &= \beta_{3} S_2 I_1  \geq 0,&
\frac{dR_2}{dt} \bigg|_{R_2=0} &= m R_1 + \frac{\mu_{12}I_{2}^2 }{1+\alpha I_{2}^2}\geq 0 &
\end{align*}

\vspace{1.5mm}
\noindent
\\ 
Thus all the above rates are non-negative on the bounding planes (given by $S_1=0$, $I_1=0$, $R_1=0$, $S_2=0$,  $I_2=0$, and $R_2=0$) of the non-negative region of the real space. So, if a solution begins in the interior of this region, it will remain inside it throughout time $t$. This  happens because the direction of the vector field is always in the inward direction on the bounding planes as indicated by the above inequalities. Hence, we conclude that all the solutions of the the system (\ref{sec2equ1})-(\ref{sec2equ6}) remain positive for any time $t>0$  provided that the initial conditions are positive. This establishes the positivity of the solutions of the system (\ref{sec2equ1})-(\ref{sec2equ6}). Next we will show that the solution is bounded.  \vspace{.25cm}

\underline{\textbf{Boundedness}}\textbf{:}
Let  $N(t) = S_1(t)+I_1(t)+R_1(t) + S_2(t)+I_2(t)+R_2(t) $ \\

Now,  
\begin{equation*}
\begin{split}
\frac{dN}{dt} & = \frac{dS_1}{dt} +  \frac{dS_2}{dI_1}{dt}+  \frac{dR_1}{dt}+ \frac{dS_2(t)}{dt}+\frac{d I_2(t)}{dt}+\frac{d R_2(t)}{dt}  \\[4pt]
& \le b_1-\mu(S_1+I_1+R_1+S_2+I_2+R_2)  \\
\end{split}
\end{equation*}

Here the integrating factor is $e^{\mu t}.$ Therefore after integration we get, \\

$$N(t)\le \frac{b_1}{\mu}$$

Thus we have shown that the system (\ref{sec2equ1})-(\ref{sec2equ6}) is positive and bounded. Therefore the biologically feasible region is given by the following set, 
\begin{equation*}
\Omega = \bigg\{\bigg(S_1(t), I_1(t), R_1(t),S_2(t), I_2(t), R_2(t)\bigg) \in \mathbb{R}^{6}_{+} : S_1(t)  + I_1(t) + R_1(t)+S_2(t)  + I_2(t) + R_2(t) \leq \frac{b_1}{\mu}, \ t \geq 0 \bigg\}
\end{equation*}
     
       \section{Equilibrium points and Reproduction Number}
     
    The system (1.1)-(1.6) admits two equilibrium namely the disease free equilibrium and the infected equilibrium. The disease free equilibrium denoted by $E_{0}$ is calculated to be,\\
     $E_{0} = (S_1^*, S_2^*,0, 0, 0,0)$  where, \\
     $$ S_{1}^*=\frac{b_1 }{(\mu + m)}$$
     $$ S_{2}^*=\frac{b_{1} m}{\mu(\mu + m)}$$
     
     The infected equilibrium is denoted by $E_1$ and is given by,\\
     $E_1=(S_1^*,I_1^*,R_1^*,S_2^*,I_2^*,R_2^*)$ where,\\
     
     $$S_1^*=\frac{(b_1-d_1-\mu-u_{11})+\delta_1 R_1^*}{(\mu+m)}$$
     $$I_1^*=\frac{\beta_1(A+B R_1^*)+\beta_2(A+B R_1^*)}{(d_1+\mu+u_{11})}$$
     $$R_1^*=\frac{u_{11}}{\mu+\delta_1+m}$$
     $$S_2^*=\frac{m(A+B R_1^*)+\delta_2 R_2^*}{\beta_3+\beta_4 I_2^*+\mu}$$
     $$I_2^*=\frac{b_1+\delta_1 R_1^*-\beta_1(A+BR_1^*)-(\mu+m)(A+BR_1^*)}{\beta_2(A+BR_1^*)}$$
     $$R_2^*=\frac{m R_1^*+\frac{u_{12}I_2^{*2}}{1+\alpha I_2^{*2}}}{(\mu+\delta_2)}$$
     
     Here,\\
     $$A=\frac{(b_1-d_1-\mu-u_{11})}{(\mu+m)}$$
     $$B=\frac{\delta_1}{(\mu+m)}$$
    
     \subsection{\textbf{Calculation of $R_0$}}
     
     The basic reproduction number which is the average number of secondary cases produced per primary case is calculated using the next generation matrix method described in \cite{diekmann2010construction}. Our system (1.1)-(1.6) has two infected states and four uninfected states. Calculating the jacobian matrix at infection free equilibrium $E_0$ considering infected states we have,
    
\begin{equation*}
J(E_0)=    
\begin{bmatrix}
\beta_1 S_1^*-d_1-\mu-u_{11} & \beta_2 S_2^*  \\[6pt]
\beta_3 S_2^*& \beta_4 S_2^* -d_2-\mu  \\[6pt]

\end{bmatrix}
\end{equation*}      
or,
$$J(E_0) = F + V$$ where,
\begin{equation*}
F =    
\begin{bmatrix}
\beta_1 S_1^* & \beta_2 S_2^*  \\[6pt]
\beta_3 S_2^*& \beta_4 S_2^*   \\[6pt]
\end{bmatrix}
\end{equation*}

\begin{equation*}
V =    
\begin{bmatrix}
-(d_1 + \mu + u_{11}) & 0  \\[6pt]
0& -d_2-\mu  \\[6pt]
\end{bmatrix}
\end{equation*}

Calculating the inverse of $V$ we get,
\begin{equation*}
V^{-1} =    
\begin{bmatrix}
\frac{-1}{(d_1 + \mu + u_{11})} & 0  \\[6pt]
0& \frac{-1}{d_2+\mu}  \\[6pt]
\end{bmatrix}
\end{equation*}
Now
\begin{equation*}
-FV^{-1} =    
\begin{bmatrix}
\beta_1 S_1^* p & \beta_2 S_1^* q \\[6pt]
\beta_3 S_2^* p& \beta_4 S_2^* q \\[6pt]
\end{bmatrix}
\end{equation*} where,
 $$p=\frac{1}{d_1+\mu+u_{11}}$$
 $$q=\frac{1}{d_2 + \mu}$$
The characterstics equation of $-FV^{-1}$ is calculated as,
$$\lambda^2 - \bigg(\beta_1 S_1^* p + \beta_4 S_2^* q\bigg)\lambda + S_1^* S_2^* pq \bigg(\beta_1 \beta_4 -\beta_2 \beta_3\bigg)$$
The spectral radius of $-FV^{-1}$ is given by,\\
\begin{equation*}
         \ \mathbf{\frac{(\beta_{1}S_{1}^*p + \beta_{4}S_{2}^*q) + \sqrt{M}}{2}}
     \end{equation*}
     where, \\
     $$M = (\beta_{1}S_{1}^*p-\beta_{4}S_{2}^*q)^2+4S_{1^*}S_{2}^*\beta_{2}\beta_{3}pq$$ 
   
   Therefore the basic reproduction number is given by,\\
     
     \begin{equation*}
         \mathbf{ R_{0}}= \mathbf{\frac{(\beta_{1}S_{1}^*p + \beta_{4}S_{2}^*q) + \sqrt{M}}{2}}
     \end{equation*}

   \newpage   
   
    \section{Numerical Simulation}
    \textbf{Parameters Values}
       In the  following table we provide the parameter values and the source from which the parameter values are taken. Since the Severe acute respiratory syndrome coronavirus 2 (SARS-CoV-2) mimics the influenza virus regarding clinical presentation, transmission mechanism, and seasonal coincidence \cite{cuadrado2020sars} and also due to the lack of data on the age structure modelling of COVID-19, we have taken some of the parameter values from the work on influenza disease \cite{kouokam2013disease}. In similar lines to (\cite{kouokam2013disease},\cite{kumar2019role}) for numerical simulation we take the values of $\delta_1$ and $\delta_2$ and $u_{11}$ and $u_{12}$ to be same. Taking parameter values from table 2, we will numerically show the asymptotic stability of $E_0$ and $E_1$ based on the values of basic reproduction number in similar lines to \cite{kouokam2013disease}.\\

    \begin{table}[ht!]
	\caption{Parameters values and their Source} 
     	
     	\centering 
     	\begin{tabular}{|l|l|l|} 
     		\hline\hline
     		
     		\textbf{Parameters} &  \textbf{Values}&  \textbf{Source} \\  
     		\hline\hline 
     	
     		$b_{1}$ & 7.192 & \cite{samui2020mathematical} \\
     		\hline\hline
     		$\delta_{1}$ & .0714 & \cite{kouokam2013disease}  \\
     		\hline\hline
     		$\beta_{1}, \beta_{2}$ &  4/3, 2 & \cite{kouokam2013disease}\\
     		
     		\hline\hline
     		$\mu$ & .062 & \cite{samui2020mathematical}\\
     	
     		\hline\hline
     		$d_{1}, d_{2}$ & .000073, .0000913 & \cite{kouokam2013disease}\\
     		 
     		\hline\hline
     		
     		$\mu_{11}$ & .1&  -\\
     		 
     	\hline\hline
     		
     		$\delta_{2}$ & .0714 & \cite{kouokam2013disease}\\
     		   
     		\hline\hline

     		$\beta_{3}, \beta_{4}$ & 4, .00000008 & \cite{kouokam2013disease} \\
     		\hline\hline
     	$ \mu_{12} $ & .1 & \cite{kumar2019role}\\
     		\hline\hline
     	$\alpha$ & 0.4 & \cite{kumar2019role}\\
     	\hline\hline
     $m$ & 0.000182 	& 	\cite{kouokam2013disease}\\
     \hline\hline
     	\end{tabular}
     \end{table} \vspace{.25cm}

  \textbf{Stability of infection free equilibrium}\\
    We will numerically show that whenever we have $R_0 < 1$ the infection free equilibrium $E_0$ is locally asymptotically stable. For the parameter values from table 3, the value of $R_0$ was calculated as 0.98 and $E_0  = (0.1157, .00039, 0,0,0,0 )$. From figure 1 we see that the disease free equilibrium $E_0$ is locally asymptotically stable. The initial values were taken as $(S1,S2,I1,I2,R1,R2)=(100,150,5,70,10,30)$
    
    \newpage
    \begin{center}
	\begin{figure}[hbt!]
		\includegraphics[height = 9cm, width = 17.5cm,angle=0]{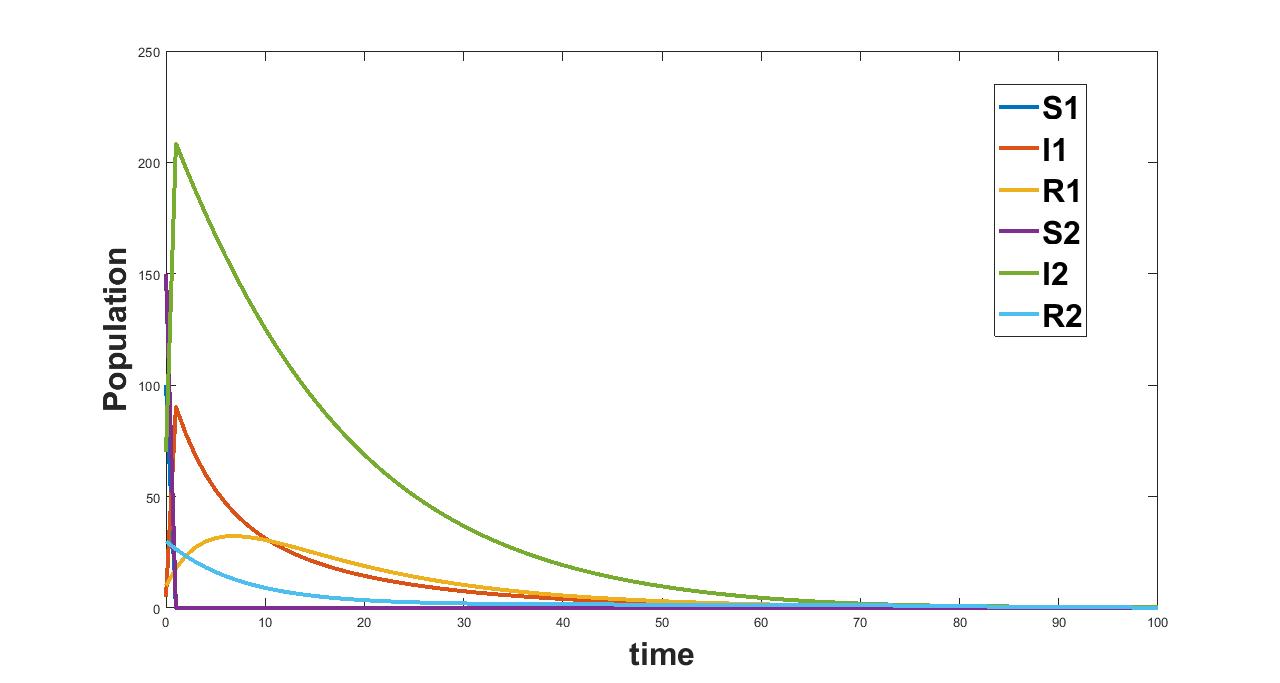} 
			\caption{Stability of $E_0$}
	
	\end{figure} 
\end{center}

\begin{table}[ht!]
	\caption{Parameters values for $R_0 < 1$} 
     	
     	\centering 
     	\begin{tabular}{|l|l|} 
     		\hline\hline
     		
     		\textbf{Parameters} &  \textbf{Values} \\  
     		\hline\hline 
     	
     		$b_{1}$ & .007192  \\
     		\hline\hline
     		$\delta_{1}$ & .0714  \\
     		\hline\hline
     		$\beta_{1}, \beta_{2}$ &  4/3, 2 \\
     		
     		\hline\hline
     		$\mu$ & .062 \\
     	
     		\hline\hline
     		$d_{1}, d_{2}$ & .000073, .0000913 \\
     		 
     		\hline\hline
     		
     		$\mu_{11}$ & .1 \\
     		 
     	\hline\hline
     		
     		$\delta_{2}$ & .0714 \\
     		   
     		\hline\hline

     		$\beta_{3}, \beta_{4}$ & 4, .00000008  \\
     		\hline\hline
     	$ \mu_{12}$ & .1 \\
     		\hline\hline
     	$	\alpha$ & 0.4  \\
     	\hline\hline
     $m$ & 0.000182 \\
     \hline\hline
     	\end{tabular}
     \end{table} \vspace{.25cm}
     
\newpage
\textbf{Stability of infected equilibrium}\\
 Here we shall numerically show that whenever the value of  basic reproduction number denoted by $R_0$ crosses unity we have the infected  equilibrium $E_1$ to be locally asymptotically stable. The value of $R_0$ was calculated to be  2.7615 and $E_1  = (10.42, 1, 0.144, 0.0041, 0.03, 0.0005)$ for the parameter values from table 4 . From figure 2 we see that the infected  equilibrium $E_1$ is locally asymptotically stable provided $R_0 > 1$. The initial values were taken as $(S_1,S_2,I_1,I_2,R_1,R_2)=(100,200,1,2,1,2)$.
    
    \begin{center}
	\begin{figure}[hbt!]
		\includegraphics[height = 9cm, width = 17.5cm]{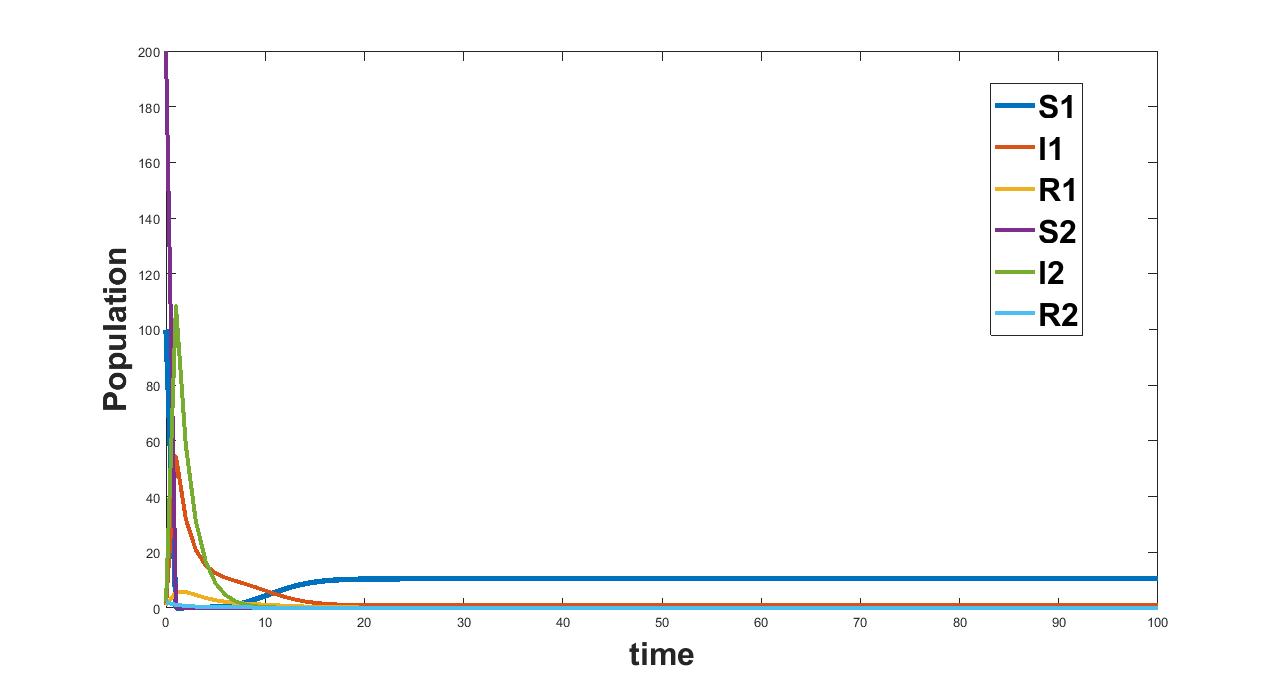} 
			\caption{Stability of $E_1$}
	
	\end{figure} 
\end{center}

\vspace{1cm}
\begin{table}[ht!]
	\caption{Parameters values for $R_0 > 1$} 
     	
     	\centering 
     	\begin{tabular}{|l|l|} 
     		\hline\hline
     		
     		\textbf{Parameters} &  \textbf{Values} \\  
     		\hline\hline 
     	
     		$b_{1}$ & 7.192  \\
     		\hline\hline
     		$\delta_{1}$ & .0714 \\
     		\hline\hline
     		$\beta_{1}, \beta_{2}$ &  .0133, 2 \\
     			\hline\hline
     		$\mu$ & .62 \\
            \hline\hline
     		$d_{1}, d_{2}$ & .000073, .0000913\\
     		\hline\hline
     			$\mu_{11}$ & .1\\
     			\hline\hline
     		$\delta_{2}$ & .0714 \\
     		   \hline\hline
     			$\beta_{3}, \beta_{4}$ & 4, .00000008  \\
     		\hline\hline
     	$ \mu_{12} , \alpha$ & .1,  .5 \\
     		\hline\hline
     		$m$ & 0.00182 	\\
     \hline\hline
     	\end{tabular}
     \end{table} \vspace{.25cm}
 
 \newpage    
 Corona virus has caused the most severe health issues for adults over the age of 60 — with particularly fatal results for those 80 years and older\cite{ch}. With the parameter values from table 2, now we will numerically illustrate this fact. From figure 3 we see that the number of infected cases in case of adult population is comparatively higher than that in case of young.

\begin{center}
	\begin{figure}[hbt!]
		\includegraphics[height = 10cm, width = 17.5cm]{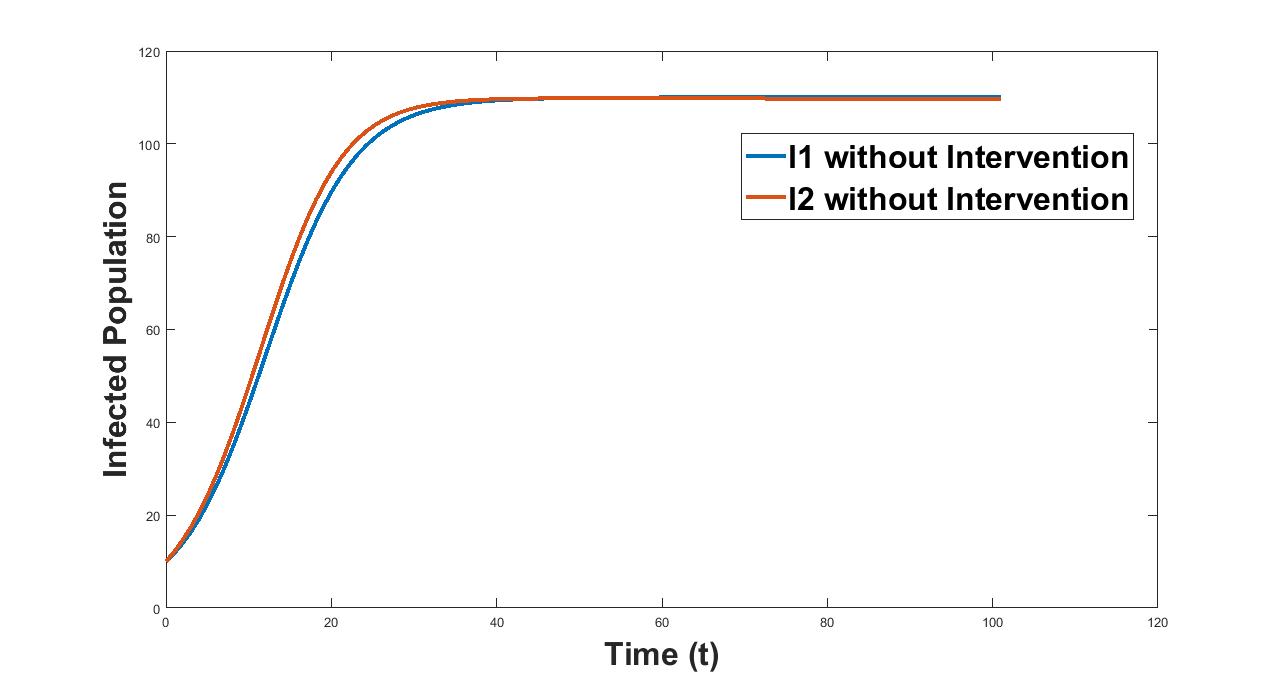} 
			\caption{Infected Populations without controls}
	
	\end{figure} 
\end{center}
\newpage
\section{Sensitivity Analysis}
In this section we check the sensitivity of all the parameters of the model (1.1)-(1.6). As each parameter is varied in different intervals, the total infected cell population, mean infected cell population and the mean square error are plotted with respect to time. These plots are used to determine the sensitivity of the parameter. The different intervals chosen are given in the following table \ref{t}. The fixed parameter values are taken from table 2. 

\begin{table}[hbt!]
		\caption{Sensitivity Analysis}
		\centering
		\label{t}
		{
			\begin{tabular}{|l|l|l|}
				\hline
				\textbf{Parameter} & \textbf{Interval} & \textbf{Step Size}  \\
				\hline
				
				$u_{11}$ & 0 to 0.5 & 0.01 
					\\ \cline{2-2}
				 & 1.5 to 2  &\\
				 \hline
				 $b_1$ & 6.5 to 7.1920 & 0.01 
					\\ \cline{2-2}
				 & 7.1920  to 8  &
				 	\\ \cline{2-2}
				 & 0.1 to 0.5  & \\
				 \hline
				 	$ m $ & 0 to 0.00182 & 0.0001 
					\\ \cline{2-3}
				 & .00182 to 1  & 0.01 \\
				 \hline
				$u_{12}$ & 0 to .5 & 0.05 
				\\ \cline{2-2} 
				& 0.5 to 2 & 
				\\ \hline
				
				$\beta_1$ & 0 to 1.33 & 0.01 
				\\ \cline{2-2}
				& 1.33 to 2 &  \\ 
				\hline 
					$\beta_2 $ & 0 to 2 & 0.01 
				\\ \cline{2-2}
				& 2 to 3 &  \\ 
				\hline 
				
					$\beta_3$ & 0 to 2.5 & 0.01 
				\\ \cline{2-2}
				& 2.5 to 5 &  \\ 
			
				\hline 
					$\beta_4$ & 0 to 0.5 & 0.01 
				\\ \cline{2-2}
				& 0.5 to 1 &  \\ 
				\hline 
				$\alpha$ & 0 to 0.5 & 0.01  
				\\ \cline{2-2}
				& 0.5 to 2 & \\
				\hline 
				$d_1$ & 0 to 0.000073  & 0.00001  
				\\ \cline{2-3}
				& 0.000073 to 1 & 0.01 \\
				\hline 
				$d_2$ & 0 to 0.0000913 & 0.00001  
				\\ \cline{2-2}
				& 0.0000913 to 2 & \\
				\hline 
				$\mu$ & 0 to 0.5 & 0.01  
				\\ \cline{2-2}
				& 0.5 to 2 & \\
				\hline
			$\delta_1$ & 0 to 0.0714 & 0.001  
				\\ \cline{2-2}
				& 0.0714 to 1 & \\
				\hline 
				$\delta_2$ & 0 to 0.0714 & 0.001  
				\\ \cline{2-2}
				& 0.0714 to 1 & \\
				\hline 
			\end{tabular}
		}
	\end{table}

	\subsubsection{Parameter ${b_1}$}

	The results related to sensitivity of $b_1$, varied in three intervals as mentioned in table \ref{t}, are given in figure (4-6). The plots of infected population for each varied value of the parameter $b_1$ per interval, the mean infected population and the mean square error are used to determine the sensitivity.
		We conclude from these plots that the parameter $b_1$ is not sensitive in interval III and sensitive in I and II. 
		
		\begin{figure}[hbt!]
			\begin{center}
				\includegraphics[width=2.2in, height=1.8in, angle=0]{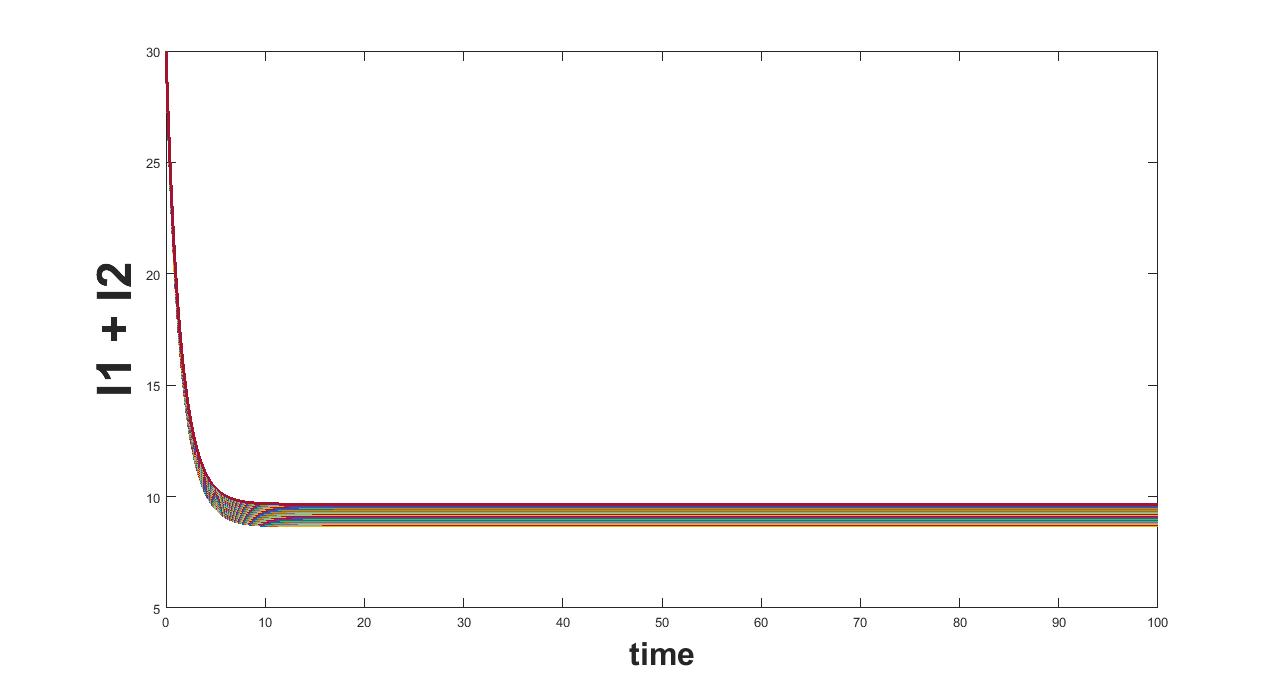}
				\hspace{-.4cm}
				\includegraphics[width=2.2in, height=1.8in, angle=0]{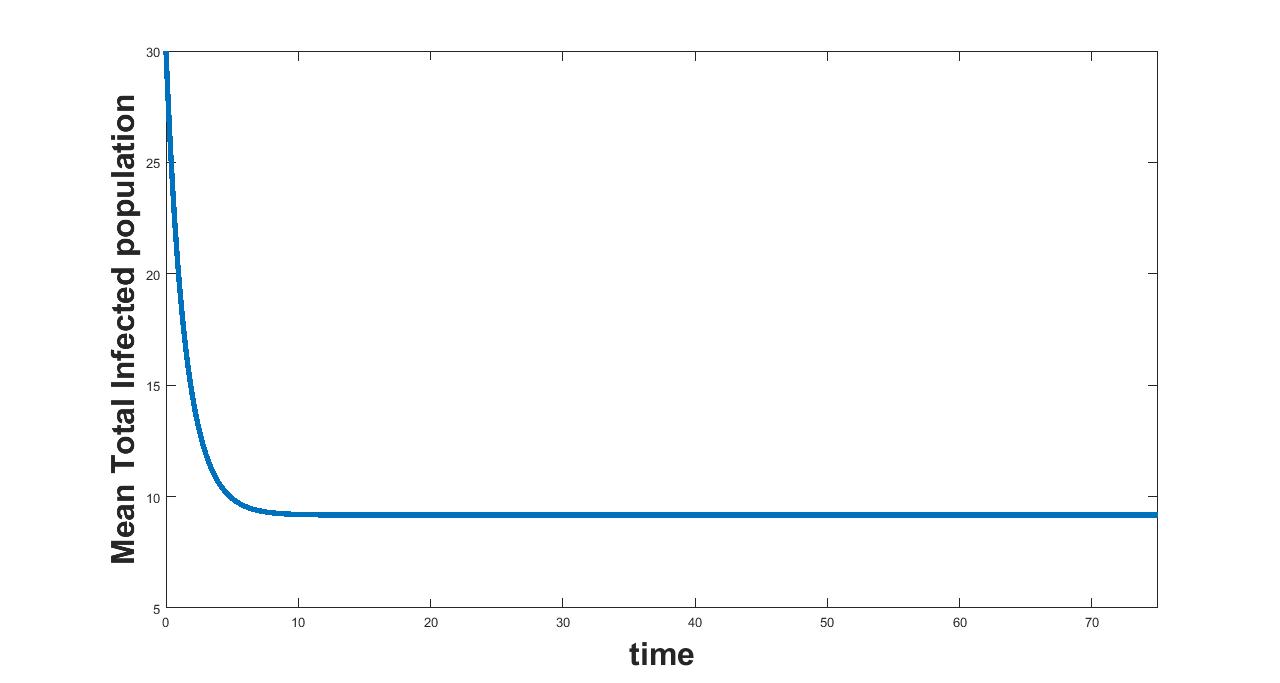}
				\hspace{-.395cm}
					\includegraphics[width=2.2in, height=1.8in, angle=0]{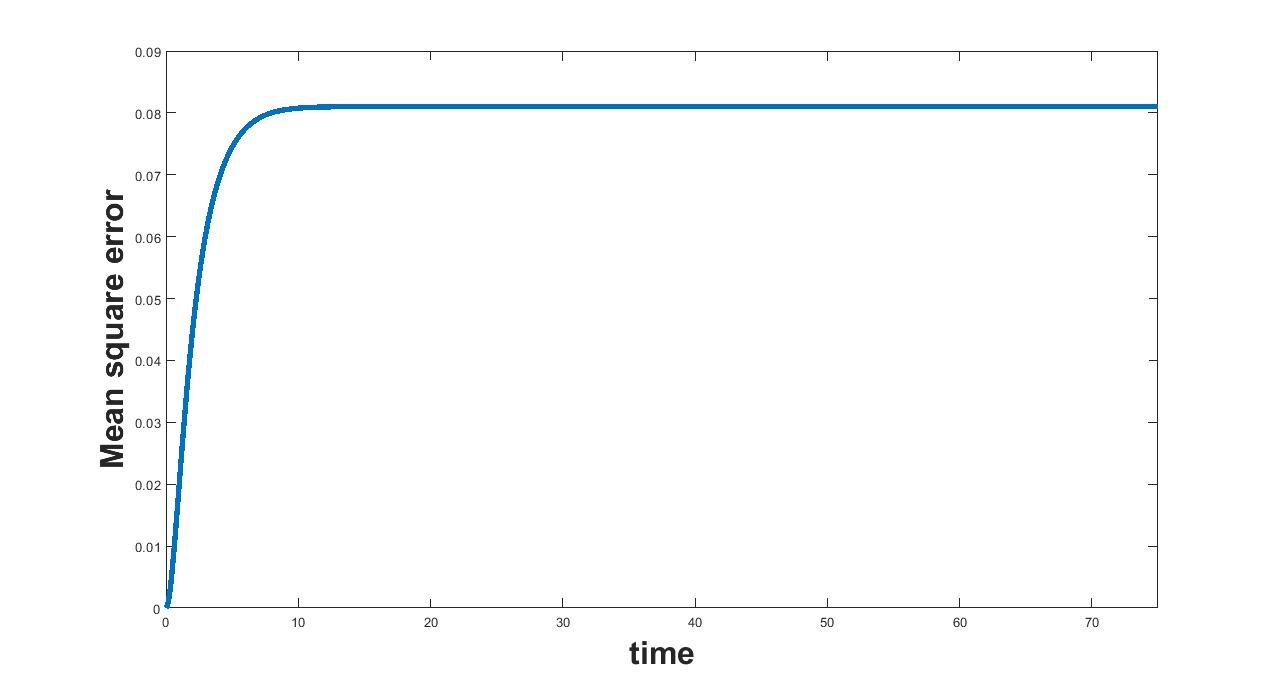}
			\caption*{(a) Interval I :  6.5 to 7.1920}
				
			\end{center}
		\end{figure}
	\begin{figure}[hbt!]
			\begin{center}
				\includegraphics[width=2.2in, height=1.8in, angle=0]{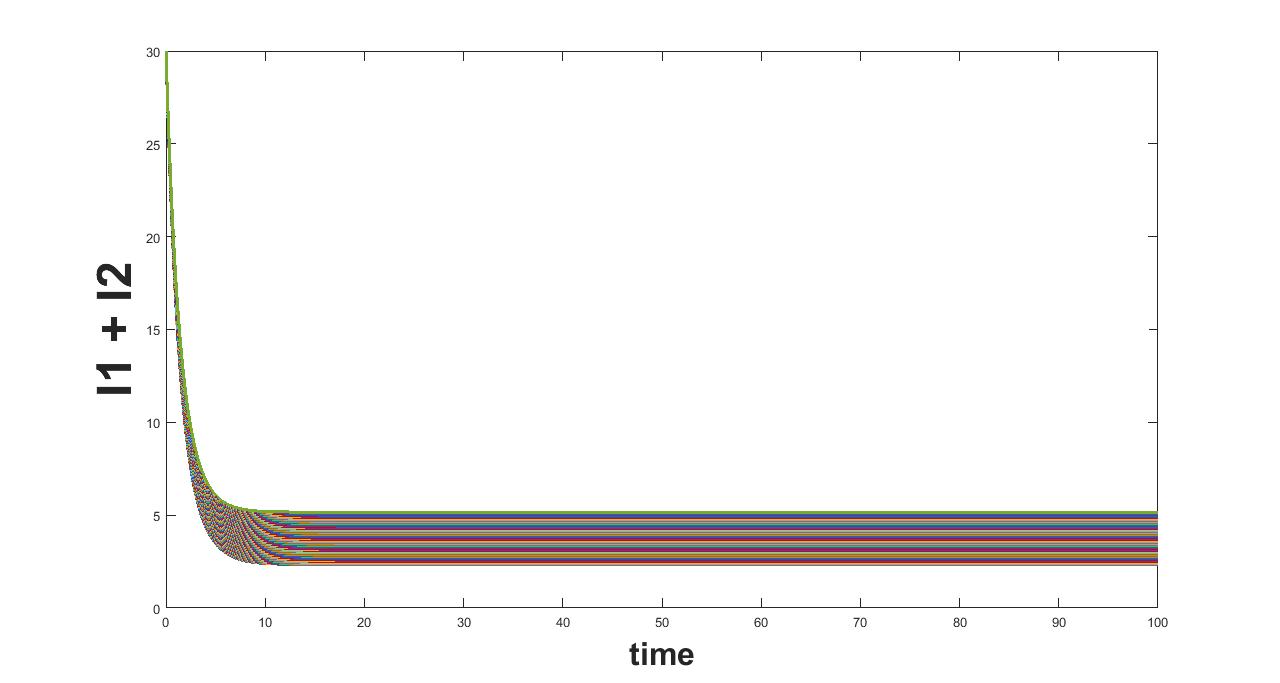}
				\hspace{-.4cm}
				\includegraphics[width=2.2in, height=1.8in, angle=0]{b1mean1}
				\hspace{-.395cm}
					\includegraphics[width=2.2in, height=1.8in, angle=0]{b1error1}
			\caption*{(b) Interval I :  7.1920 to 8}
				
			\end{center}
		\end{figure}

		\begin{figure}[hbt!]
			\begin{center}
				\includegraphics[width=2.2in, height=1.8in, angle=0]{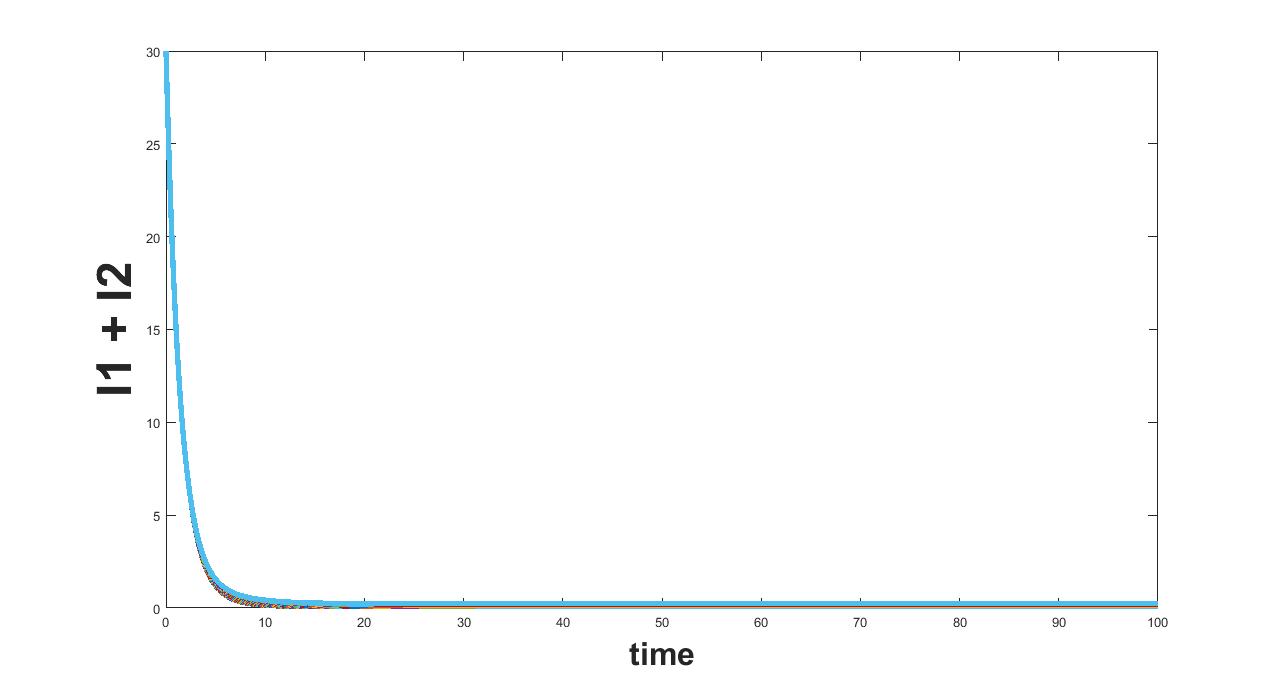}
				\hspace{-.4cm}
				\includegraphics[width=2.2in, height=1.8in, angle=0]{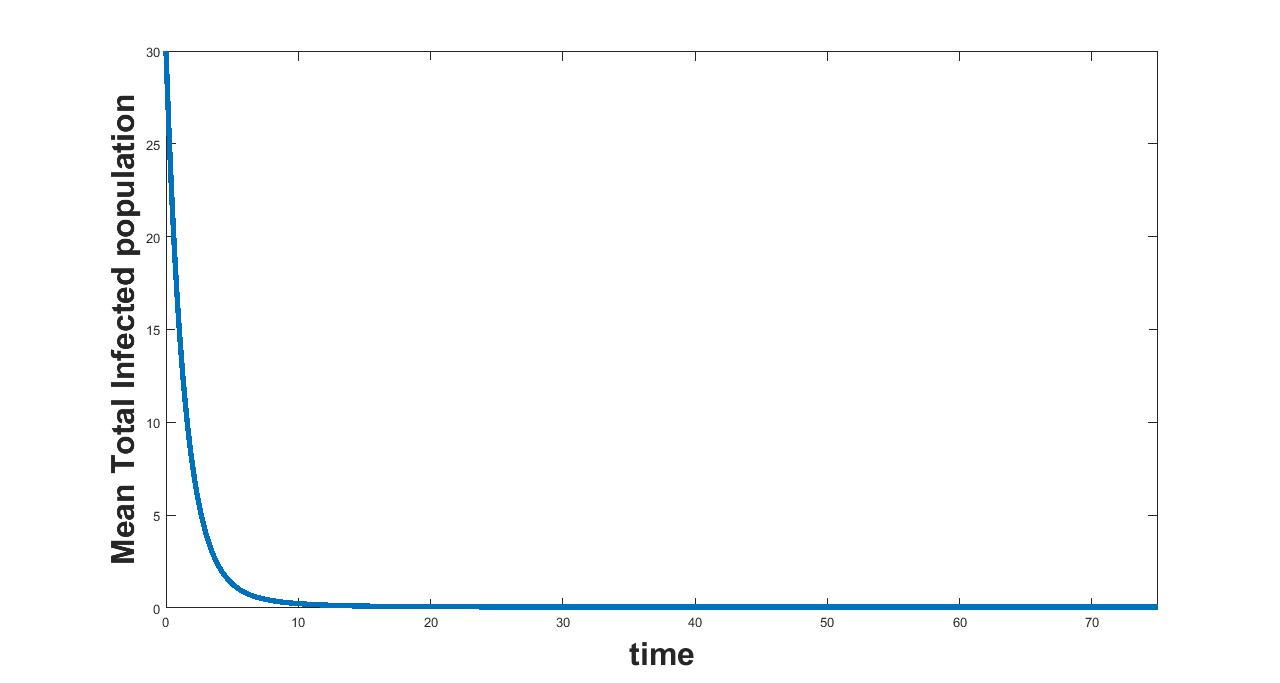}
				\hspace{-.395cm}
					\includegraphics[width=2.2in, height=1.8in, angle=0]{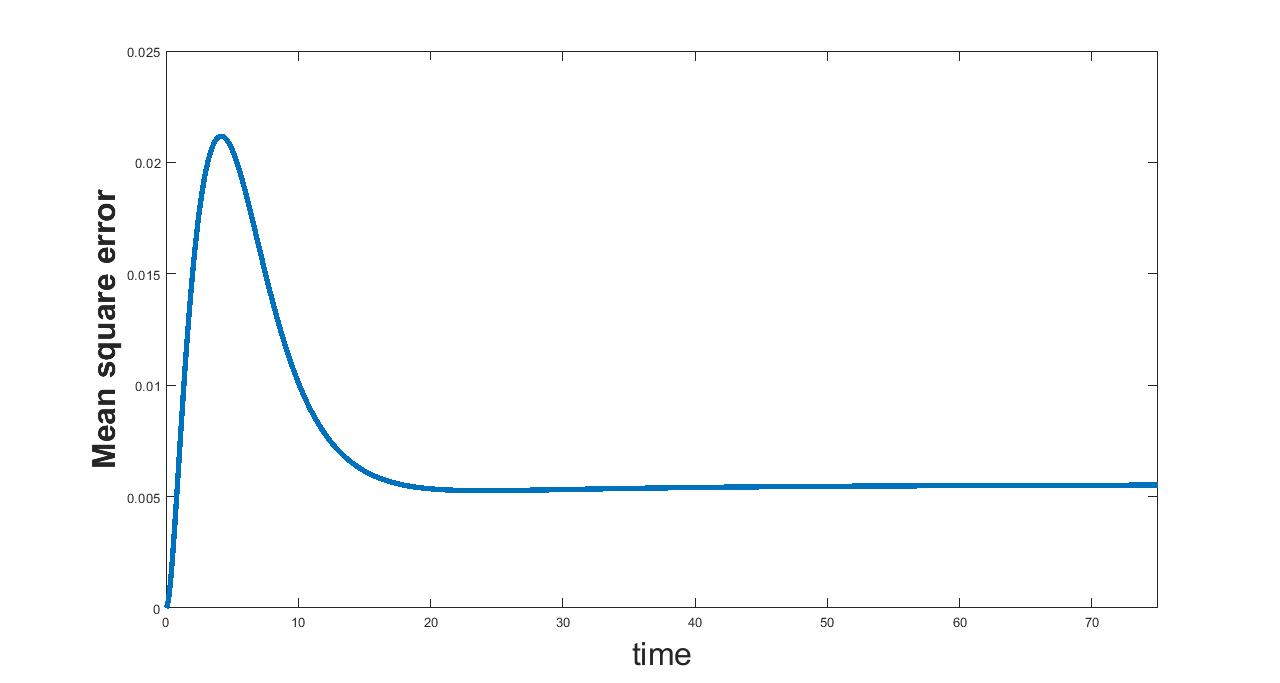}
			\caption*{(c) Interval I :  0.1 to 0.5}
				\label{b1}
			\end{center}
		\end{figure}

		\subsubsection{Parameter $\boldsymbol{\mu}$}

		The results related to sensitivity of $\mu$, varied in two intervals as mentioned in table \ref{t}, are given in figure (7-8). The plots of infected population for each varied value of the parameter $\mu$ per interval, the mean infected population and the mean square error are used to determine the sensitivity. We conclude from these plots that the parameter $\mu$ is sensitive in interval I and insensitive in II.
		\begin{figure}[hbt!]
			\begin{center}
				\includegraphics[width=2.2in, height=1.8in, angle=0]{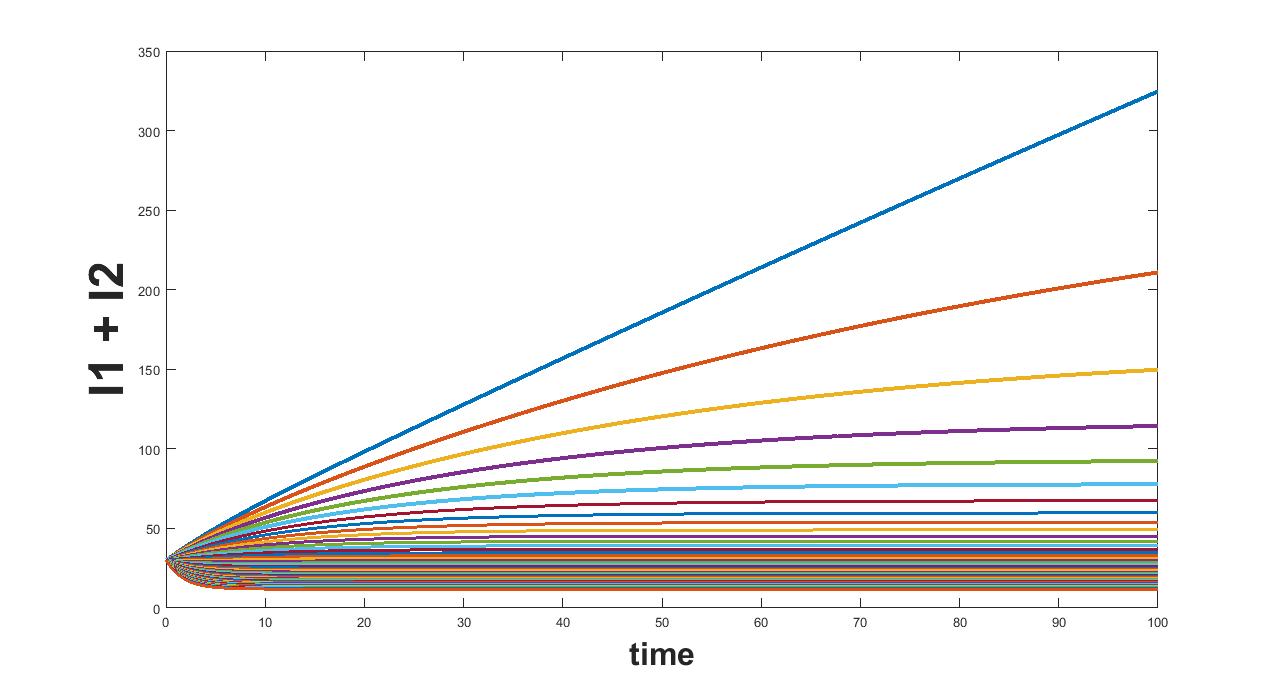}
				\hspace{-.4cm}
				\includegraphics[width=2.2in, height=1.8in, angle=0]{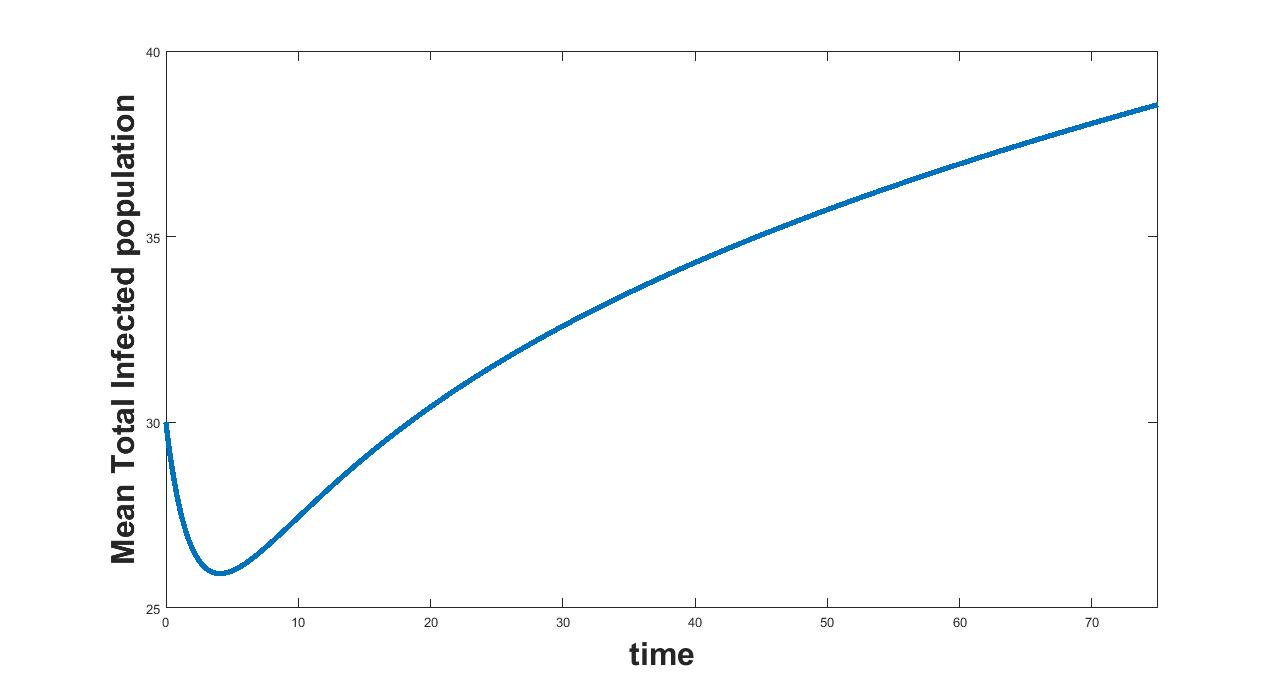}
				\hspace{-.395cm}
					\includegraphics[width=2.2in, height=1.8in, angle=0]{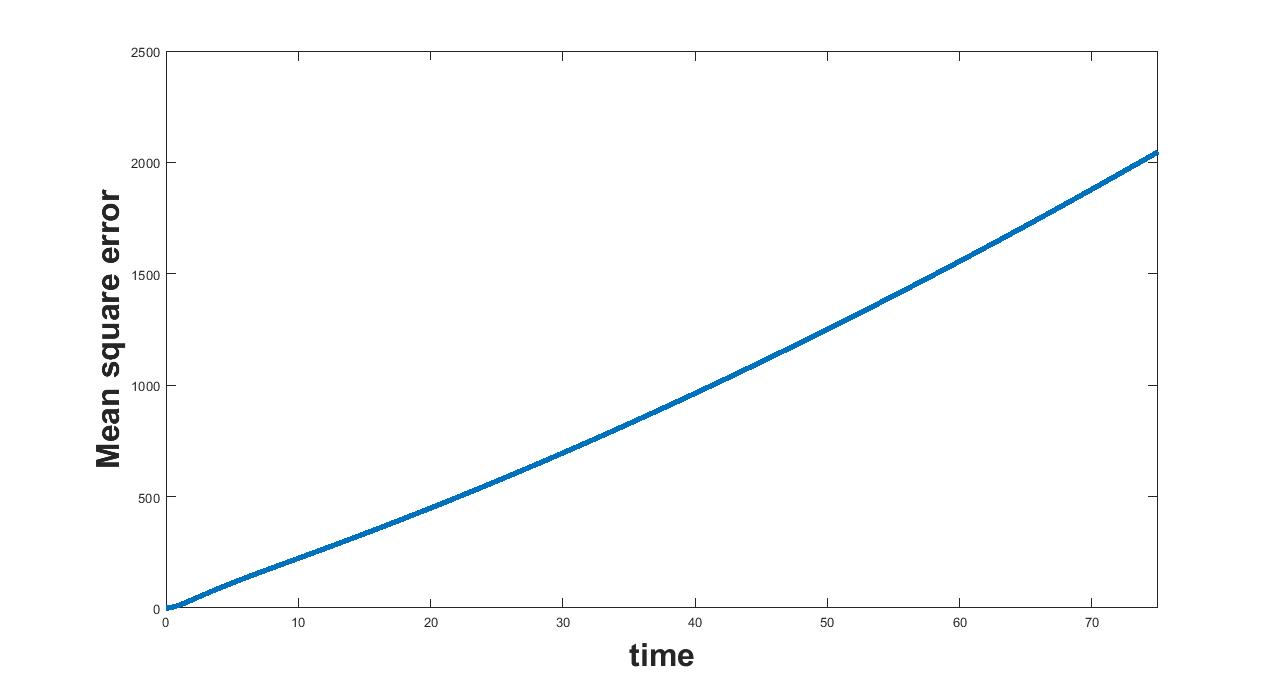}
			\caption*{(a) Interval I :  0 to 0.5}
				
			\end{center}
		\end{figure}
	\begin{figure}[hbt!]
			\begin{center}
				\includegraphics[width=2.2in, height=1.8in, angle=0]{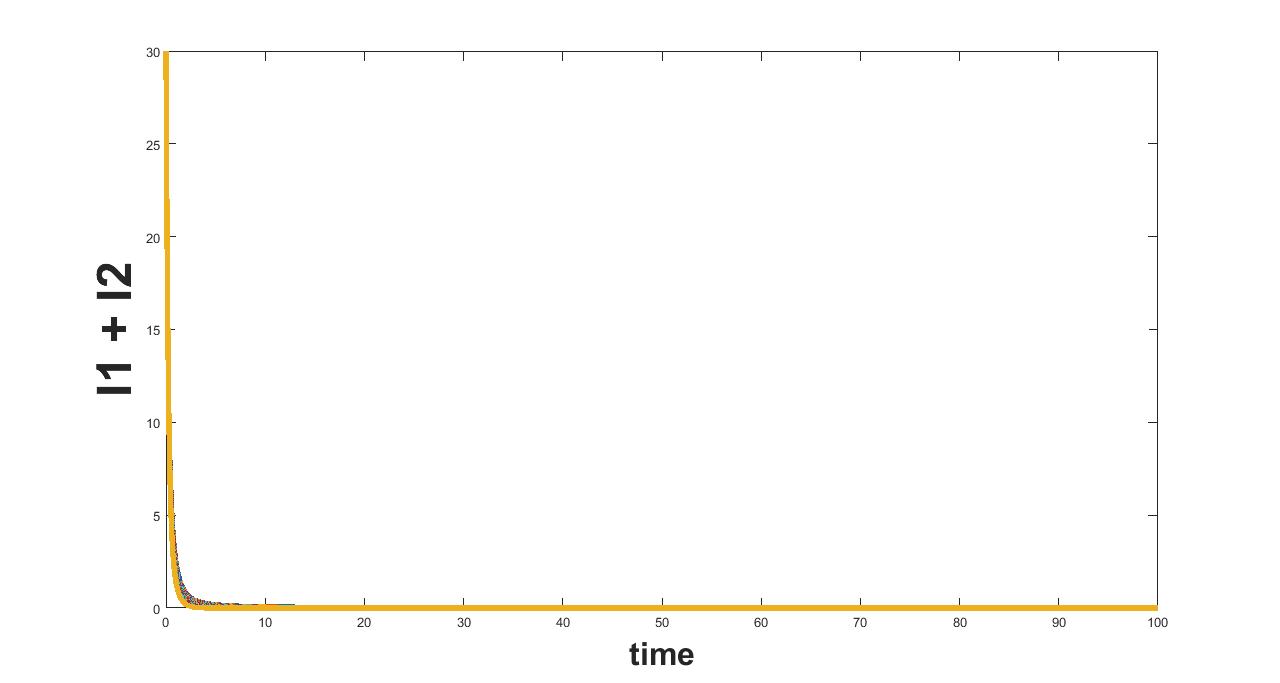}
				\hspace{-.4cm}
				\includegraphics[width=2.2in, height=1.8in, angle=0]{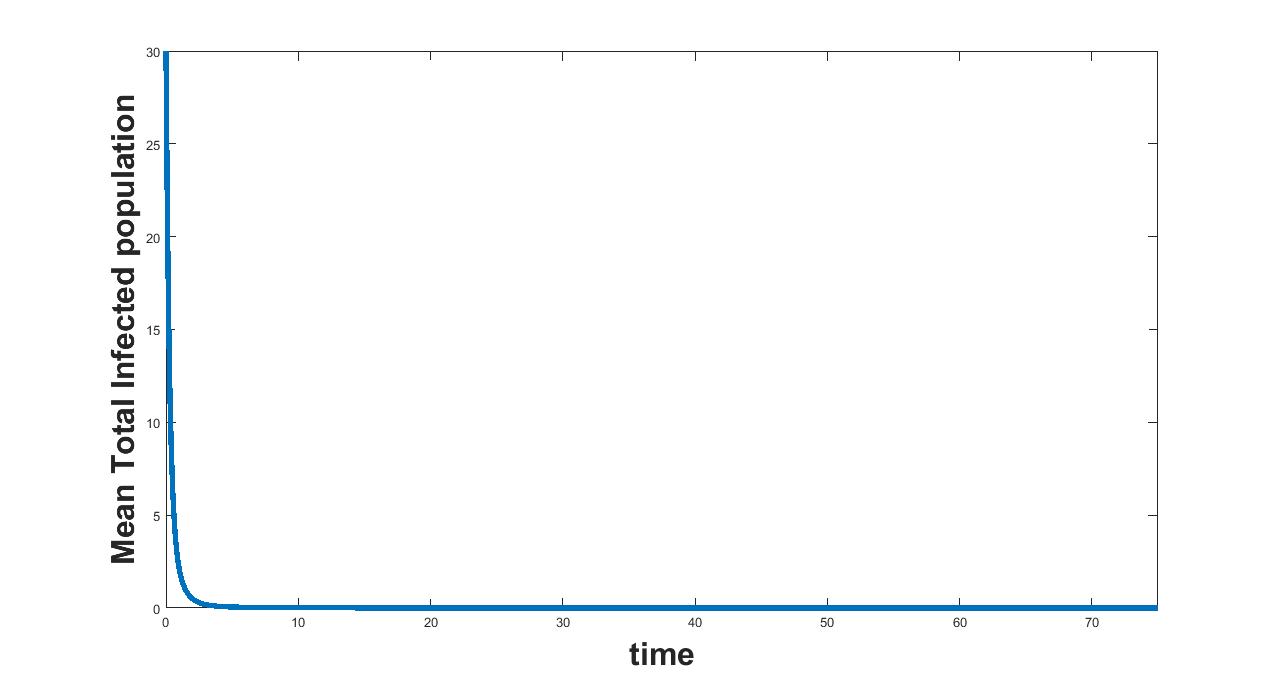}
				\hspace{-.395cm}
					\includegraphics[width=2.2in, height=1.8in, angle=0]{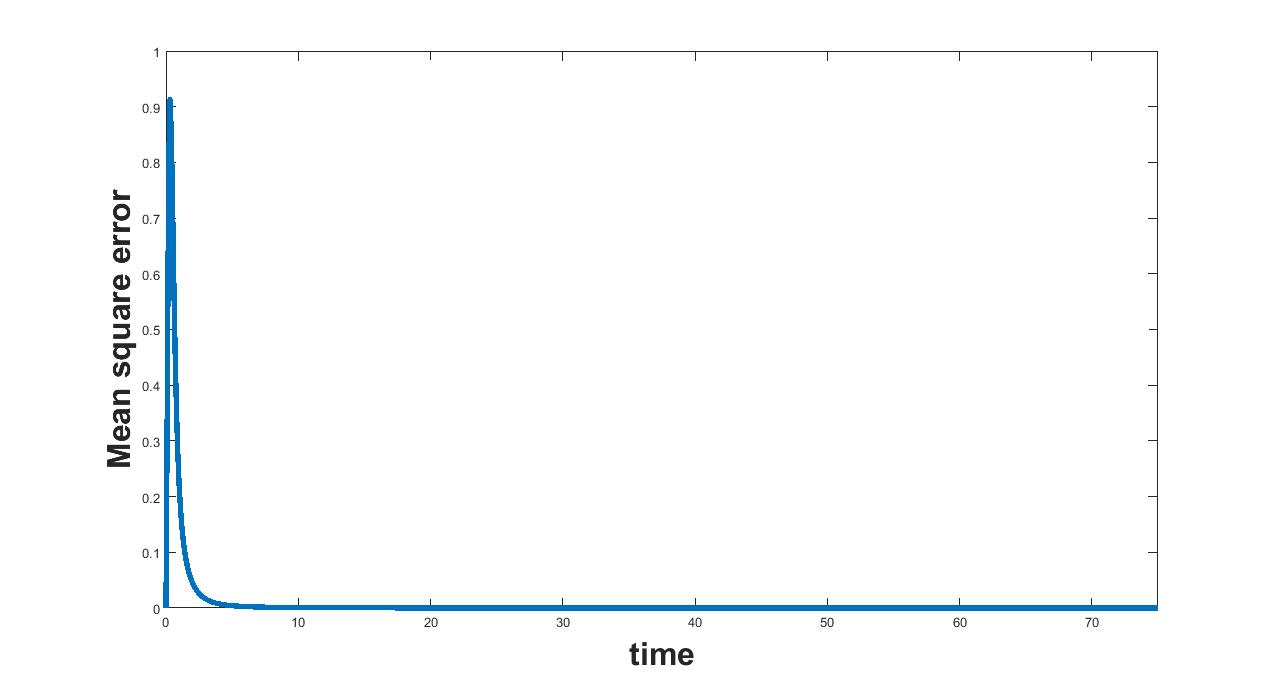}
			\caption*{(b) Interval I :  0.5 to 2}

				\label{mu}
			\end{center}
		\end{figure}

		\subsubsection{Parameter $\boldsymbol{\beta_{1}}$}

		The results related to sensitivity of $\beta_1$, varied in two intervals as mentioned in table \ref{t}, are given in figure (9-10). The plots of infected population for each varied value of the parameter $\beta_1$ per interval, the mean infected population and the mean square error are used to determine the sensitivity. We conclude from these plots that the parameter $\beta_1$ is sensitive in interval I and insensitive in II. In similar lines, the sensitivity analysis is done for other parameters. The results are summarized in table 6. The corresponding plots are given in Appendix - A owing to the brevity of the manuscript.
		\begin{figure}[hbt!]
			\begin{center}
				\includegraphics[width=2.2in, height=1.8in, angle=0]{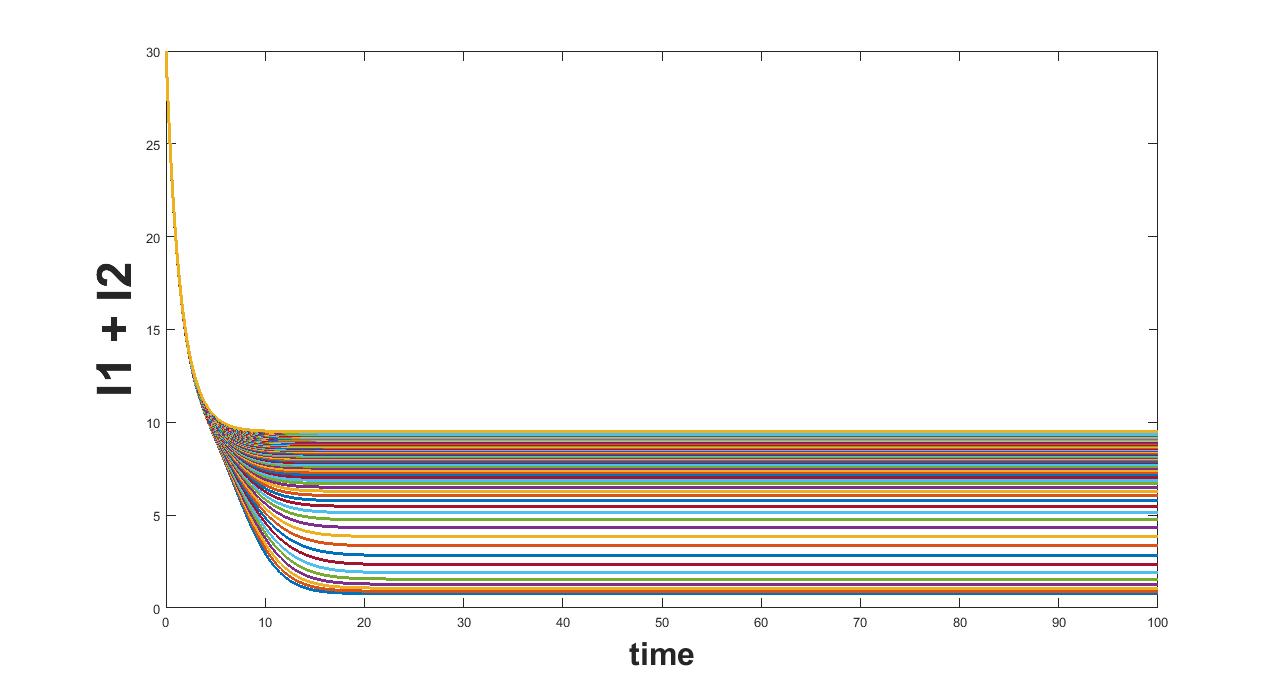}
				\hspace{-.4cm}
				\includegraphics[width=2.2in, height=1.8in, angle=0]{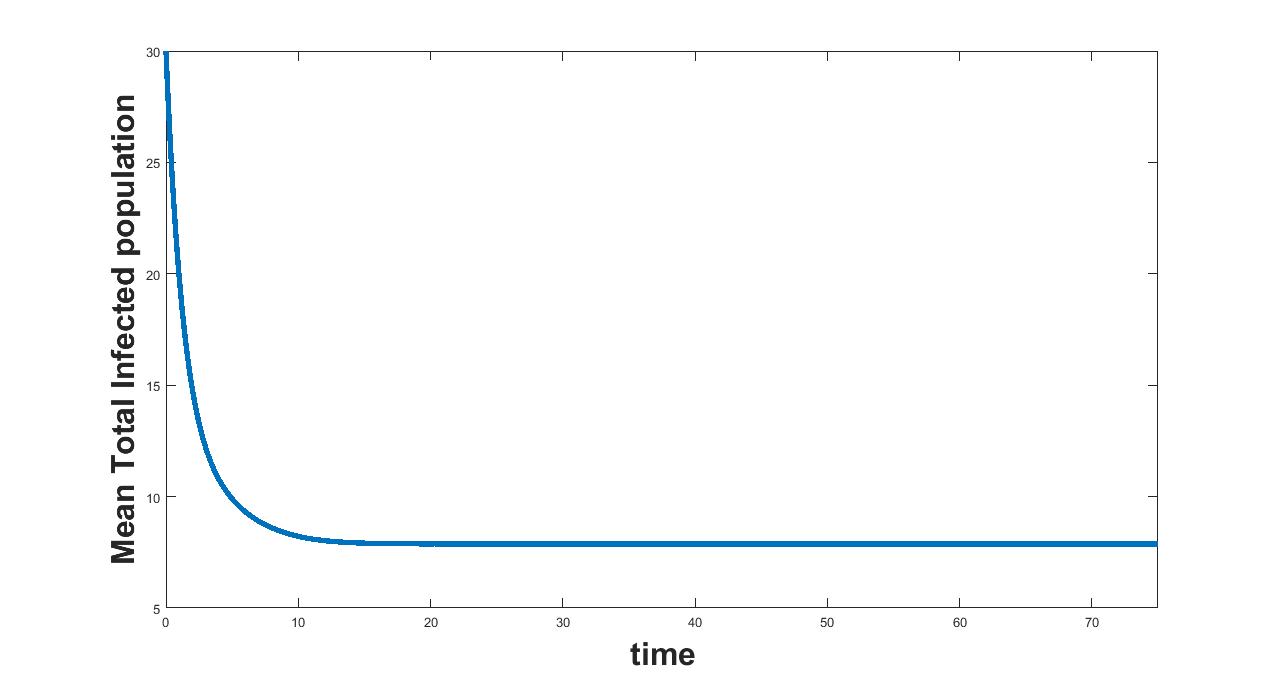}
				\hspace{-.395cm}
					\includegraphics[width=2.2in, height=1.8in, angle=0]{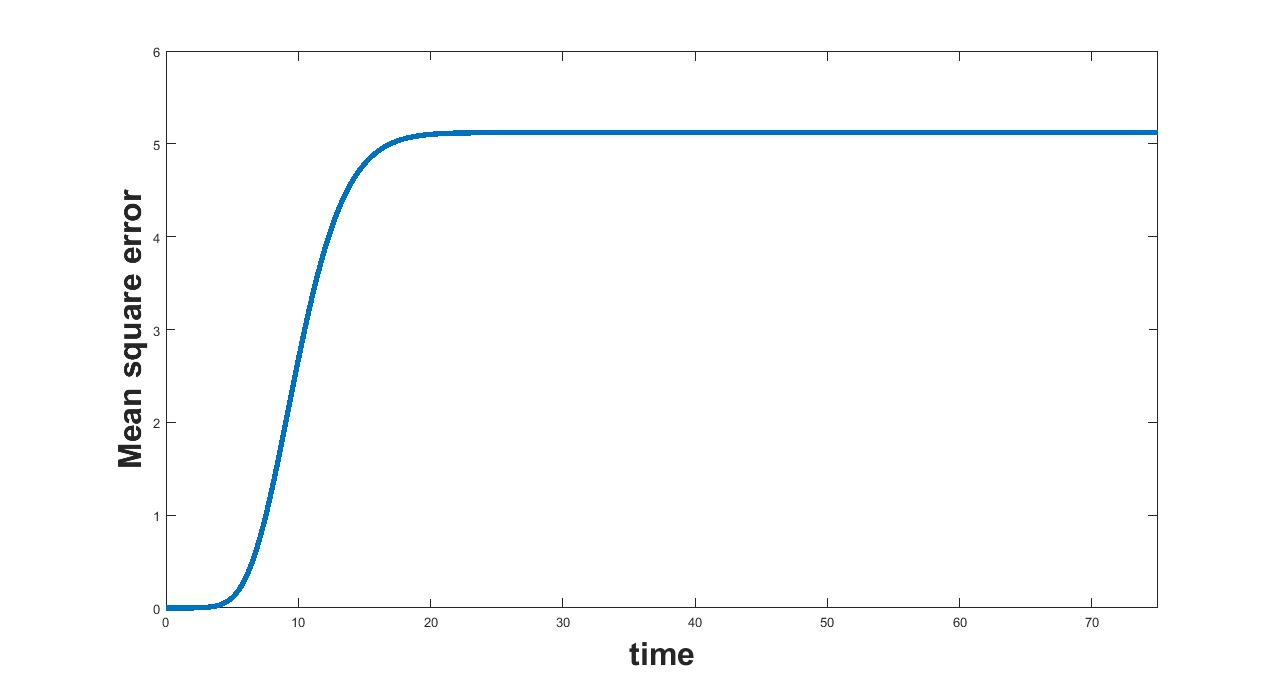}
			\caption*{(a) Interval I :  0 to 1.33}
				
			\end{center}
		\end{figure}
		\newpage
	\begin{figure}[hbt!]
			\begin{center}
				\includegraphics[width=2.2in, height=1.8in, angle=0]{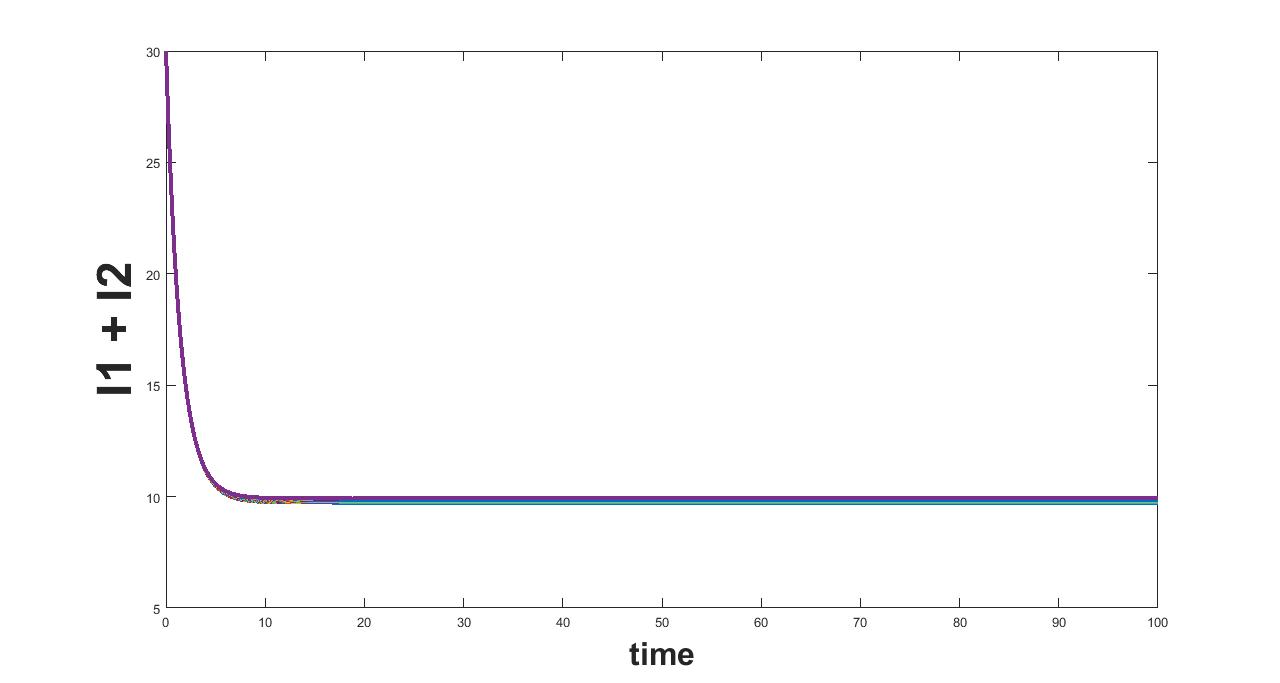}
				\hspace{-.4cm}
				\includegraphics[width=2.2in, height=1.8in, angle=0]{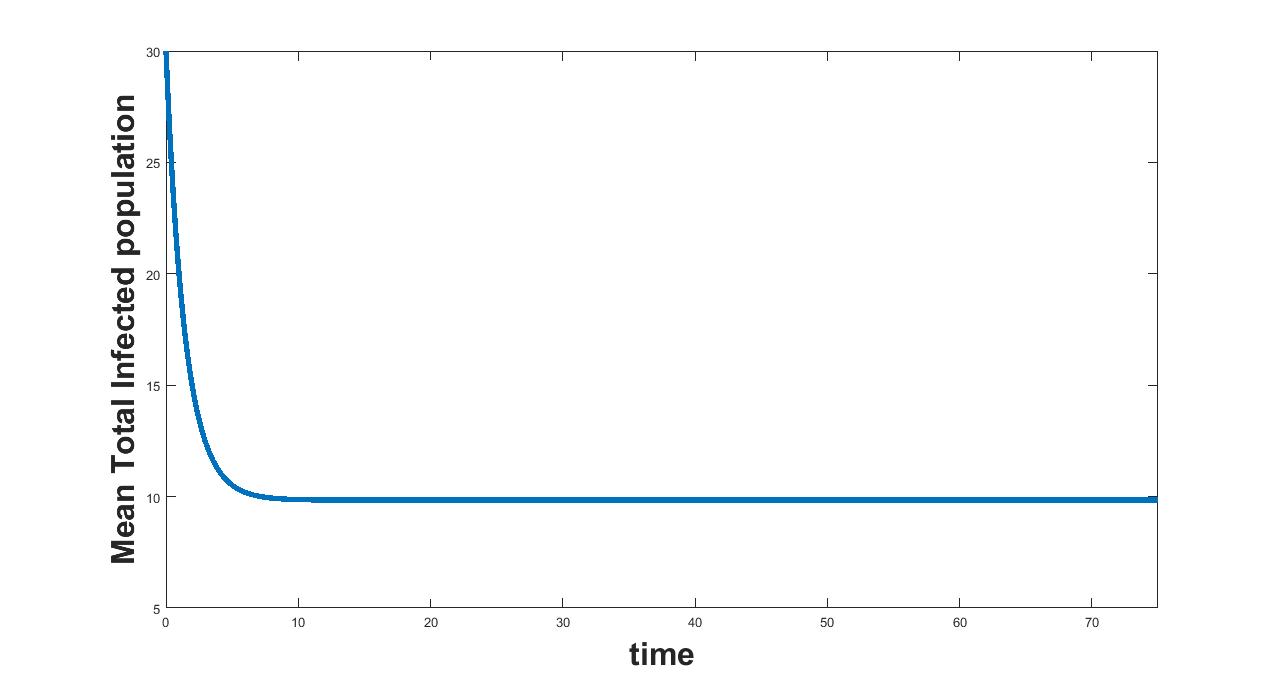}
				\hspace{-.395cm}
					\includegraphics[width=2.2in, height=1.8in, angle=0]{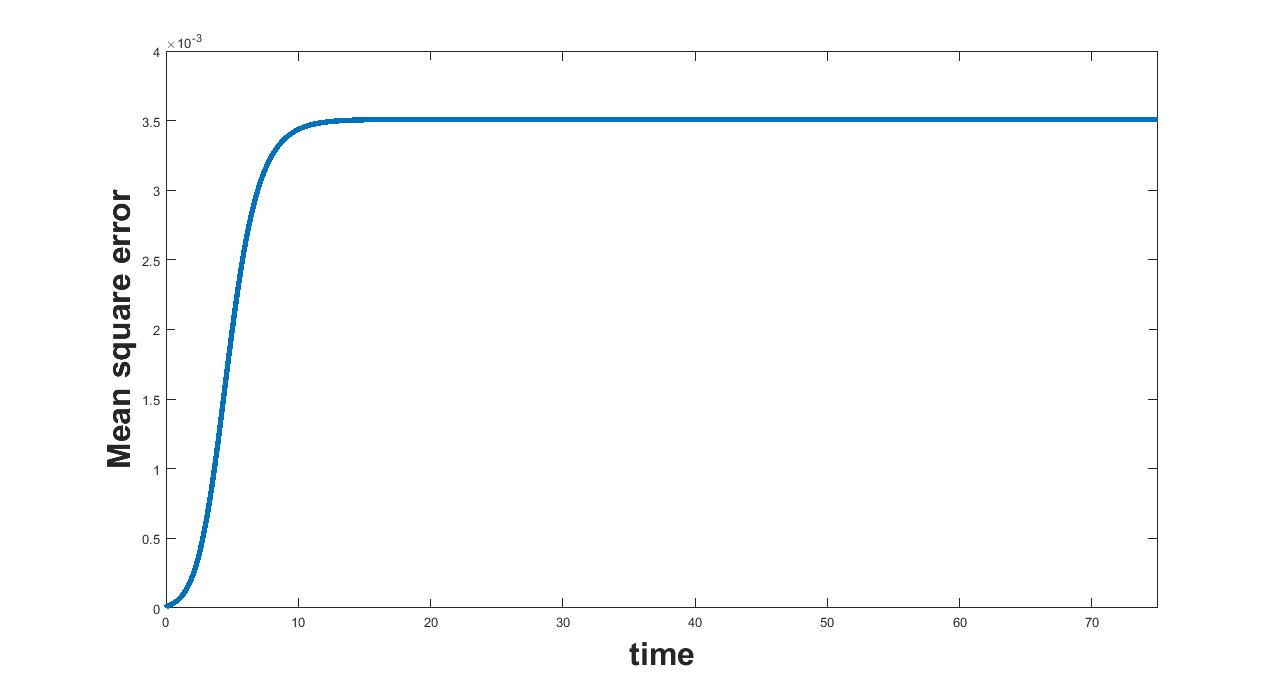}
			\caption*{(b) Interval I :  1.33 to 2}

				\label{beta1}
			\end{center}
		\end{figure}
			
	\begin{table}[hbt!]
		\caption{Summary of Sensitivity Analysis}
		\centering
		\label{sen_anl}
		{
			\begin{tabular}{|l|l|l|}
				\hline
				\textbf{Parameter} & \textbf{Interval} & \textbf{Step Size}  \\
				\hline
				
				$u_{11}$ & 0 to 0.5 & \checkmark
					\\ \cline{2-3}
				 & 1.5 to 2  & $\times$\\
				 \hline
				 $b_1$ & 6.5 to 7.1920 & \checkmark
					\\ \cline{2-3}
				 & 7.1920  to 8  & \checkmark
				 	\\ \cline{2-3}
				 & 0.1 to 0.5  & $\times$\\
				 \hline
				 	$ m $ & 0 to 0.00182 & $\times$
					\\ \cline{2-3}
				 & .00182 to 1  & $\times$ \\
				 \hline
				$u_{12}$ & 0 to .5 & $\times$
				\\ \cline{2-3} 
				& 0.5 to 2 & $\times$
				\\ \hline
				
				$\beta_1$ & 0 to 1.33 & \checkmark 
				\\ \cline{2-3}
				& 1.33 to 2 & $\times$ \\ 
				\hline 
					$\beta_2 $ & 0 to 2 & $\times$
				\\ \cline{2-3}
				& 2 to 3 &  $\times$\\ 
				\hline 
				
					$\beta_3$ & 0 to 2.5 & $\times$
				\\ \cline{2-3}
				& 2.5 to 5 &  $\times$\\ 
			
				\hline 
					$\beta_4$ & 0 to 0.5 & $\times$ 
				\\ \cline{2-3}
				& 0.5 to 1 &  $\times$\\ 
				\hline 
				$\alpha$ & 0 to 0.5 &   $\times$
				\\ \cline{2-3}
				& 0.5 to 2 &  $\times$ \\
				\hline 
				$d_1$ & 0 to 0.000073  & $\times$ 
				\\ \cline{2-3}
				& 0.000073 to 1 & \checkmark\\
				\hline 
				$d_2$ & 0 to 0.0000913 & $\times$ 
				\\ \cline{2-3}
				& 0.0000913 to 2 & $\times$\\
				\hline 
				$\mu$ & 0 to 0.5 & \checkmark
				\\ \cline{2-3}
				& 0.5 to 2 & $\times$\\
			 
				\hline 
			$\delta_1$ & 0 to 0.0714 & $\times$ 
				\\ \cline{2-3}
				& 0.0714 to 1 & $\times$\\
				\hline 
				$\delta_2$ & 0 to 0.0714 & $\times$  
				\\ \cline{2-3}
				& 0.0714 to 1 &$\times$ \\
				\hline 
			\end{tabular}
		}
	\end{table}

\newpage
\section{OPTIMAL CONTROL PROBLEM} \vspace{.25cm}

	In this section, we will formulate an optimal control problem to see the role of treatment in reducing the number of infection. The controls that we considered are:
	
	\begin{enumerate}
		\item  The first control that we consider is the  treatment based recovery of young infected individuals as a result of which an infected young individual recovers. This treatment could be any common drugs interventions that are recommended in the treatment of COVID-19. Many potential vaccines for COVID-19 are being studied and some are under clinical trials. Drugs such as remdesivir, favipiravir, ivermectin, lopinavir/ritonavir, mRNA-1273, phase I trial (NCT04280224) and AVT technology are being used as therapeutic agents by different countries for treating Covid-19  \cite{10,8,9,tu2020review}. We denote this intervention by the control variable  $u_{11}$.

	    \item The second control that we consider is the treatment based recovery of adult Population. We consider a limited treatment as $\frac{\mu_2 I_2^{2}}{1+\alpha I_{2}^2}$ where $\mu_{2}(t)$ is a constant treatment rate. Here $\frac{\mu_2 I_2^{2}}{1+\alpha I_{2}^2}$ represents Holling type III response where the initial treatment is less because of less infectives and as the disease spreads, the treatment also increases and non-linearly saturates due to the limitation in the medical facilities.

	\end{enumerate}
	
	The set of all admissible controls is given by \\
	
	$U = \left\{(u_{11}(t),u_{12}(t)) : u_{11}(t) \in [0,u_{11} max] , u_{12}(t) \in [0,u_{12} max] ,t \in [0,T] \right\}$
	
	Without medical interventions, $u_{11} $ and $ \; u_{12},\;$  are just constant parameters.
	
 In order to reduce the complexity of the problem here we choose to model the control efforts via a linear combination of the quadratic terms. Also when the objective function is quadratic with respect to the control, diﬀerential equations arising from optimization have a known solution. Other functional forms sometimes lead to systems of diﬀerential equations that are diﬃcult to solve (\cite{djidjou2020optimal}, \cite{lee2010optimal}). Based on these we now propose and  define the optimal control problem with the goal to reduce the cost functional defines as follows,
	
	\begin{equation}
		J(u_{11}(t),u_{12}(t)) = \int_{0}^{T} (A_{1}u_{11}(t)^2+A_{2}u_{12}(t)^2+I_1(t)+I_2(t)) dt   
		\label{oobj} 
	\end{equation} 
	
	where $ u=(u_{11}(t),u_{12}(t)) \in U$ \\
	
	subject to the system 
	
\begin{eqnarray}
	\begin{aligned}
   	\frac{dS_{1}}{dt}& =&  b_{1}+ \delta_{1}R_{1} \ - \beta_{1} S_{1}I_{1} - \beta_{2} S_{1}I_{2} -\mu S_{1} -m S_{1}  \\
   	\frac{dI_{1}}{dt} &=&  \beta_{1} S_{1}I_{1} + \beta_{2} S_{1}I_{2} \ - d_{1}I_{1}   \ - \mu I_{1}-\mu_{11}(t)I_{1}  \\ 
   	\frac{dR_{1}}{dt} &=&   \mu_{11}(t)I_{1} \ -  \mu R_{1}-\delta_{1}R_{1} -m R_{1} \label{9}\\
  \frac{dS_{2}}{dt}& =& m S_{1} + \delta_{2}R_{2} \ - \beta_{3} S_{2}I_{1} - \beta_{4} S_{2}I_{2} -\mu S_{2}  \\
   	\frac{dI_{2}}{dt} &=&  \beta_{3} S_{2}I_{1} + \beta_{4} S_{2}I_{2} \ - d_{2}I_{2}   \ - \mu I_{2}-\frac{\mu_{12}(t)I_{2}^2}{1+\alpha I_{2}^2}  \\ 
   	\frac{dR_{2}}{dt} &=&  m R_{1}+\frac{\mu_{12}(t)I_{2}^2}{1+\alpha I_{2}^2}  \ -  \mu R_{2}-\delta_{2}R_{2} 
   	\end{aligned}
   \end{eqnarray} 
	
	The integrand of the cost function (6.1), denoted by 
	$$L(S,I,V,u_{11},u_{12}) = (I_1(t)+I_2(t)+A_{1}u_{11}(t)^2+A_{2}u_{12}(t)^2)$$
	
	is called the Lagrangian or the running cost.
	
	Here, the cost function represents the number of infected population throughout the observation period, and the overall cost of implementation of the treatments. Effectively, we want to minimize the infected population and the  cost. Here, $A_1$ and $A_2$ are positive weight constants which not only balance units of integrand but also related cost.
	
	The admissible solution set for the Optimal Control Problem (6.1)-(6.2) is given by
	
	$\Omega = \left\{ (S_1, I_1, R_1, S_2, I_2, R_2,  u_{11}, u_{12})\; | \; S_1,  I_1 , S_2, I_2, R_1,  R_2  \text{that satisfy }(8)\right\}$

	
	
	{\textbf{EXISTENCE OF OPTIMAL CONTROL}}\vspace{.25cm}

	We will show the existence of optimal control functions that minimize the cost functions within a finite time span $[0,T]$ showing that we satisfy the conditions stated in Theorem 4.1 of \cite{Wendell}.
	
	\begin{thm}
		There exists a 2-tuple of optimal controls $(u_{11}^{*}(t) , u_{12}^{*}(t))$ in the set of admissible controls U such that the cost functional is minimized i.e., 
		
		$$J[u_{11}^{*}(t) , u_{12}^{*}(t)] = \min_{(u_{11}^{*} , u_{12}^{*} ) \in U} \bigg \{ J[u_{11},u_{12}]\bigg\}$$ 
		corresponding to the optimal control problem (6.1)-(6.2).
	\end{thm}
	
	
	
	\begin{proof}
		
		 In order to show the existence of optimal control functions, we will show that the following conditions are satisfied : 
		
		\begin{enumerate}
			\item  The solution set for the system (6.2) along with bounded controls must be non-empty, $i.e.$, $\Omega \neq \phi$.
			
			\item  U is closed and convex and system should be expressed linearly in terms of the control variables with coefficients that are functions of time and state variables.
			
			\item The Lagrangian L should be convex on U and $L(S,I,V,u_{11},u_{12}) \geq g(u_{11},u_{12})$, where $g(u_{11},u_{12})$ is a continuous function of control variables such that $|(u_{11},u_{12})|^{-1} g(u_{11},u_{12}) \to \infty$ whenever  $|(u_{11},u_{12})| \to \infty$, where $|.|$ is an $l^2(0,T)$ norm.
		\end{enumerate}

		Now we will show that each of the conditions are satisfied : 
		
		1. From Positivity and boundedness of solutions of the system(6.2), all solutions are bounded for each bounded control variable in $U$.
		
		Also,the right hand side of the system (6.2) satisfies Lipschitz condition with respect to state variables. 
		
		Hence, using the positivity and boundedness condition and the existence of solution from Picard-Lindelof Theorem\cite{makarov2013picard}, we have satisfied condition 1.
		
		2. $U$ is closed and convex by definition. Also, the system (6.2) is clearly linear with respect to controls such that coefficients are only state variables or functions dependent on time. Hence condition 2 is satisfied.
		
		3. Choosing $g(u_{11},u_{12}) = c(u_{11}^{2}+u_{12}^{2})$ such that $c = min\left\{A_{1},A_{2}\right\}$, we can satisfy the condition 3.
		
		Hence there exists a control 2-tuple $(u_{11}^{*},u_{12}^{*})\in U$ that minimizes the cost function (6.1).
	\end{proof}

	\textbf{CHARACTERIZATION OF OPTIMAL CONTROL}\vspace{.25cm}
	
	We will obtain the necessary conditions for optimal control functions using the Pontryagin's Maximum Principle \cite{liberzon2011calculus} and also obtain the characteristics of the optimal controls.
	
	The Hamiltonian for this problem is given by 
	
	$$H(S,I,V,u_{11},u_{12},u_{2},\lambda) := L(S_1,I_1,R_1,S_2,I_2,R_2, u_{11},u_{12}) + \lambda_{1} \frac{\mathrm{d} S_1}{\mathrm{d} t} +\lambda _{2}\frac{\mathrm{d} S_2}{\mathrm{d} t}+ \lambda _{3} \frac{\mathrm{d} I_1}{\mathrm{d} t} + \lambda _{4}\frac{\mathrm{d} I_2}{\mathrm{d} t}+ \lambda _{5} \frac{\mathrm{d} R_1}{\mathrm{d} t}+\lambda_6\frac{\mathrm{d} R_2}{\mathrm{d} t}$$
	
	Here $\lambda$ = ($\lambda_{1}$,$\lambda_{2}$,$\lambda_{3}$,$\lambda_{4}$,$\lambda_{5}$,$\lambda_{6}$) is called co-state vector or adjoint vector.
	
	Now the Canonical equations that relate the state variables to the co-state variables are  given by 
	
	\begin{equation}
	\begin{aligned}
	 \frac{\mathrm{d} \lambda _{1}}{\mathrm{d} t} &= -\frac{\partial H}{\partial S_1}\\
	 \frac{\mathrm{d} \lambda _{2}}{\mathrm{d} t} &= -\frac{\partial H}{\partial S_2}\\
	 \frac{\mathrm{d} \lambda _{3}}{\mathrm{d} t} &= -\frac{\partial H}{\partial I_1}\\
	 \frac{\mathrm{d} \lambda _{4}}{\mathrm{d} t} &= -\frac{\partial H}{\partial I_2}\\
	 \frac{\mathrm{d} \lambda _{5}}{\mathrm{d} t} &= -\frac{\partial H}{\partial R_1}\\
	 \frac{\mathrm{d} \lambda _{6}}{\mathrm{d} t} &= -\frac{\partial H}{\partial R_2}
	\end{aligned}
	\end{equation}

	
	Substituting the Hamiltonian value gives the canonical system 
	
	\begin{equation}
	\begin{aligned}
	\frac{\mathrm{d} \lambda _{1}}{\mathrm{d} t} &= \lambda _{1}(\beta_1 I_1+\beta_2 I_2 + m+\mu)+\lambda _{2} m + \lambda_3 (\beta_1 I_1+\beta_2 I_2)\\
	\frac{\mathrm{d} \lambda _{2}}{\mathrm{d} t} &= \lambda _{2}(\beta_3 I_1+\beta_4 I_2+\mu)-\lambda _{4} (\beta_3 I_1+\beta_4 I_2)\\
	\frac{\mathrm{d} \lambda _{3}}{\mathrm{d} t} &= -1 + \lambda_1 \beta_1 S_1 + \lambda_2 \beta_3 S_2 - \lambda_3(\beta_1 S_1-d_1-\mu-\mu_{11}(t))-\lambda _{4}\beta_3 S_2-\lambda _{5} \mu_{11}(t)\\
	\frac{\mathrm{d} \lambda _{4}}{\mathrm{d} t} &= \lambda _{1}\beta_2 S_1+\lambda _{2} \beta_4 S_2 - \lambda_3 \beta_2 S_1 -\lambda_4 \bigg(\beta_4 S_2 - d_2 -\mu -\frac{2\mu_{12}(t) I_2}{1+ \alpha I_{2}^2}\bigg)-\lambda_6 \bigg(\frac{2\mu_{12}(t) I_2}{1+ \alpha I_{2}^2}\bigg)\\
	\frac{\mathrm{d} \lambda _{5}}{\mathrm{d} t} &= -\lambda _{1} \delta_1+\lambda _{5} (\mu + \delta_1 + m) - \lambda_6 m\\
	\frac{\mathrm{d} \lambda _{6}}{\mathrm{d} t} &= -\lambda_2 \delta_2 + \lambda_6 (\mu + \delta_2)
	\end{aligned}
	\end{equation}
	
	along with transversality conditions
	$ \lambda _{1} (T) = 0, \  \lambda _{2} (T) = 0, \  \lambda _{3} (T) = 0.\lambda _{4} (T) = 0, \  \lambda _{5} (T) = 0, \  \lambda _{6} (T) = 0. $
	
	Now, to obtain the optimal controls, we will use the Hamiltonian minimization condition 
	$ \frac{\partial H}{\partial u_{i}}$ = 0 , at  $u_{i} = u_{i}^{*}$  for i = 11, 12 .
	
	Differentiating the Hamiltonian and solving the equations, we obtain the optimal controls as 
	
	\begin{eqnarray*}
	u_{11}^{*} &=& \min\bigg\{ \max\bigg\{\frac{(\lambda _{3}-\lambda_5)I_1}{2A_{1}},0 \bigg\}, u_{11}max\bigg\}\\
	u_{12}^{*} &= &\min\bigg\{ \max\bigg\{\frac{(\lambda _{4}-\lambda_6)I_2^2}{2A_{2}(1+\alpha I_{2}^2)},0 \bigg\}, u_{12}max\bigg\}
		\end{eqnarray*}
	 
     \section{Simulations for Optimal Control Problem }
     
     	In this section, we perform numerical simulations to understand the role of the treatments in reducing the infection. Entire simulations is done using MATLAB software.
	
The various combinations of controls considered are: 
	
	1. Implementation of  treatment only for young population
	
	2. Implementation of  treatment only for adult population. \vspace{.25cm}
	
	3. Implementation of  treatment for both young and the adult population. \vspace{.25cm}
	
	For our simulations, we have taken the total number of days as $T = 100$ and the other fixed parameters from table 2.
	
	We first solve the state system numerically using Fourth Order Runge-Kutta method in MATLAB without any interventions along with the initial values as $(S_1,S_2,I_1,I_2,R_1,R_2)=(100,100,10,10,5,5)$
	
	Now, to simulate the system with controls, we use the Forward-Backward Sweep method stating with the initial values of controls and solve the state system forward in time. Following this we solve the adjoint state system backward in time due to the transversality conditions, using the optimal state variables and initial values of optimal control.
	
	Now, using the values of adjoint state variables, the values of optimal control are updated and with these updated control variables, we go through this process again. We continue this till the convergence criterion is met \cite{liberzon2011calculus}. In similar lines to \cite{zamir2020non} the positive weights chosen for objective coefficients are $A_{1}$ = .0001, $A_{2}$ = .005. We have chosen the weights related to the treatment of adults  10 times more than that of young ones because the cost related to the treatment of adult population is generally higher compared to the cost of treatment of young. 
	
	First, we solve system (6.2) in the absence of the controls $(u_{11} = 0, u_{12} = 0)$
along with same initial population size. The corresponding count of the infectives $(I_1,I_2)$ are shown in figures 11 and 12 with blue color. From figure 11 we see that the number of infected adult cases considering the implementation of controls are less than the case considering no controls and best result is obtained when the control $u_{12}$ is considered alone as the number of infected adult population is least compared to the other two cases.

\begin{center}
	\begin{figure}[hbt!]
		\includegraphics[height = 10cm, width = 17.5cm]{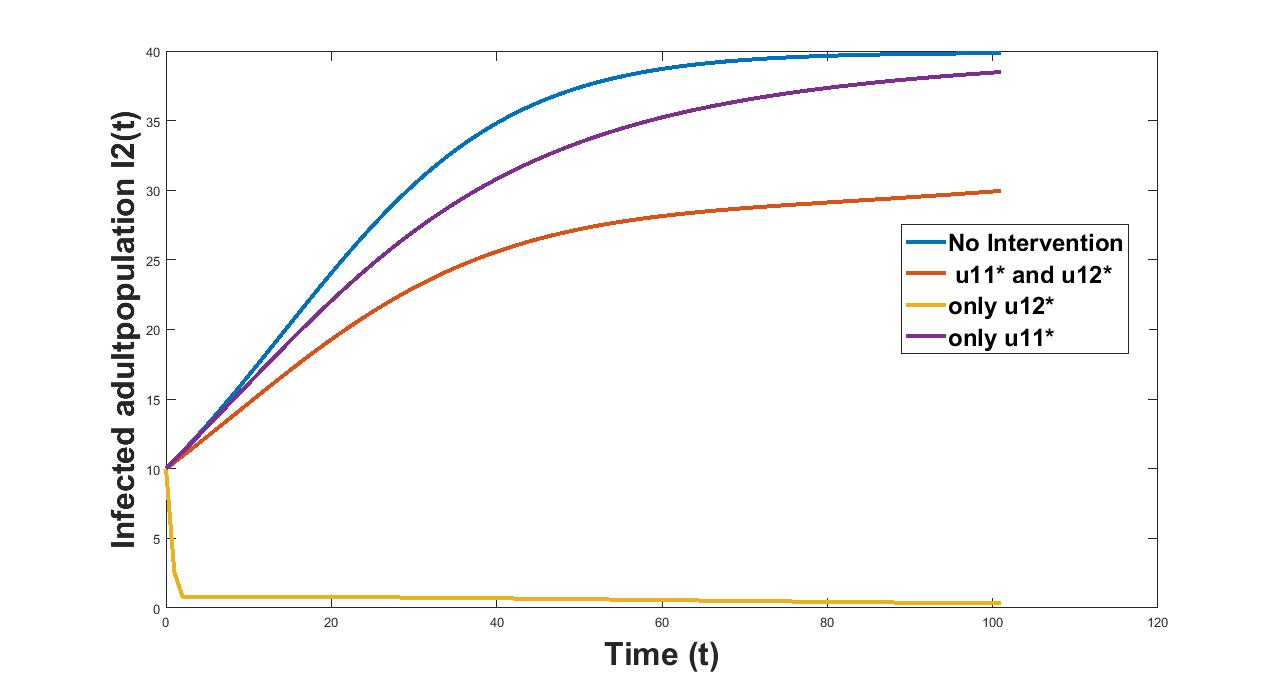} 
		\caption{$I2$ under optimal controls u$_{11}^{*}$, u$_{12}^{*}$}
	
	\end{figure} 
\end{center}

\newpage
We also calculate the average values of infected adult population in table 7 over the time period $t=100$ to support the above fact.

\begin{table}[ht!]
\caption{Table depicting the average values of the  infected adult  population. }
    \centering 
    \begin{tabular}{|l|l|l|l|} 
    \hline
    \textbf{Control Combinations} &  \textbf{Avg Infected Population(I2)} \\
    \hline
    
    $u_{11}^*=0, u_{12}^{*}$ &  0.7185 \\
     \hline
    $u_{11}^*, u_{12}^{*}$ &  24.5626 \\
     \hline
    $u_{11}^*, u_{12}^{*}=0$ &    30.1420 \\
     \hline
    
     $u_{11}^*=u_{12}^{*} = 0 $&    32.6219 \\
    \hline
   
    \end{tabular} 
     \label{t3}
\end{table} \vspace{.4cm}	

In figure 12 we simulate the infected young population over the time considering different control combinations. We see from the figure that the number of infectives is higher when no control is implemented and as treatment is implemented the number of infectives reduces and reduces the maximum when both $u_{11}$ and $u_{12}$ is considered followed by considering only $u_{11}$. 

We also calculate the average values of the infected young population over the time period considered. From table 8 we see that the average value of infected young population is least in case when both the controls are implemented.

\newpage
\begin{center}
	\begin{figure}[hbt!]
		\includegraphics[height = 10cm, width = 17.5cm]{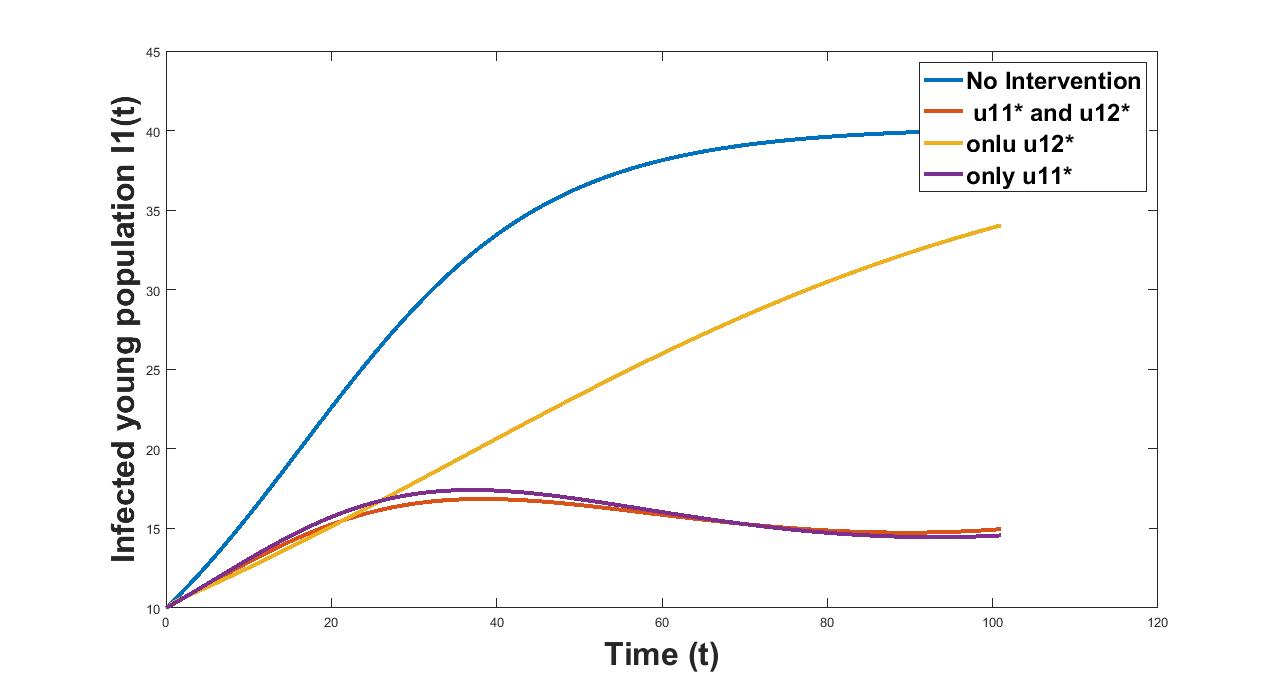} 
	\caption{$I1$ under optimal controls u$_{11}^{*}$, u$_{12}^{*}$}
	\end{figure} 
\end{center}

\begin{table}[ht!]
 \caption{Table depicting the average values of the  infected young population. }
    \centering 
    \begin{tabular}{|l|l|l|l|} 
    \hline
    \textbf{Control Combinations} &  \textbf{Avg Infected Population(I1)} \\
    \hline
    $u_{11}^*, u_{12}^{*} $&     15.1129 \\
     
    \hline
    $u_{11}^*, u_{12}^{*}=0$ &    15.2773 \\
     \hline
    \hline
    $u_{11}^*=0, u_{12}^{*} $ &    22.9808 \\

    \hline
     $u_{11}^*=u_{12}^{*} = 0$ &       31.9349 \\
    \hline
    
    \end{tabular} 
    \label{t3}
\end{table} \vspace{.4cm}

We also simulate the recovered population over the time with and without controls in figure 13 and figure 14. From Figure 13 we see that the recovered young population is highest when the control $u_{11}$ is considered. Whereas the recovered young population remains almost constant at constant value 5 in case of second control $u_{12}$  consideration. From figure  13 we also see that the recovered population  curve considering $u_{12}$ alone and the curve without controls coincide with each other. The reason for this could be that the second control $u_{12}$ features only in the second age group(adult) population considered in the model.

From figure 14 we see that the recovered adult population is highest considering the control $u_{12}$ whereas it is least in case of no control consideration and implementation of only $u_{11}$. Similar to the previous case we see that the recovered population  curve considering $u_{11}$ alone and the curve without controls coincide with each other and the reason for this could be that the control $u_{11}$ features only in the first age group(young) population considered in the model. In figure 15 we plot the cumulative infected population over the time considering different controls and we see that  the cumulative count of the infectives is minimum  when the second control $u_{12}$ is considered followed by considering both the controls. Average values of the recovered population are calculated for young and adult population considering different  combinations  of the controls in table 9 and  table 10.

\begin{center}
	\begin{figure}[hbt!]
		\includegraphics[height = 10cm, width = 17.5cm]{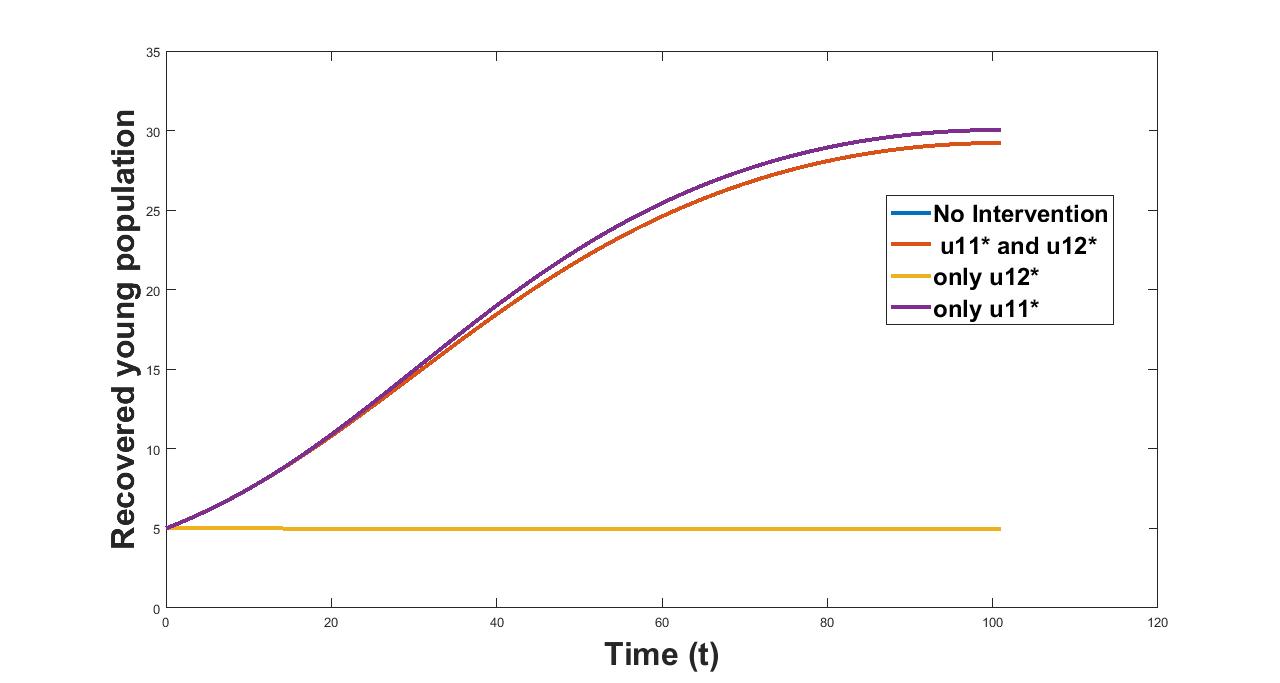} 
	\caption{$R1$ under optimal controls u$_{11}^{*}$, u$_{12}^{*}$}
	\end{figure} 
\end{center}

\vspace{1cm}
\begin{table}[ht!]
\caption{Table depicting the average values of the  Recovered young population. }
    \centering 
    \begin{tabular}{|l|l|l|l|} 
    \hline
    \textbf{Control Combinations} &  \textbf{Avg Recovered Population(R1)} \\
   
     \hline
    
     $u_{11}^*=0, u_{12}^{*}$ &         34.9995 \\
     \hline
   
    $u_{11}^*, u_{12}^{*}$ &           10.4040 \\
     \hline
    
      $u_{11}^*, u_{12}^{*}=0$ &       4.9834 \\
       \hline
    
     $u_{11}^*=u_{12}^{*} = 0 $&         4.9834 \\
    
    \hline
    \end{tabular} 
     \label{t3}
\end{table} 

\vspace{1.5cm}
\begin{table}[ht!]
 \caption{Table depicting the average values of the  Recovered adult population. }
    \centering 
    \begin{tabular}{|l|l|l|l|} 
    \hline
    \textbf{Control Combinations} &  \textbf{Avg Recovered Population(R2)} \\

    \hline
     $u_{11}^*, u_{12}^{*}=0 $&       20.4898 \\
     
    \hline
    $u_{11}^*, u_{12}^{*}$ &        19.9415 \\
    
    \hline
    $u_{11}^*=0, u_{12}^{*}$ &       4.9833 \\

    \hline
     $u_{11}^*=u_{12}^{*} = 0 $&     4.9833\\
    
    \hline
    \end{tabular} 
    \label{t3}
\end{table} \vspace{.4cm}

\begin{center}
	\begin{figure}[hbt!]
		\includegraphics[height = 10cm, width = 17.5cm]{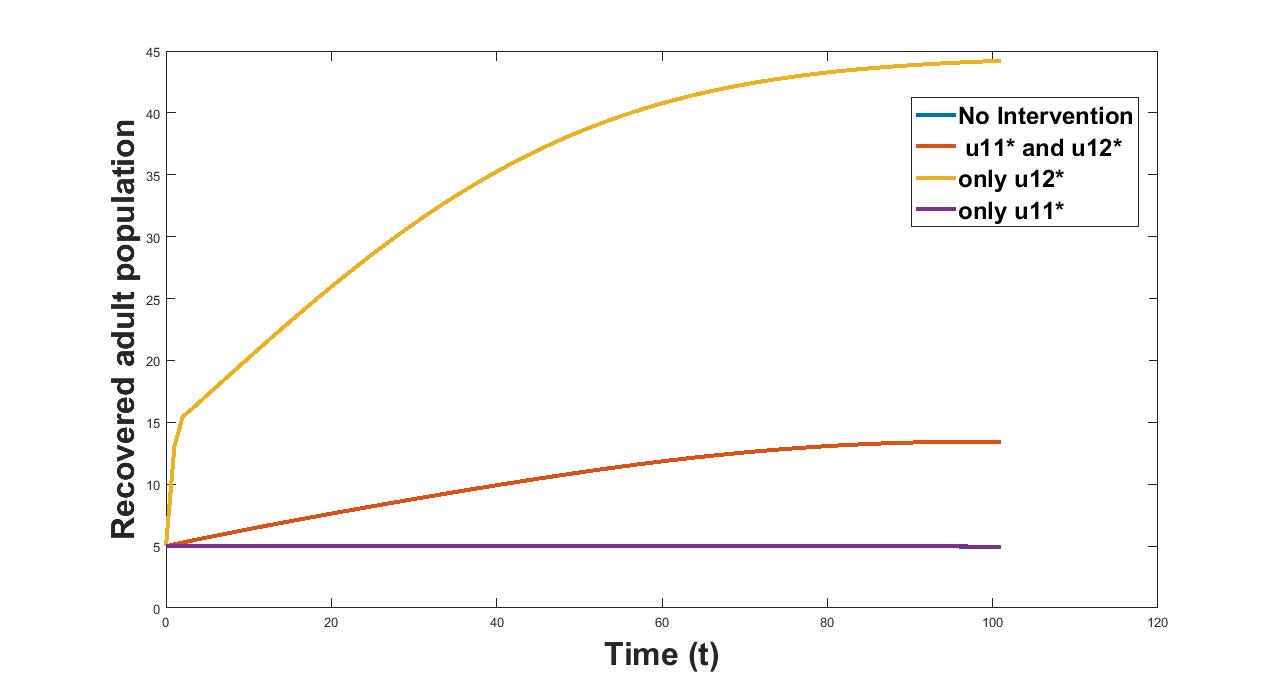} 
	\caption{$R2$ under optimal controls u$_{11}^{*}$, u$_{12}^{*}$}
	\end{figure} 
\end{center}

\begin{center}
	\begin{figure}[hbt!]
		\includegraphics[height = 10cm, width = 17.5cm]{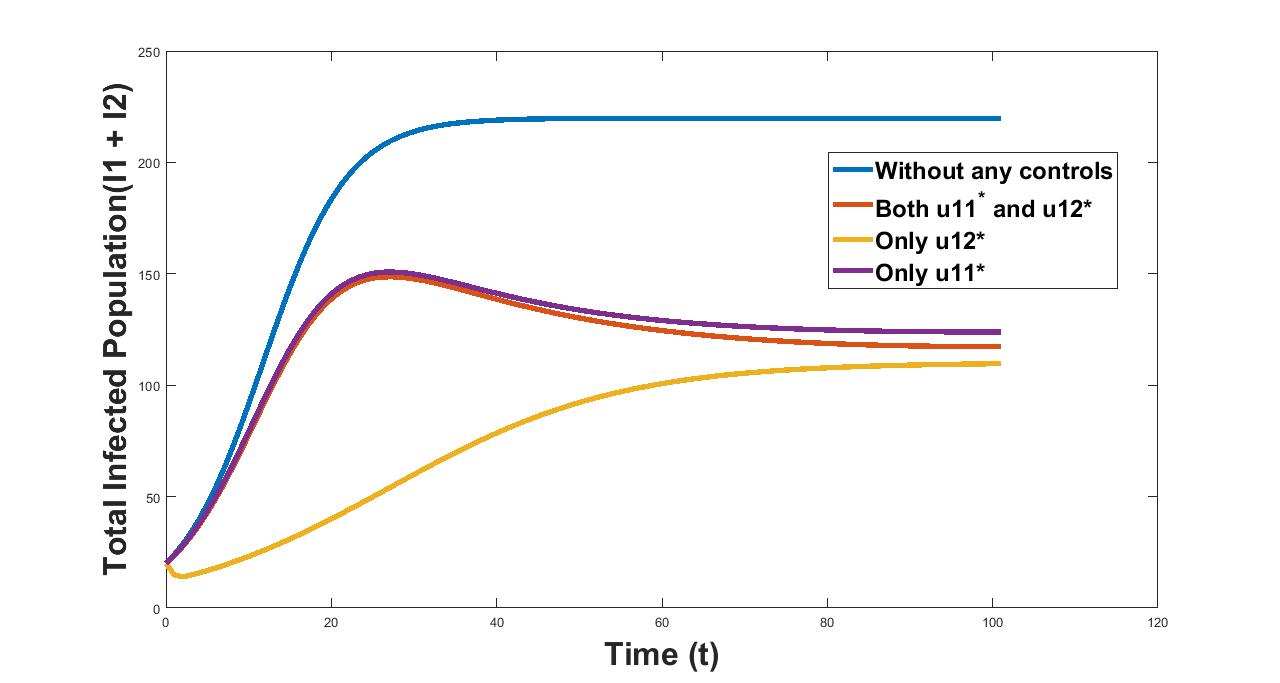} 
	\caption{I1 + I2 under different controls}
	\end{figure} 
\end{center}

\newpage
\section{The Effect of $R_0$ on Optimal Controls and on Disease Burden}
This section devotes to study the dynamics of disease and the effect of optimal controls on the
disease when the basic reproduction number $(R_0)$ varies. Since the severity of the epidemic characterized by the high epidemic peaks which is measured by the higher values of $R_0$, therefore we will observe the prevalence of the cumulative
count of the disease  by varying the basic reproduction
number. For this purpose, we will consider the basic reproduction number $R_0$ in absence of any
controls i.e. $u_{11} = 0 $ and $u_{12} = 0$ as follows for the system (1.1)-(1.6):

\begin{equation*}
         \mathbf{ R_{0}}= \mathbf{\frac{(\beta_{1}S_{1}^*p + \beta_{4}S_{2}^*q) + \sqrt{M}}{2}}
     \end{equation*} where,
 $$p=\frac{1}{d_1+\mu}$$
 $$q=\frac{1}{d_2 + \mu}$$

We choose the parameters same as taken above for the comparative study to explore the impact
of $R_0$ on optimal controls and on the disease. For this purpose, the  total infected population is
plotted to depict the severity of the disease by varying the degree of transmissibility as reflected by the basic reproduction number $R_0$ for different control strategies. 

In figure (16,17,18) we plot the effect of $R_0$ on the infected young population, infected adult and the cumulative infected population considering different controls.  Our findings suggest that when the epidemic is mild $(R_0 \in (1; 1.5))$, the second control $u_{12}$ treatment
works better than the control $u_{11}$ as given in yellow and violet colored curves
in figure 16. Whereas when severity of the epidemic increases $(R_0 \in (1.5; 7))$, the effect of the controls $u_{11}$ and $u_{11}$ and $u_{12}$ together is found to be better than the effect of $u_{12}$ alone. From figure 17 and figure 18 we see that the treatment $u_{12}$ works better in keeping the severity very low throughout the range of the basic reproduction number $(R_0 \in (1; 7))$ as shown in yellow colored curve in the figure 17 followed by the combined effect of both the controls. In figure  19 we plot the cumulative recovered population and see that with the control $u_{12}$, we have the highest recovered population throughout the range of the basic reproduction number $(R_0 \in (1; 7))$.

	\begin{figure}[hbt!]
	\begin{center}
		\includegraphics[width=4in, height=2.3in, angle=0]{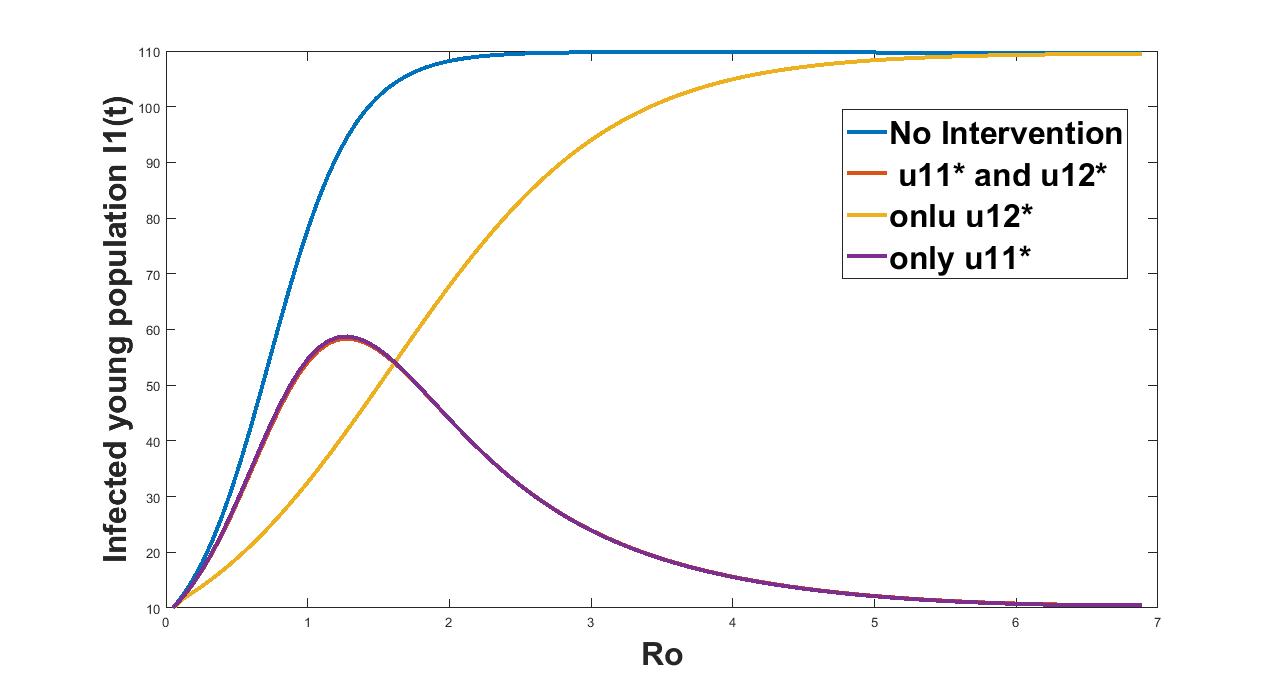} 
	\caption{Effect of $R_0$ on $I_1$  under different controls}
	\end{center}
	\end{figure} 

	\begin{figure}[hbt!]
	\begin{center}
		\includegraphics[width=4in, height=2.3in, angle=0]{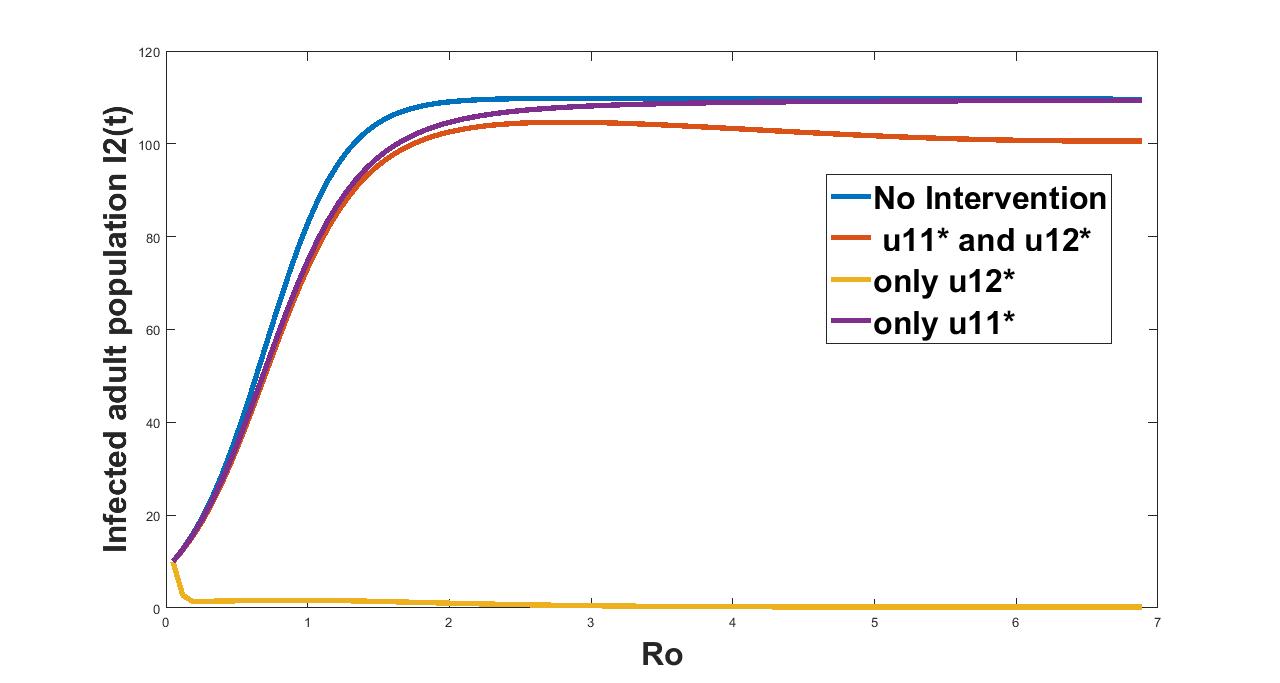}
	\caption{Effect of $R_0$ on $I_2$ under different controls}
	\end{center}
	\end{figure} 

	\begin{figure}[hbt!]
	\begin{center}
		\includegraphics[width=4in, height=2.3in, angle=0]{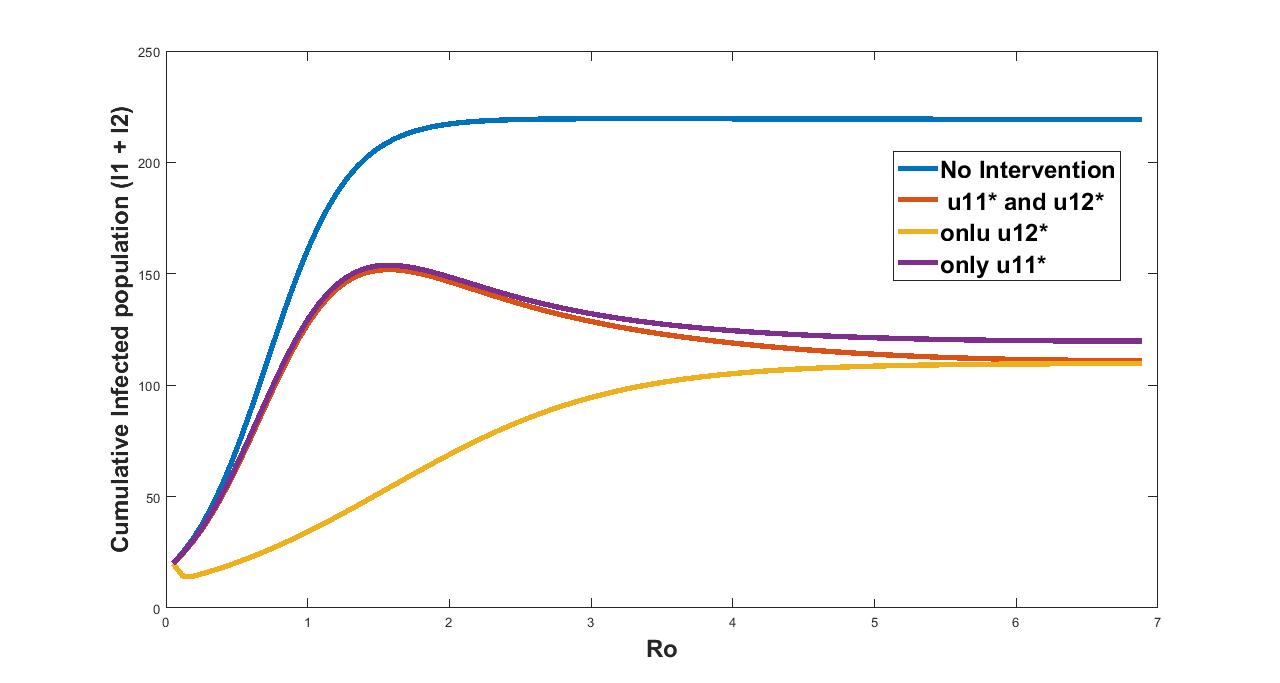}
	\caption{Effect of $R_0$ on cumulative infected population under different controls}
	\end{center}
	\end{figure} 
	\begin{figure}[hbt!]
	\begin{center}
		\includegraphics[width=4in, height=2.3in, angle=0]{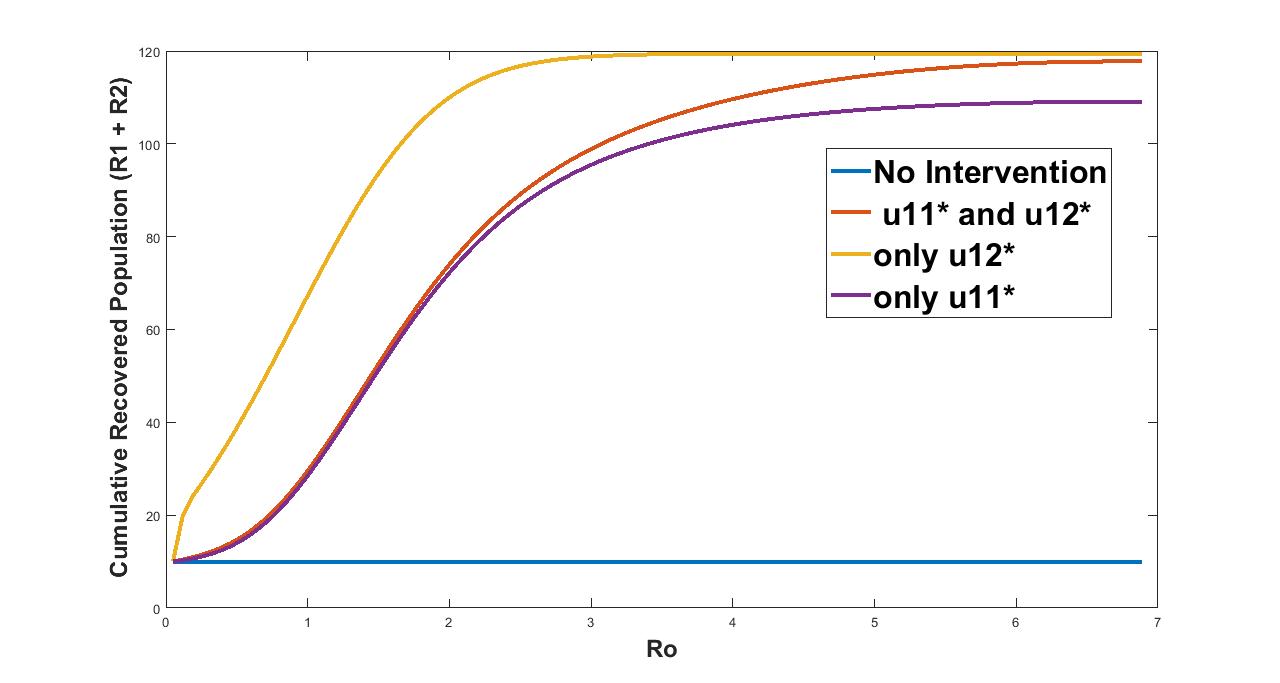}
	\caption{Effect of $R_0$ on cumulative recovered population under different controls}
	\end{center}
	\end{figure} 

\newpage
\section{The Effect of $\alpha$ and $R_0$ on Optimal Controls and  Disease Burden}

In this section we focus on studying the effect of saturation level in treatment by varying the value of $\alpha$ on the disease burden for mild as well as severe epidemic for each designed control strategy. In figure 20 we plot the effect of $\alpha$ on the cumulative count  of infected population considering both the controls $u_{11}$ and $u_{12}$. We observed from figure 20 that the cumulative count of the disease grows significantly as the level of saturation increases for entire range of $R_0$. In figure 21 we consider the control $u_{12}$ alone and study the effect of $\alpha$  by varying the value of $\alpha$ on the disease burden. Here again we see that the cumulative count of the disease grows significantly as the level of saturation increases but as epidemic becomes severe $(R_0 > 5)$ the effect of $\alpha$ remains the same as the cumulative count of infection saturates.

	\begin{figure}[hbt!]
	\begin{center}
		\includegraphics[width=4in, height=2.3in, angle=0]{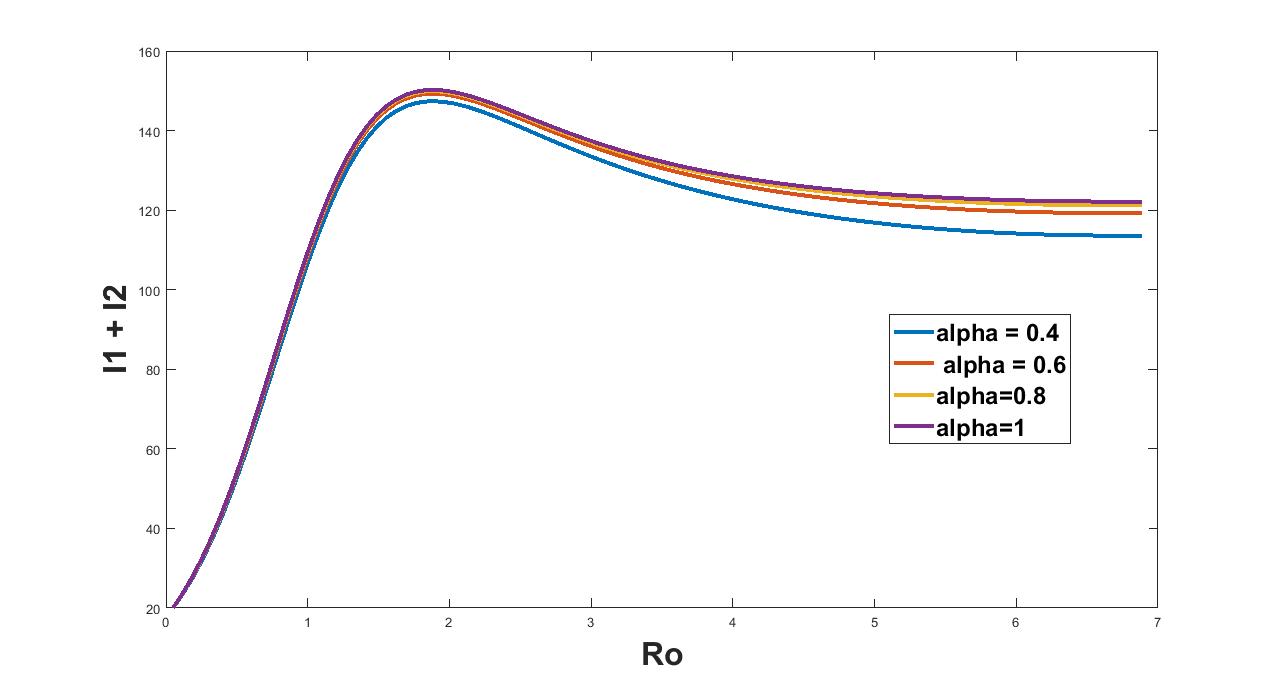}
	\caption{Effect  of $\alpha$ on the cumulative disease burden with $u_{11}^*$ and $u_{12}^*$}
	\end{center}
	\end{figure} 
\begin{figure}[hbt!]
	\begin{center}
		\includegraphics[width=4in, height=2.3in, angle=0]{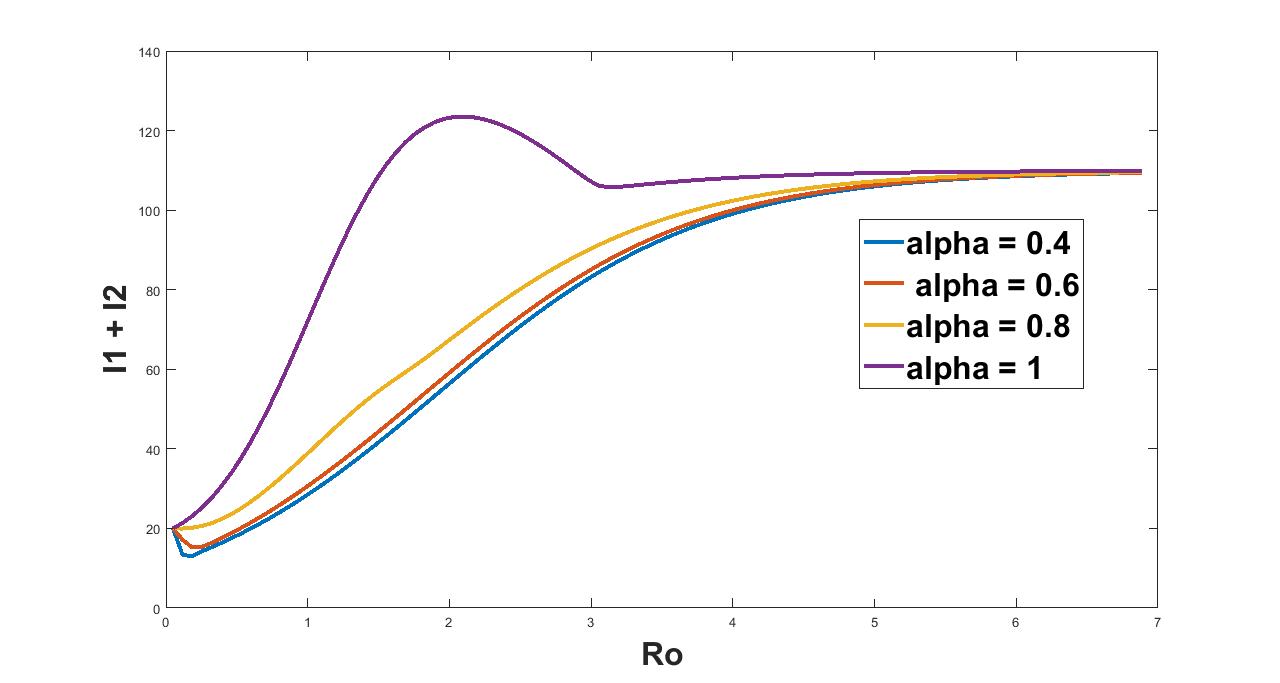}
	\caption{Effect  of $\alpha$ on the cumulative disease burden with  $u_{12}^*$}
	\end{center}
	\end{figure} 

\newpage
\section{Discussion and Conclusion}

Coronavirus disease 2019 (COVID-19) is a contagious respiratory and vascular disease caused by severe acute respiratory syndrome coronavirus 2 (SARS-CoV-2). The Coronavirus caused disease COVID-19 has been declared a pandemic by WHO. As on 20 November 2020, around 55 million have been affected and around 1.3 million have lost their worldwide [3]. 

In this study, initially  we have proposed a non-linear age structured compartmental model in which the population is divided into two age groups the young (age $\leq$30) and the adult ones ($\geq$30) and in each of to age groups we have three compartments, namely susceptible denoted by $S_1$(young) and $S_2$(adult), infected denoted by $I_1$(young) and $I_2$(adult) and recovered denoted by $R_1$(young) and $R_2$(adult). We have used Holling type III recovery rate function of infected adult individuals wherein the treatment provided initially is less, owing to less infected individuals, and as the epidemic progresses, the treatment increases accordingly and non-linearly saturates due to limitations on medical facilities. Later we proposed an  optimal  control problem to investigate the role of the control strategies namely the role of treatments in reducing the infection.

From the Stability analysis we conclude that the infection free equilibrium remains asymptotically stable whenever $R_0 < 1$ and as $R_0$ crosses unity we have the infected equilibrium to be stable.
From the sensitivity analysis of the parameters of the model we conclude that the parameters $u_{11}, b_1, \beta_1, d_1$ and $\mu$ are the only sensitive parameters in some intervals as discussed in table 6.
Findings from  the Optimal Control studies suggests that the infection among the adult population(age $\geq 30)$ is least considering the second control $u_{12}$ whereas, when both the controls $u_{11}$ and $u_{12}$ are considered together the infectives is minimum in case of young populations(age $ \leq 30$). The cumulative infected population reduced the maximum when the second control was considered followed by considering both the controls together.

 We have also studied the effect of $R_0$ and $\alpha$ on the disease burden considering different control strategies. Our findings suggests that when the epidemic is mild $(R_0 \in (1; 1.5))$, the second control $u_{12}$ treatment
works better than the control $u_{11}$ as given in yellow and violet colored curves
in figure 16. Whereas when severity of the epidemic increased $(R_0 \in (1.5; 7))$, the effect of the controls $u_{11}$ and $u_{11}$ and $u_{12}$ together was found to be better than the effect of $u_{12}$ alone. From figure 17 and figure 18 we see that the treatment $u_{12}$ works better in keeping the severity very low throughout the range of the basic reproduction number $(R_0 \in (1; 7))$ as shown in yellow colored curve in the figure 17 followed by the combined effect of both the controls.

From the study on the effect of $\alpha$ on disease burden we saw that the cumulative count of the disease grew significantly as the level of saturation increased for entire range of $R_0$ considering both the controls together. Whereas when we considered the second control alone we saw that the cumulative count of the disease grew significantly as the level of saturation increased but as epidemic became severe $(R_0 > 5)$ the effect of $\alpha$ remained the same.

\newpage
\bibliographystyle{amsplain}
\bibliography{references}

\newpage
	\section{appendix - A}
	{\flushleft{  \textbf{SENSITIVITY PLOTS FOR OTHER PARAMETERS} }}


		\subsection{Parameter $\boldsymbol{u_{11}}$}
		
\begin{figure}[hbt!]
\begin{center}
\includegraphics[width=2.2in, height=1.8in, angle=0]{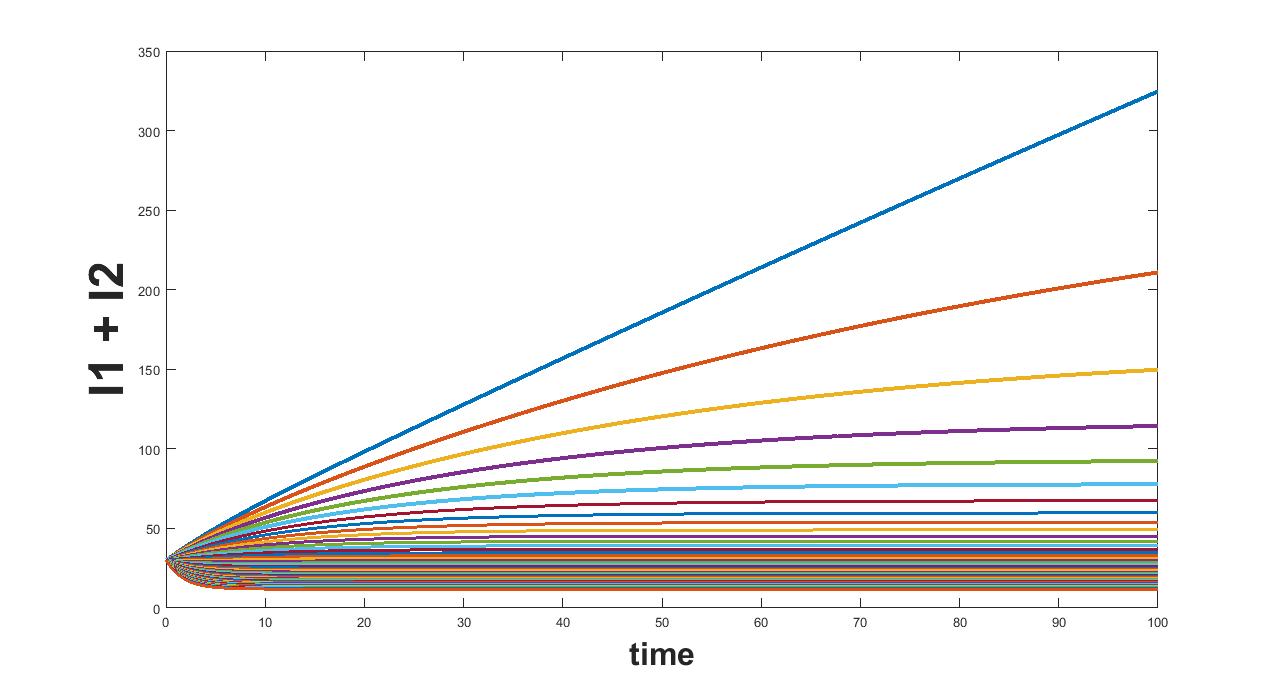}
\hspace{-.4cm}
\includegraphics[width=2.2in, height=1.8in, angle=0]{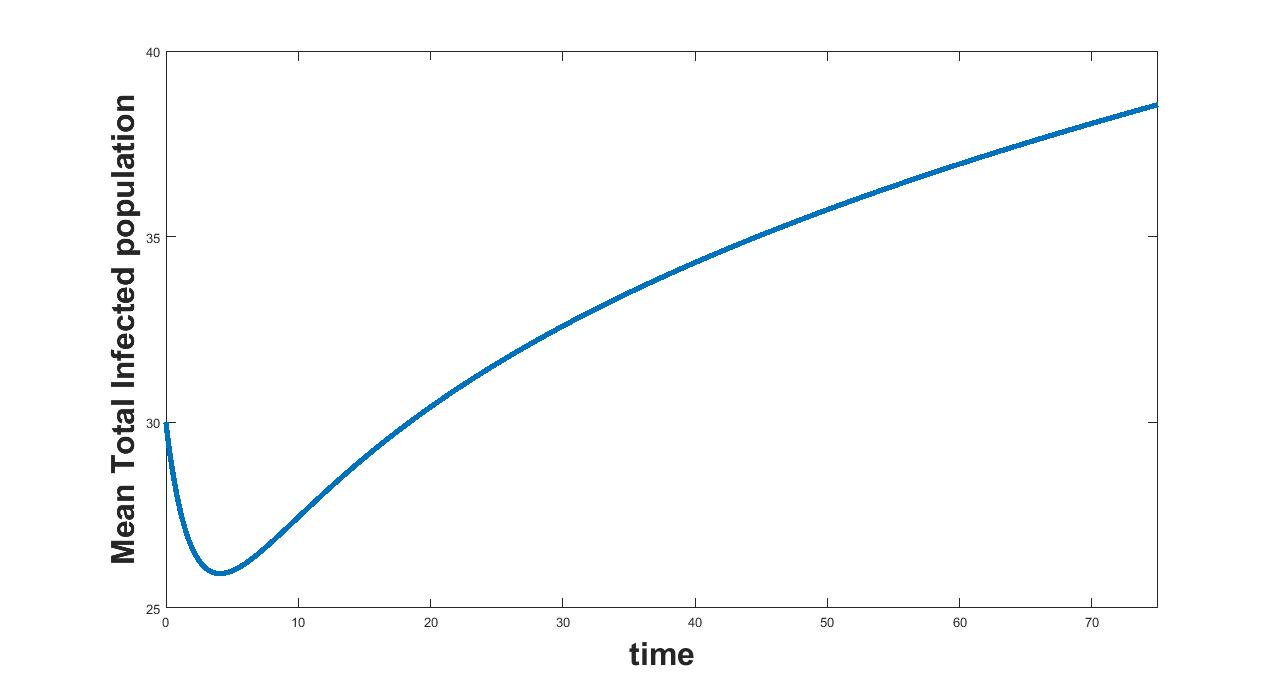}
\hspace{-.395cm}
\includegraphics[width=2.2in, height=1.8in, angle=0]{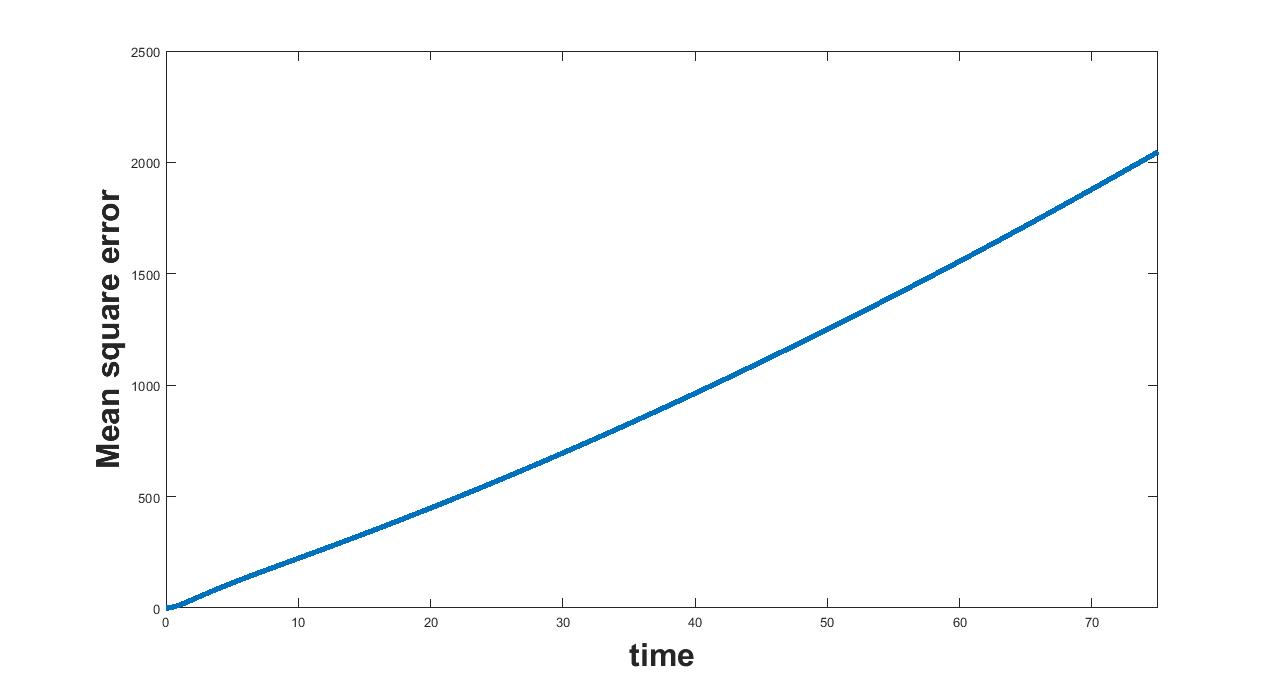}
\caption*{(a) Interval I : 0 to 0.5}
\end{center}
\end{figure}

\vspace{-3mm}

\begin{figure}[hbt!]
\begin{center}
\includegraphics[width=2.2in, height=1.8in, angle=0]{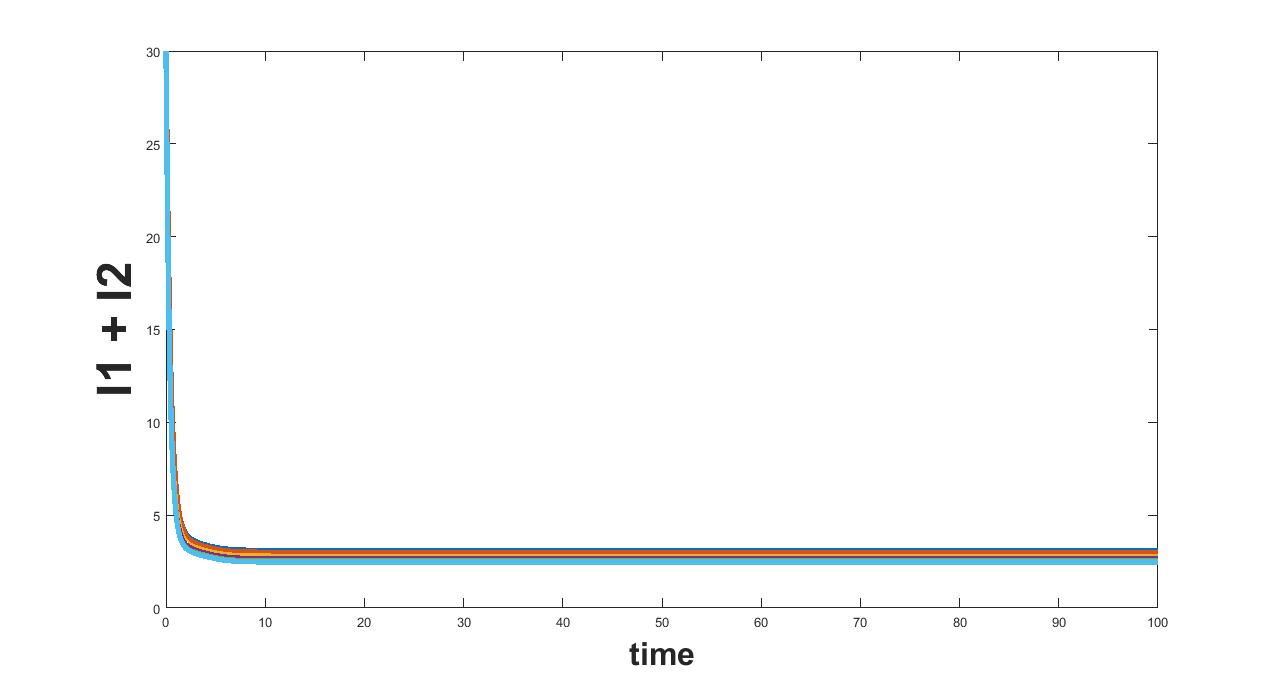}
\hspace{-.4cm}
\includegraphics[width=2.2in, height=1.8in, angle=0]{u11mean2}
\hspace{-.395cm}
\includegraphics[width=2.2in, height=1.8in, angle=0]{u11error1}
\caption*{(b) Interval II : 1.5 to 2}

\vspace{5mm}
\caption{Sensitivity Analysis of $u_{11}$}
\label{sen_beta}
\end{center}
\end{figure}
		
		\newpage
		\subsection{Parameter $\boldsymbol{\beta_2}$}
		
\begin{figure}[hbt!]
\begin{center}
\includegraphics[width=2.2in, height=1.8in, angle=0]{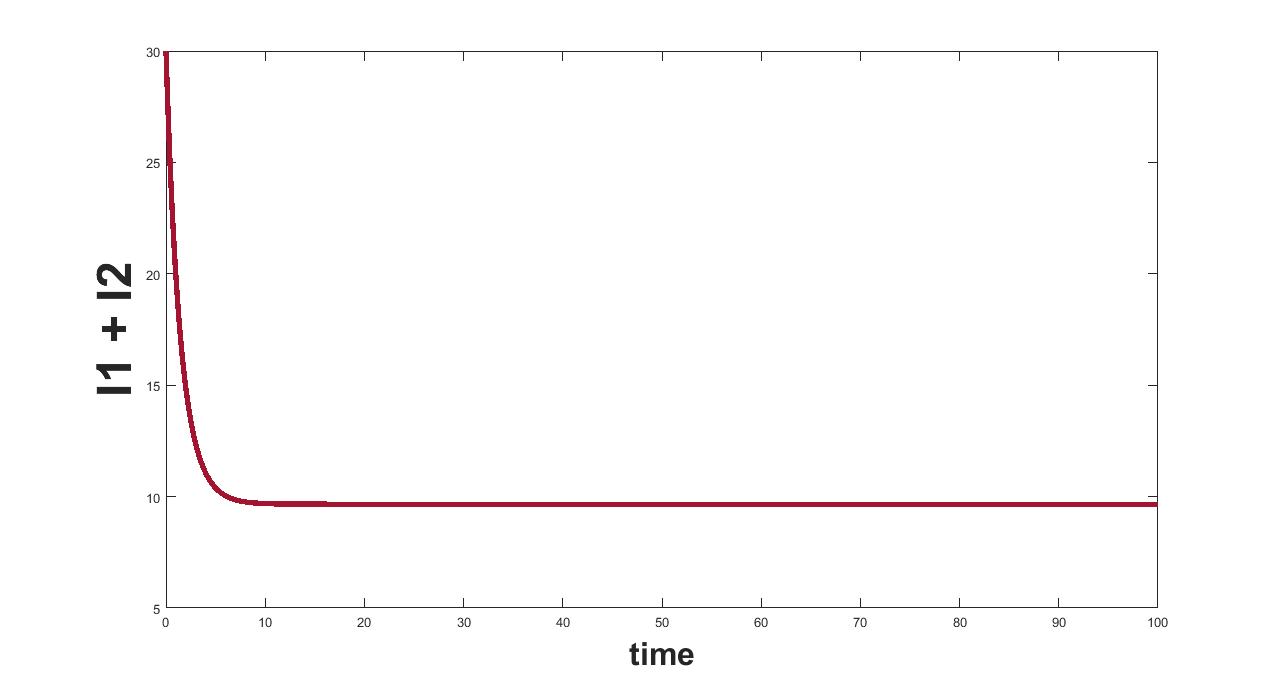}
\hspace{-.4cm}
\includegraphics[width=2.2in, height=1.8in, angle=0]{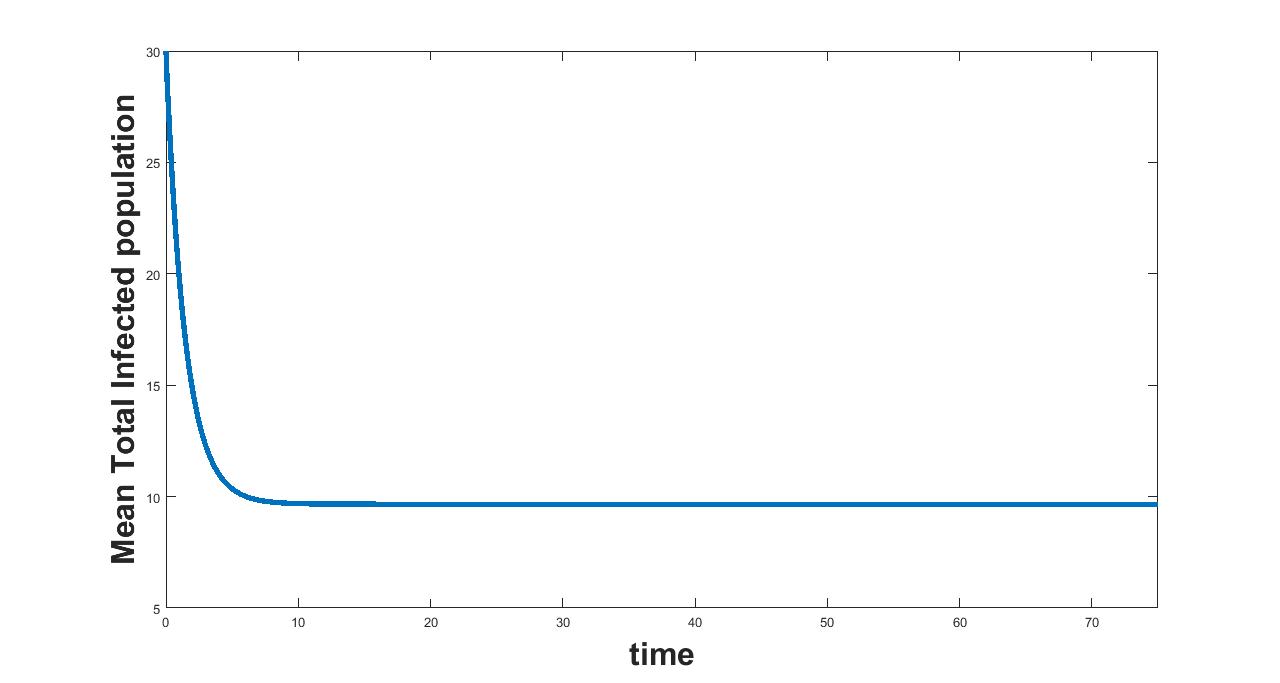}
\hspace{-.395cm}
\includegraphics[width=2.2in, height=1.8in, angle=0]{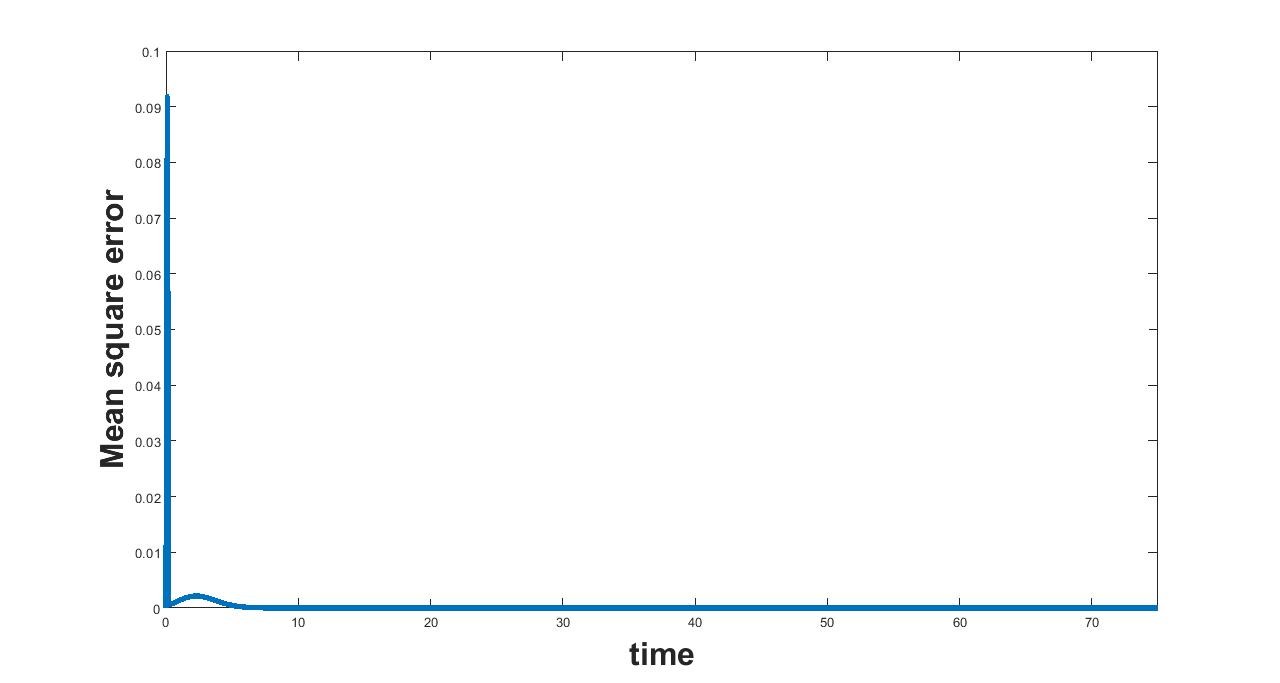}
\caption*{(a) Interval I : 0 to 2}
\end{center}
\end{figure}

\vspace{-3mm}

\begin{figure}[hbt!]
\begin{center}
\includegraphics[width=2.2in, height=1.8in, angle=0]{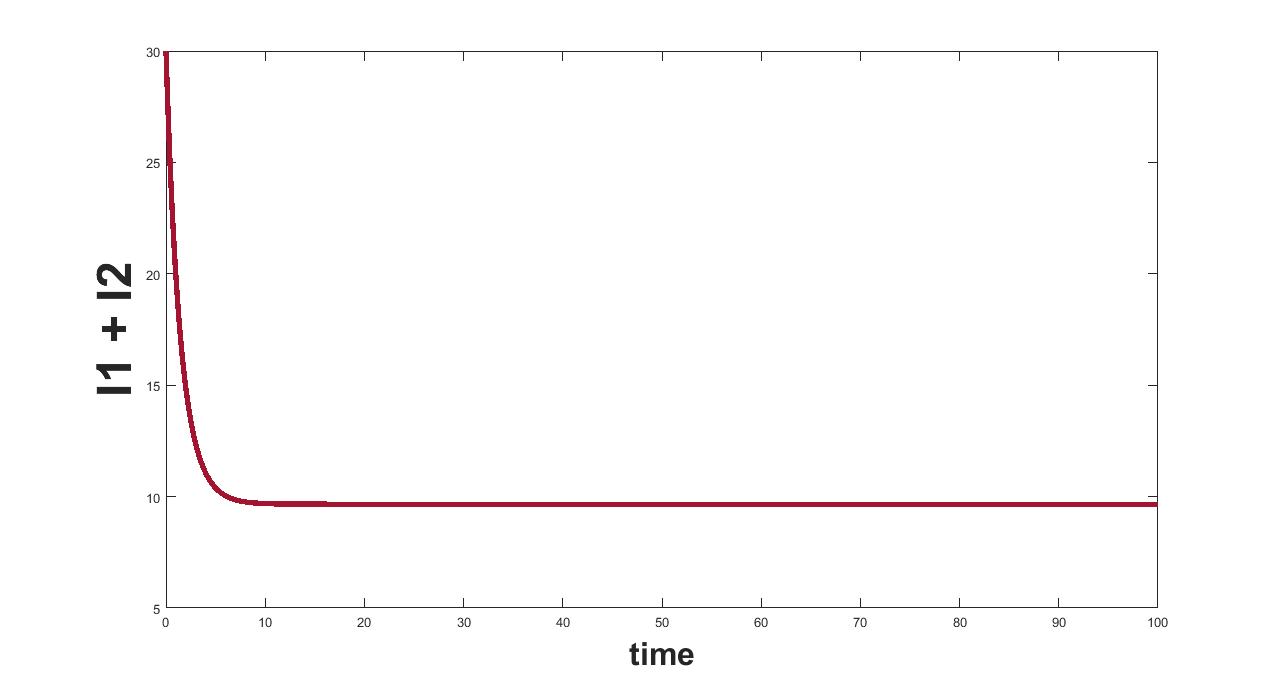}
\hspace{-.4cm}
\includegraphics[width=2.2in, height=1.8in, angle=0]{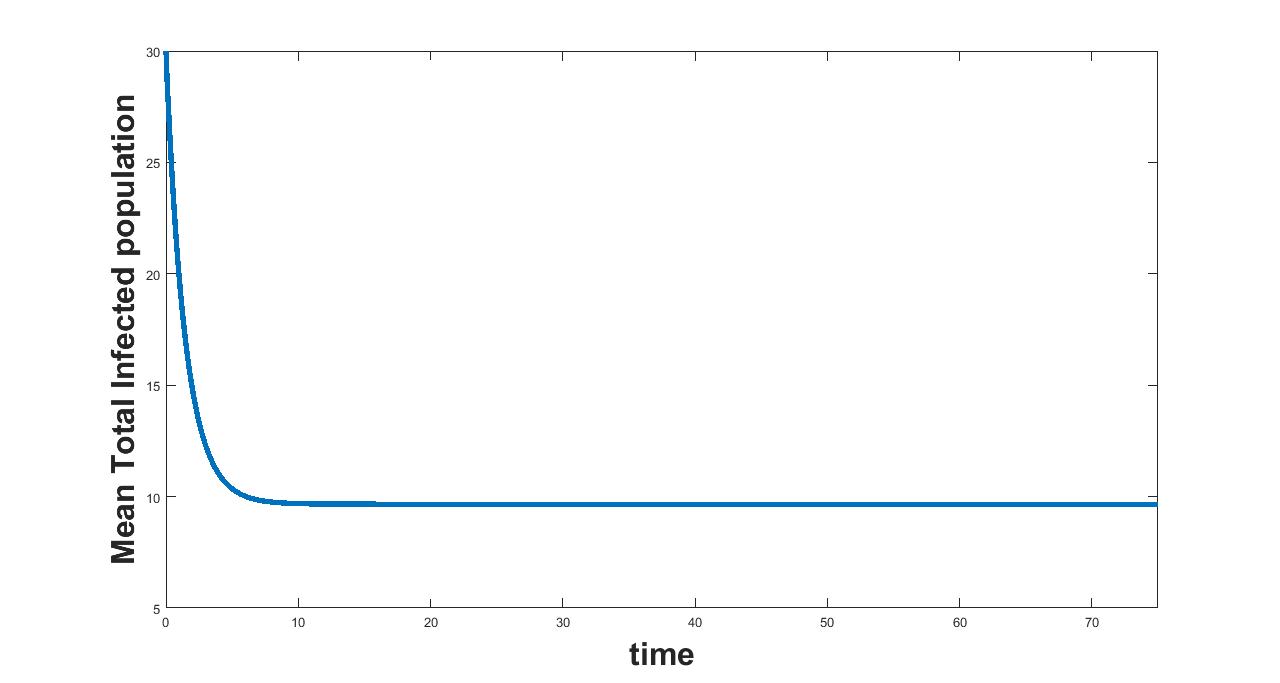}
\hspace{-.395cm}
\includegraphics[width=2.2in, height=1.8in, angle=0]{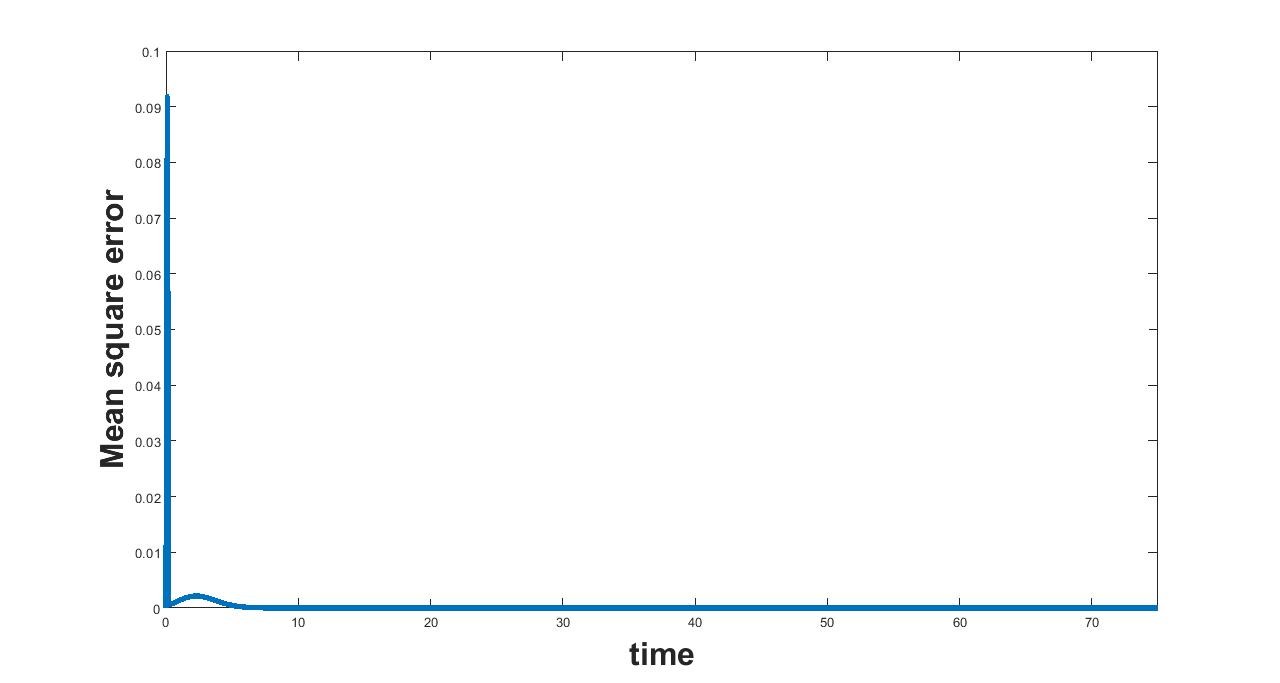}
\caption*{(b) Interval I : 2 to 3}

\vspace{5mm}
\caption{Sensitivity Analysis of $\beta_2$}
\label{sen_beta}
\end{center}
\end{figure}

		
		\subsection{Parameter $\boldsymbol{\beta_3}$}
		
\begin{figure}[hbt!]
\begin{center}
\includegraphics[width=2.2in, height=1.8in, angle=0]{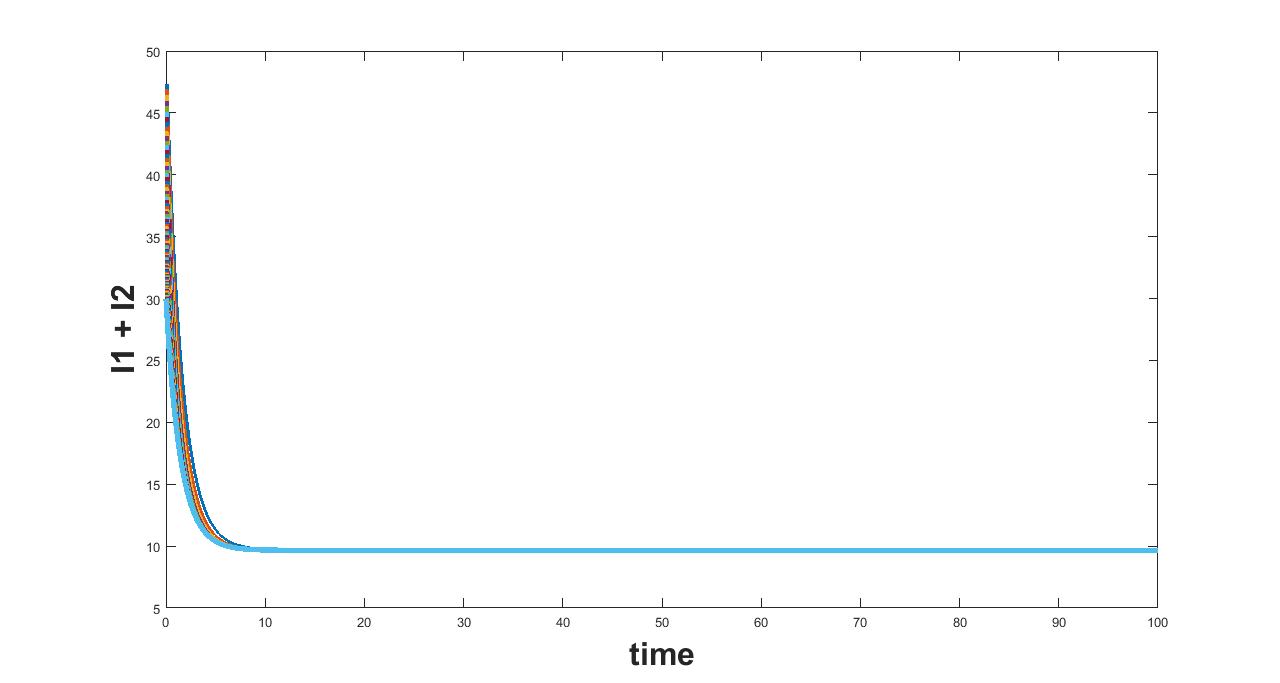}
\hspace{-.4cm}
\includegraphics[width=2.2in, height=1.8in, angle=0]{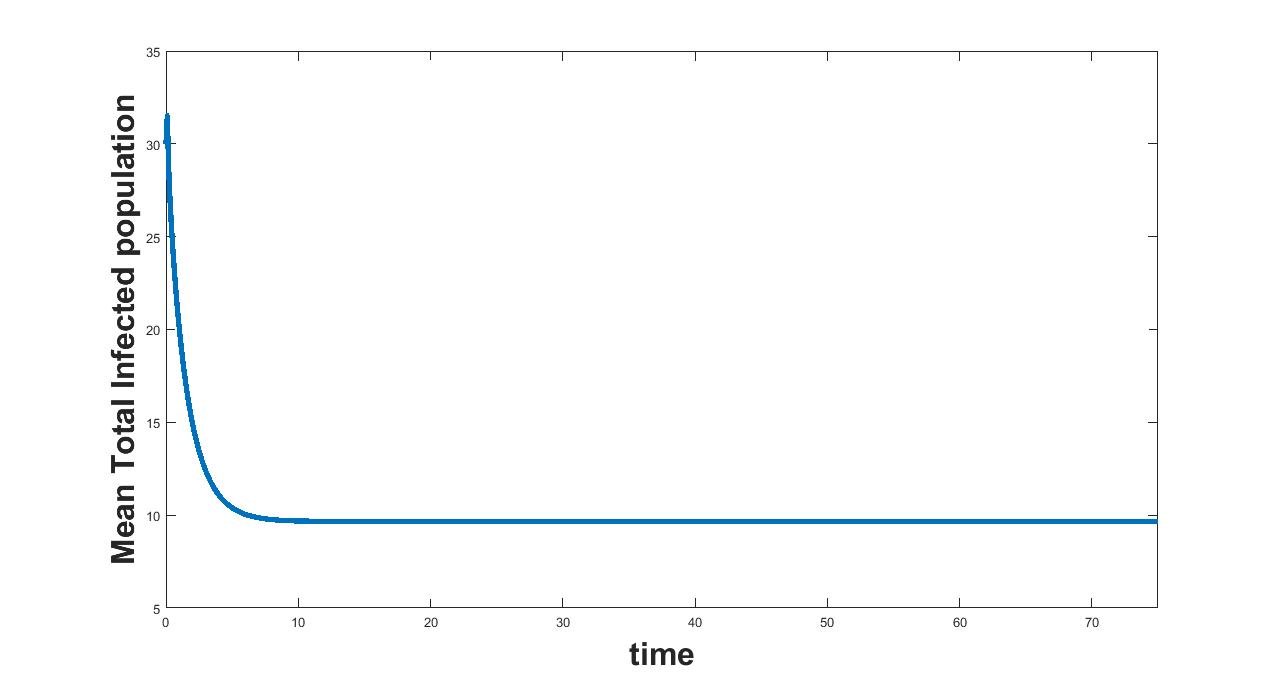}
\hspace{-.395cm}
\includegraphics[width=2.2in, height=1.8in, angle=0]{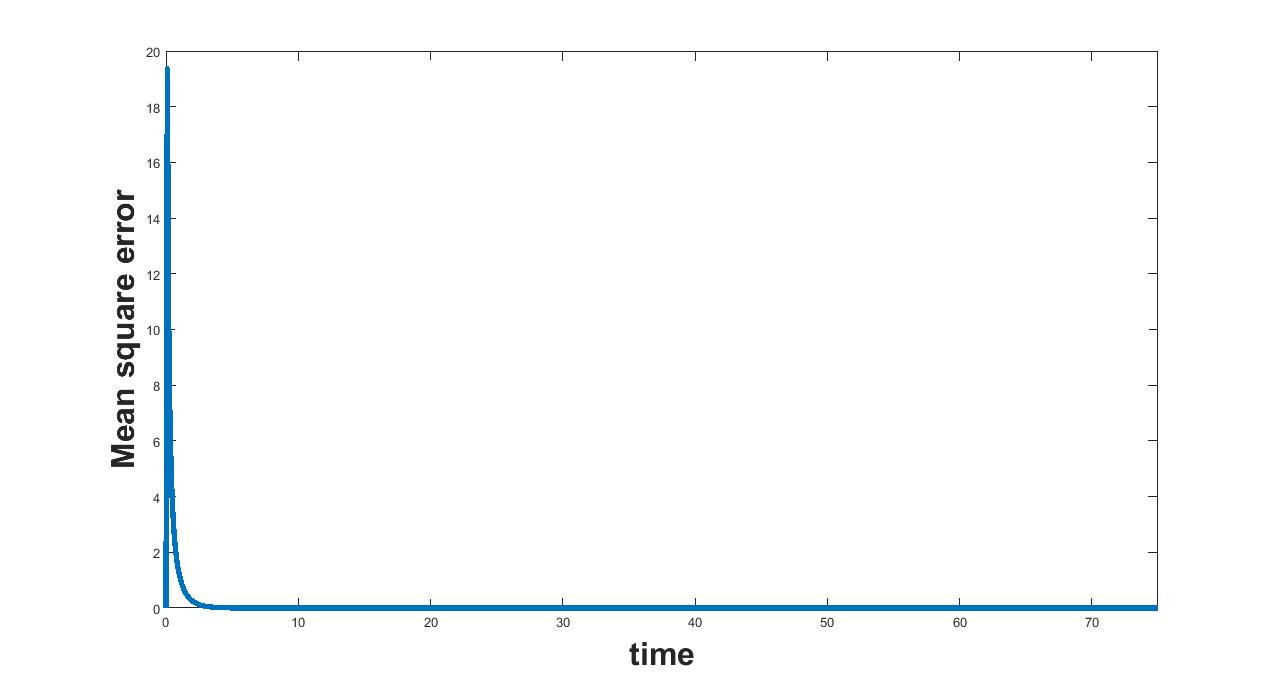}
\caption*{(a) Interval I : 0 to 2.5}
\end{center}
\end{figure}

\vspace{-3mm}

\begin{figure}[hbt!]
\begin{center}
\includegraphics[width=2.2in, height=1.8in, angle=0]{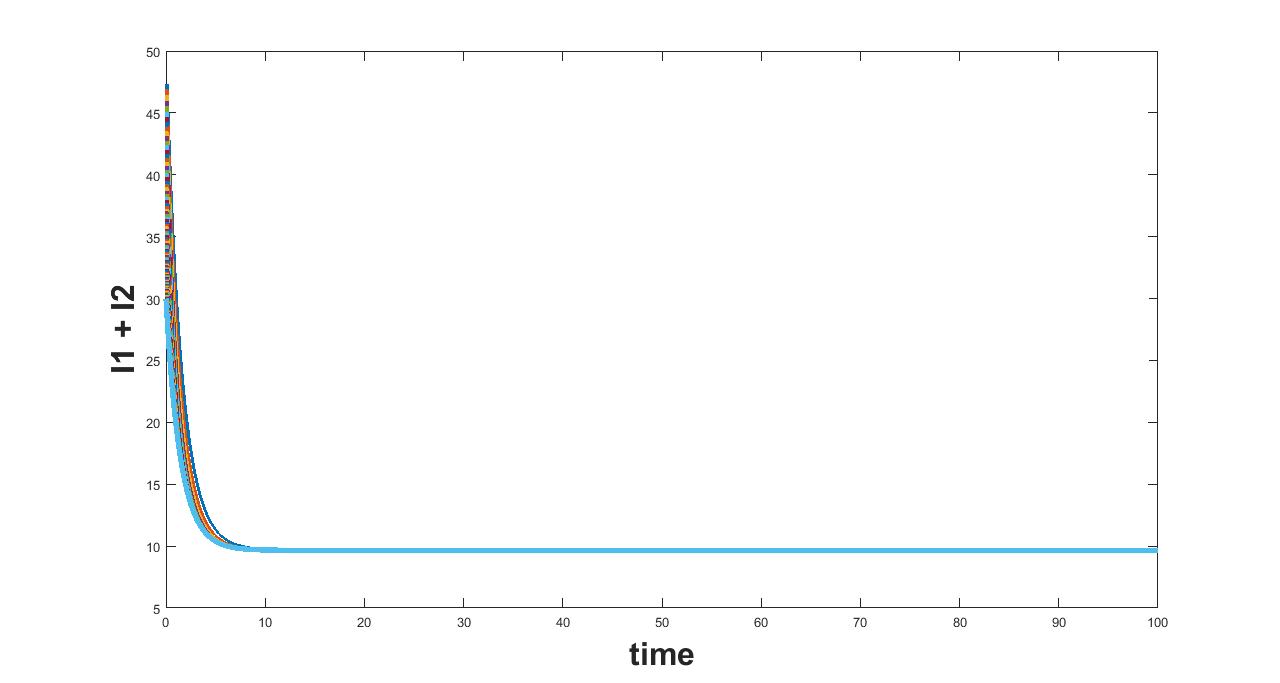}
\hspace{-.4cm}
\includegraphics[width=2.2in, height=1.8in, angle=0]{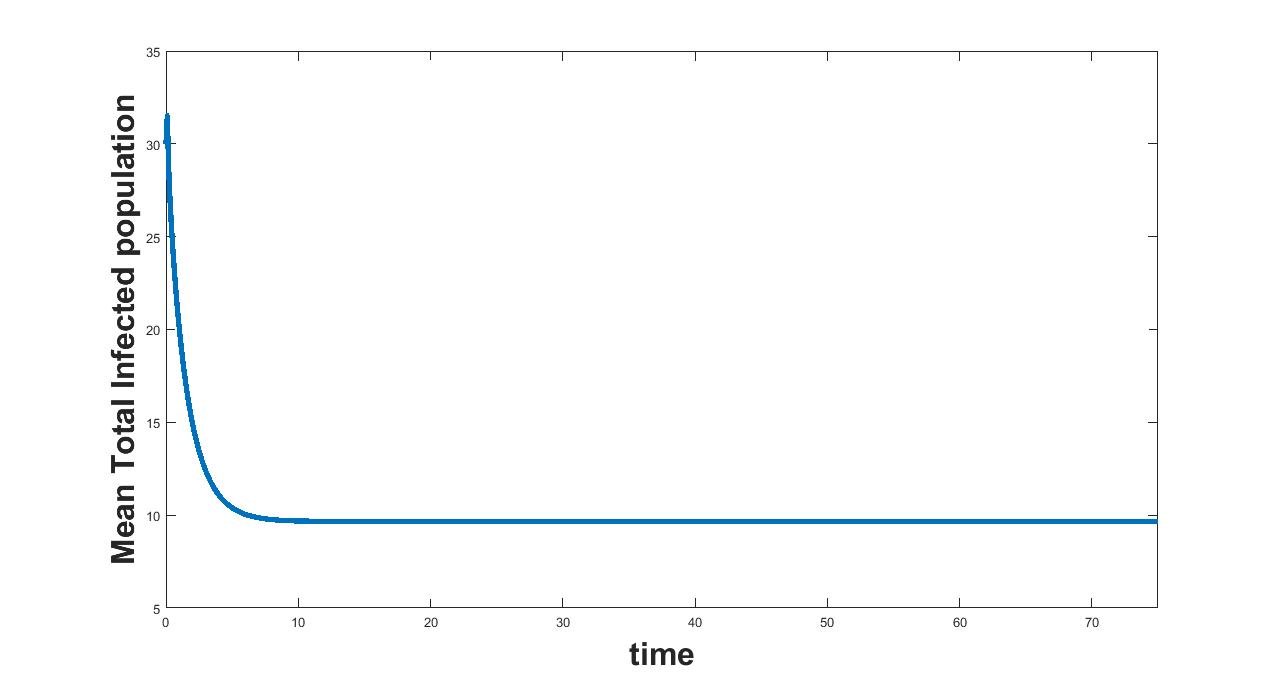}
\hspace{-.395cm}
\includegraphics[width=2.2in, height=1.8in, angle=0]{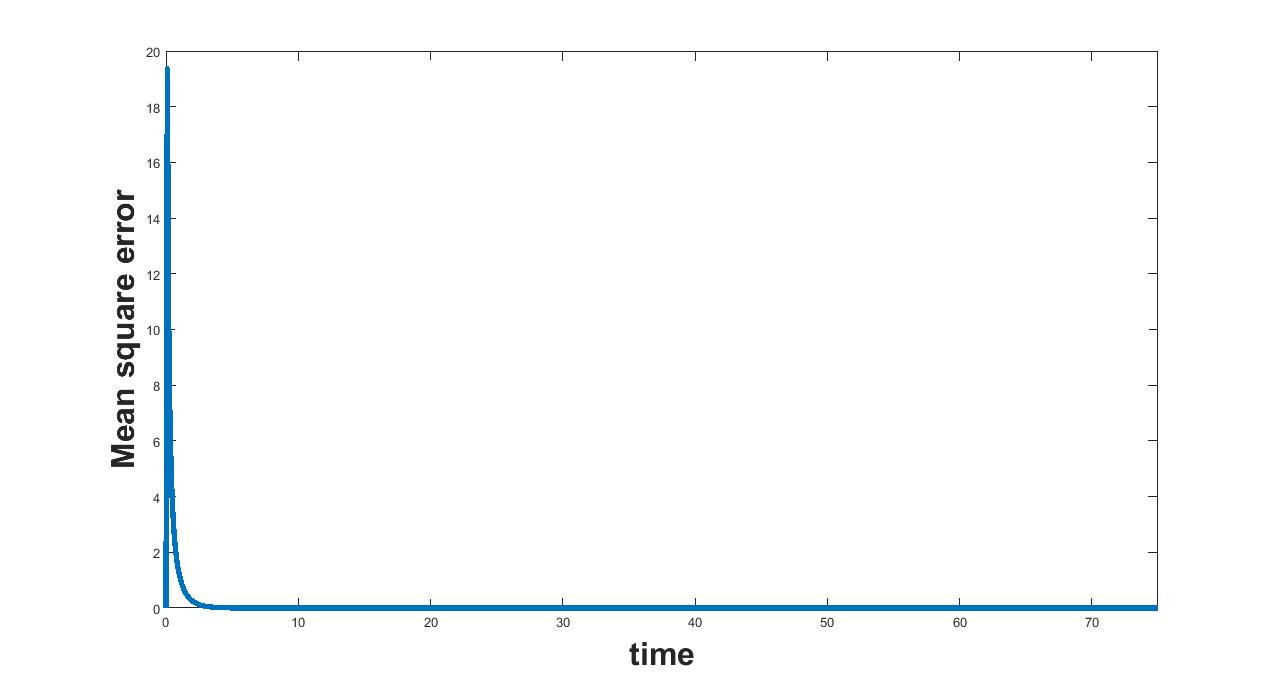}
\caption*{(b) Interval I : 2.5 to 5}

\vspace{5mm}
\caption{Sensitivity Analysis of $\beta_3$}
\label{sen_beta}
\end{center}
\end{figure}

\newpage
		\subsection{Parameter $\boldsymbol{\beta_4}$}
		
\begin{figure}[hbt!]
\begin{center}
\includegraphics[width=2.2in, height=1.8in, angle=0]{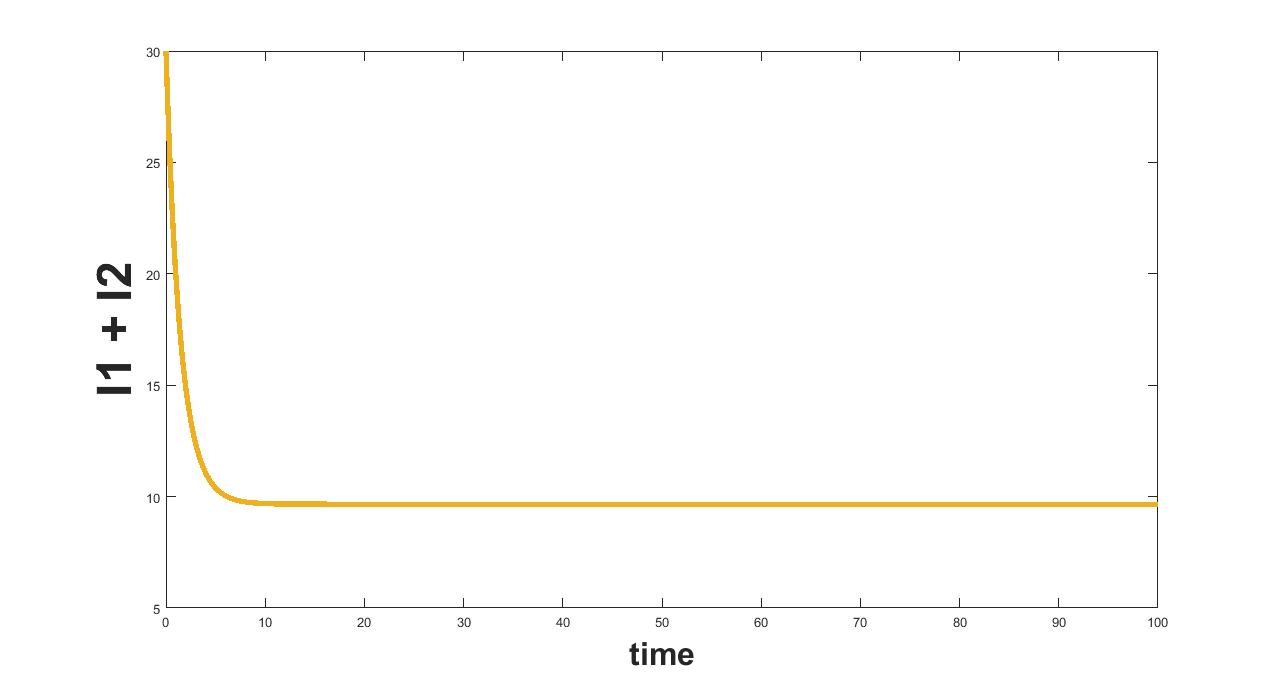}
\hspace{-.4cm}
\includegraphics[width=2.2in, height=1.8in, angle=0]{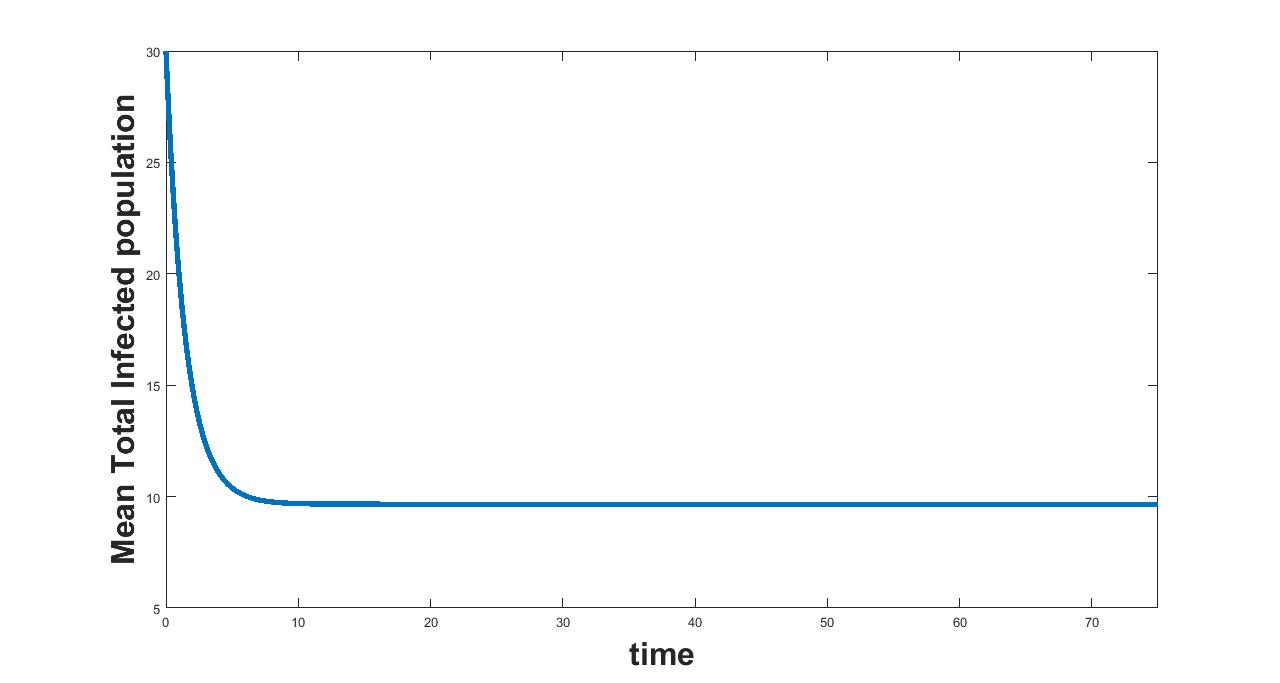}
\hspace{-.395cm}
\includegraphics[width=2.2in, height=1.8in, angle=0]{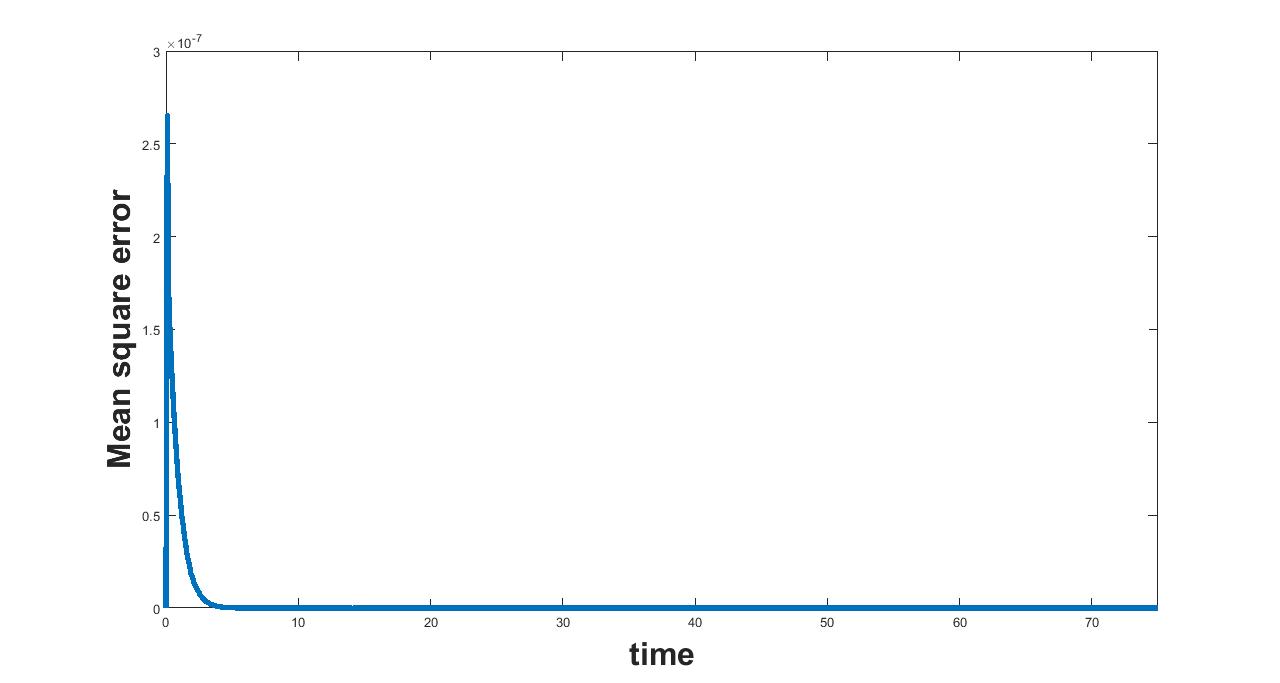}
\caption*{(a) Interval I : 0 to .5}
\end{center}
\end{figure}

\vspace{-3mm}

\begin{figure}[hbt!]
\begin{center}
\includegraphics[width=2.2in, height=1.8in, angle=0]{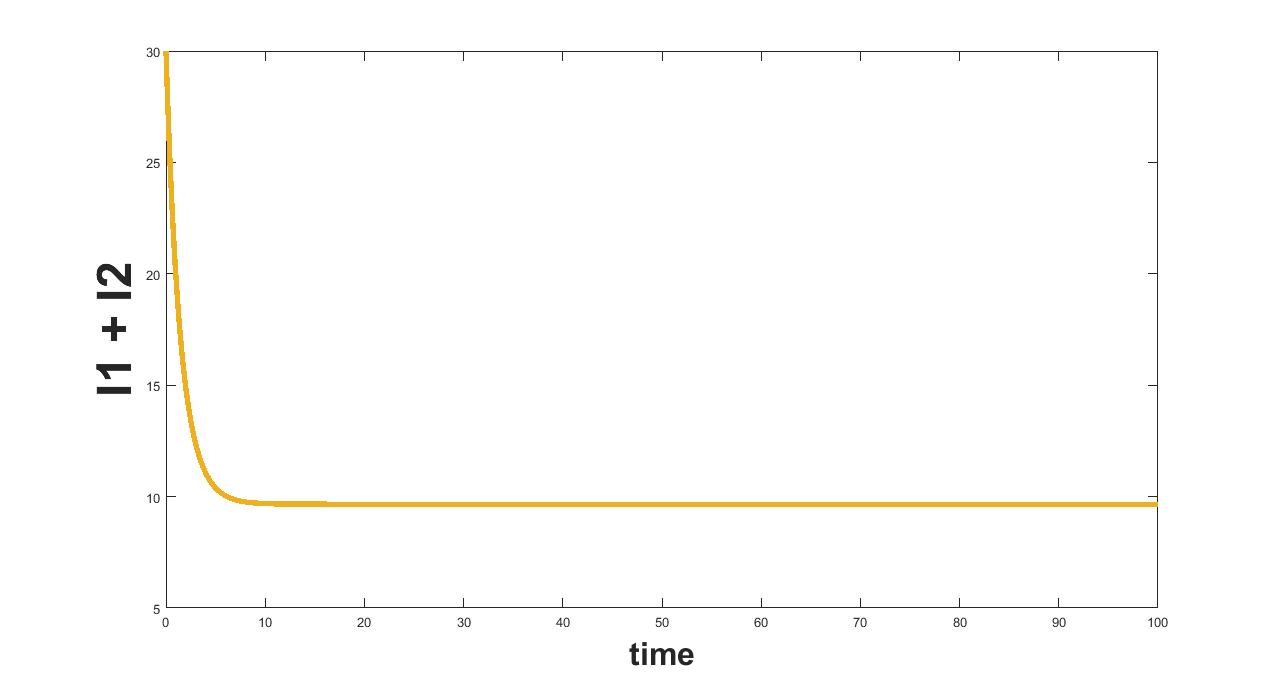}
\hspace{-.4cm}
\includegraphics[width=2.2in, height=1.8in, angle=0]{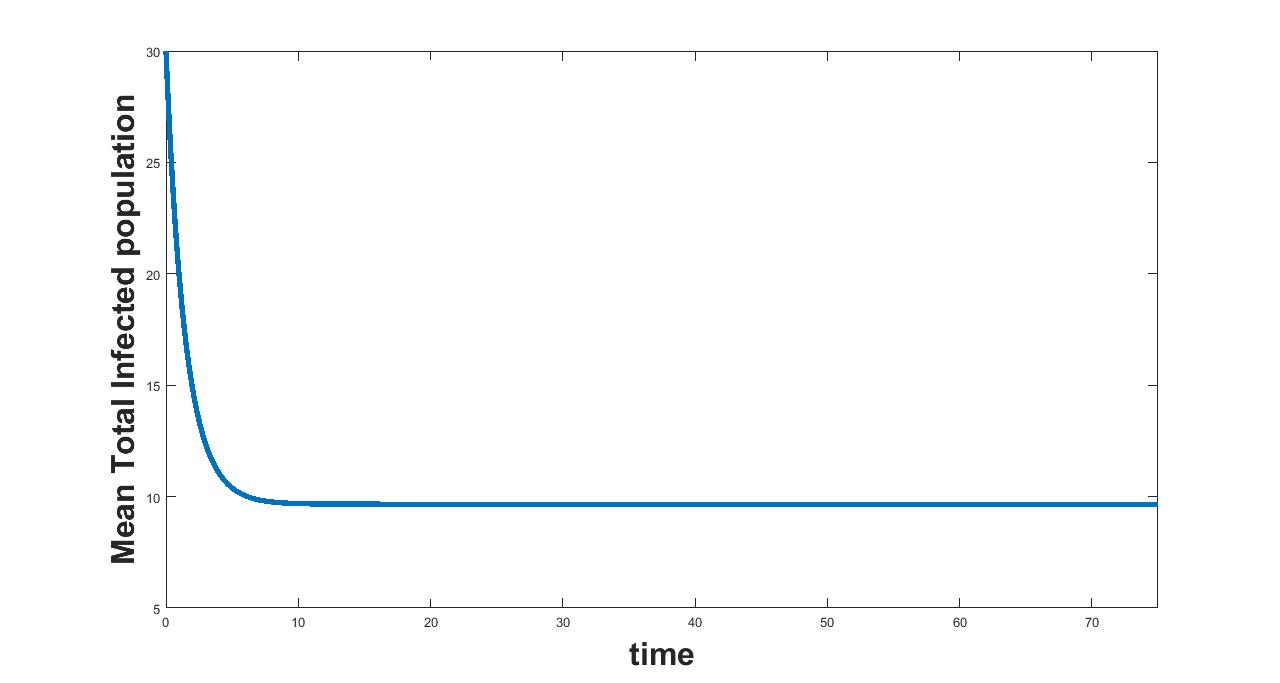}
\hspace{-.395cm}
\includegraphics[width=2.2in, height=1.8in, angle=0]{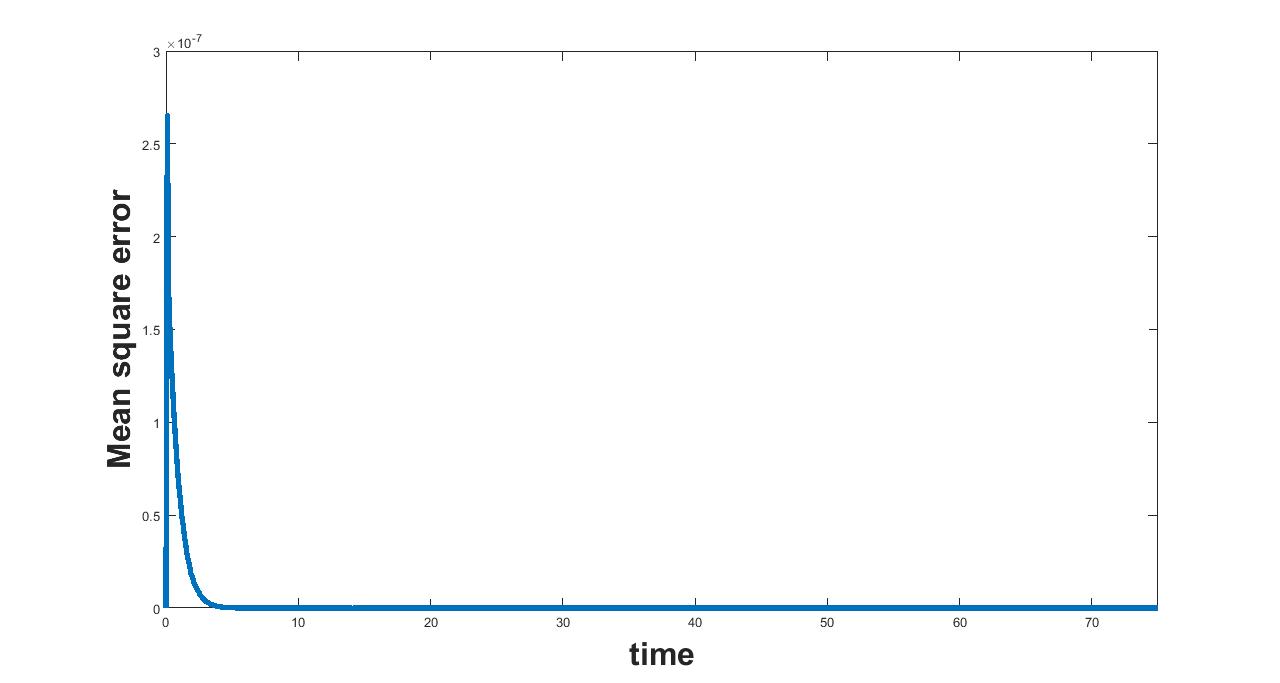}
\caption*{(a) Interval I : 0.5 to 1}

\vspace{5mm}
\caption{Sensitivity Analysis of $\beta_4$}
\label{sen_beta}
\end{center}
\end{figure}

		\newpage
		\subsection{Parameter $\boldsymbol{d_1}$}
		
\begin{figure}[hbt!]
\begin{center}
\includegraphics[width=2.2in, height=1.8in, angle=0]{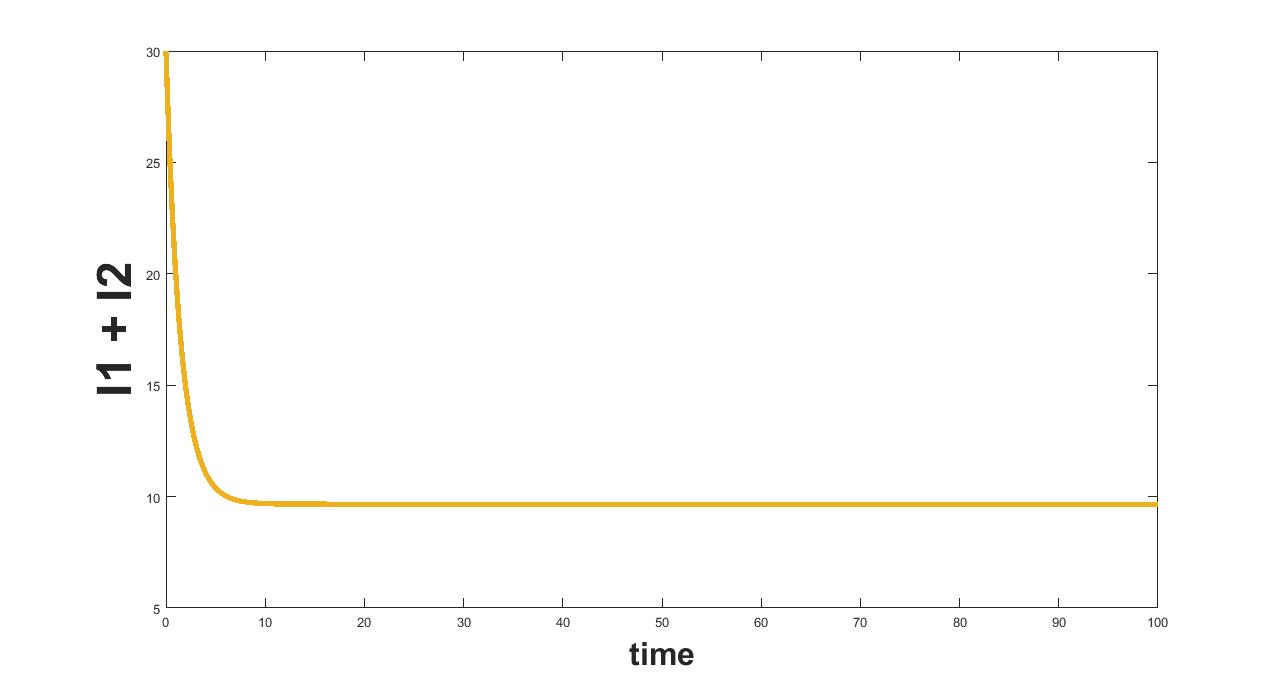}
\hspace{-.4cm}
\includegraphics[width=2.2in, height=1.8in, angle=0]{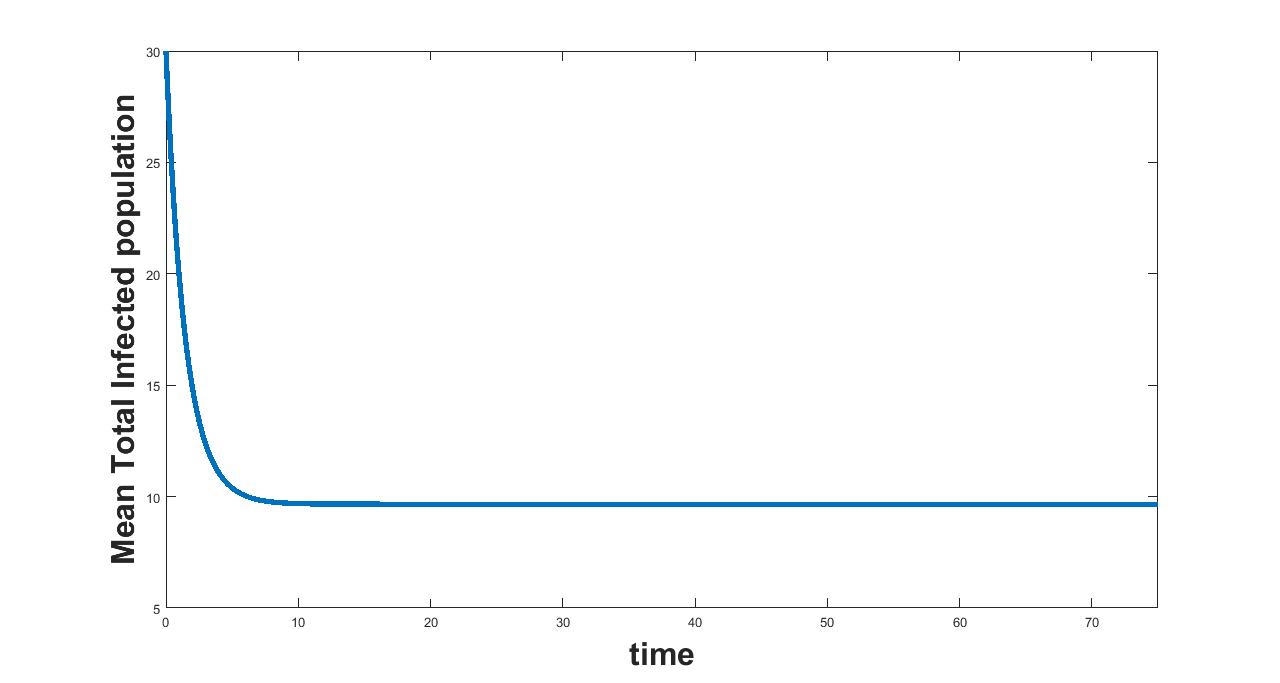}
\hspace{-.395cm}
\includegraphics[width=2.2in, height=1.8in, angle=0]{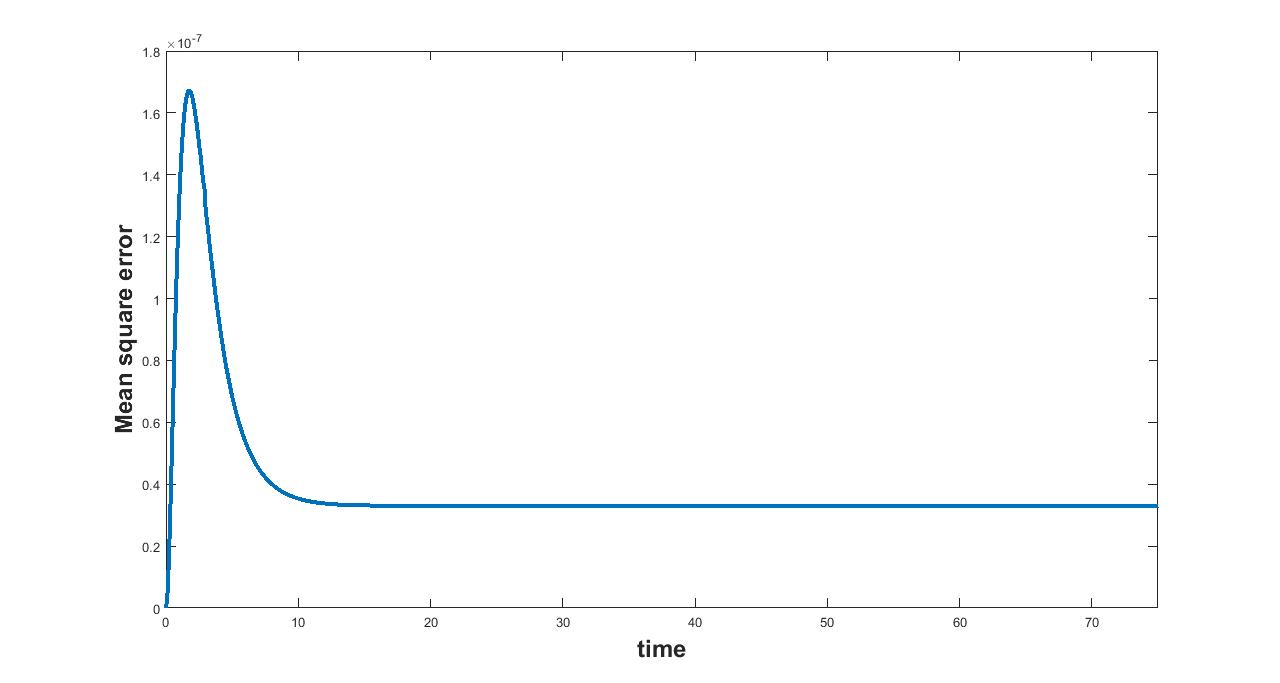}
\caption*{(a) Interval I : 0 to 0.000073 }
\end{center}
\end{figure}

\vspace{-3mm}

\begin{figure}[hbt!]
\begin{center}
\includegraphics[width=2.2in, height=1.8in, angle=0]{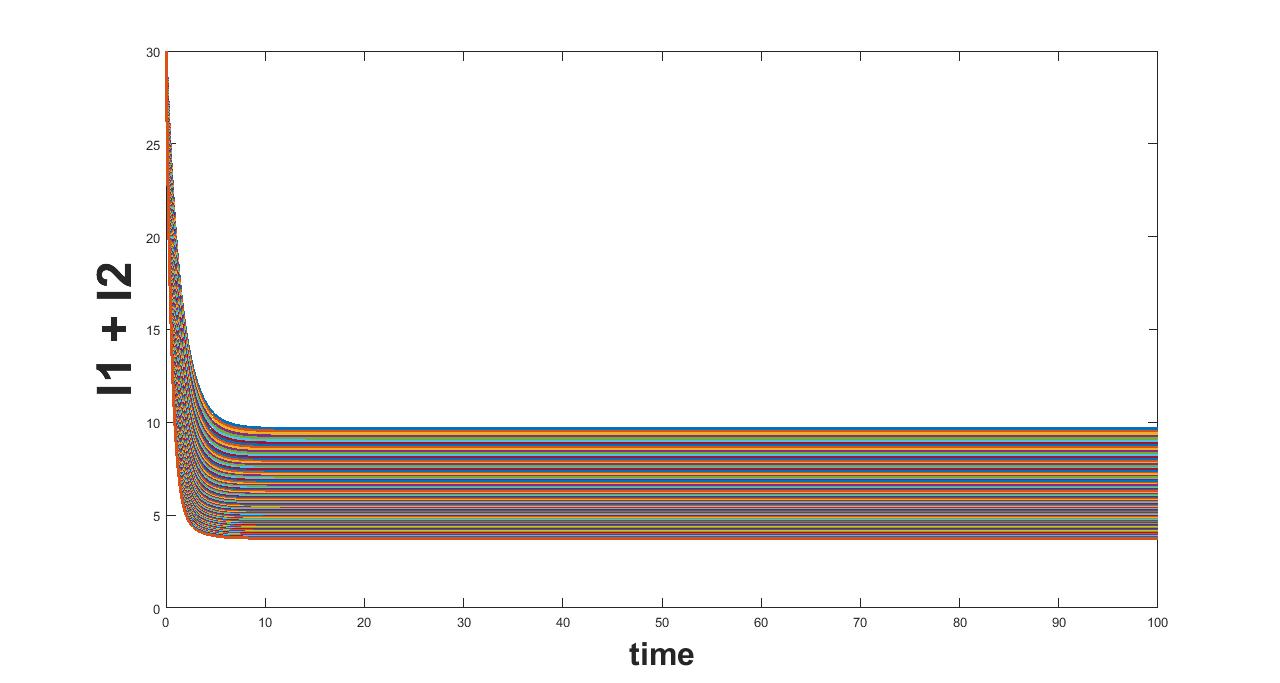}
\hspace{-.4cm}
\includegraphics[width=2.2in, height=1.8in, angle=0]{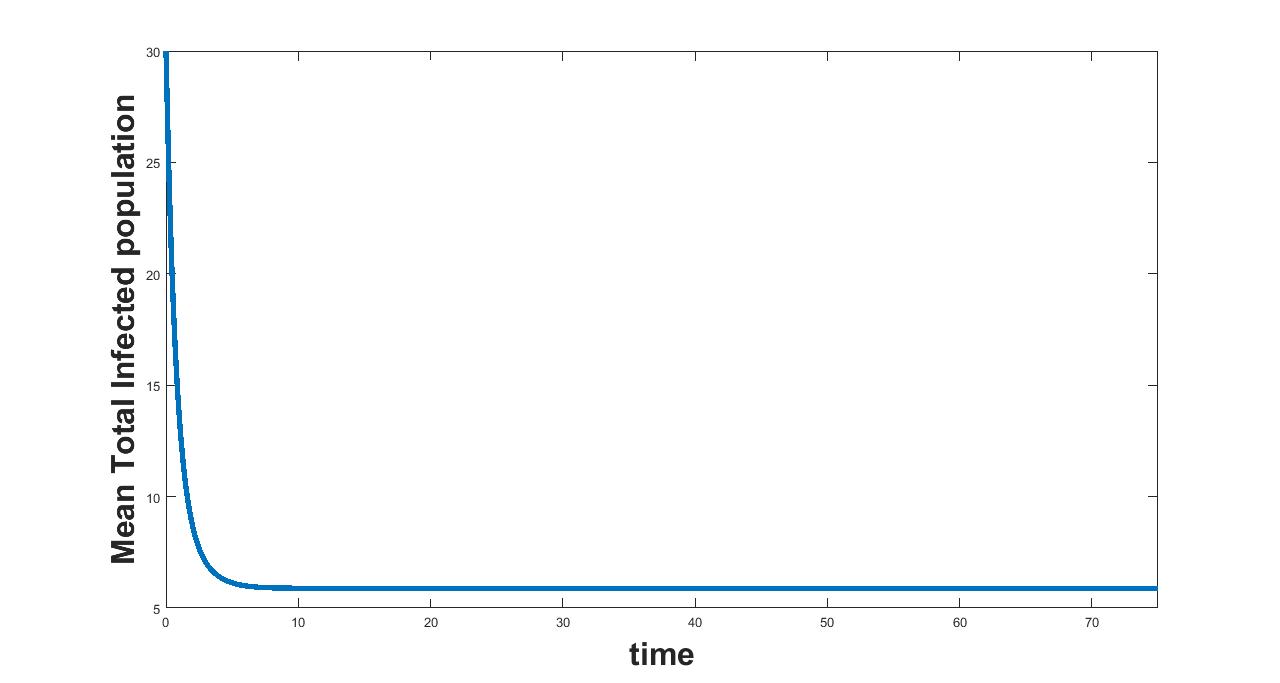}
\hspace{-.395cm}
\includegraphics[width=2.2in, height=1.8in, angle=0]{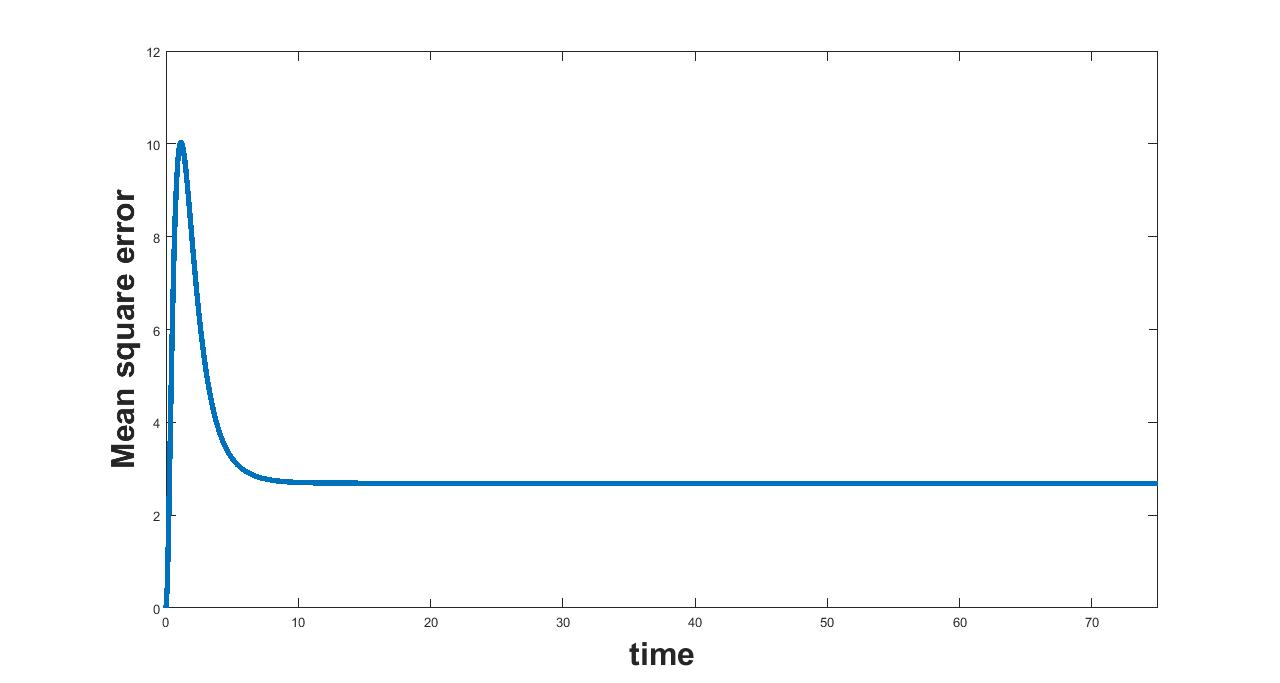}
\caption*{(b) Interval II : .000073 to 1}
\vspace{-3mm}

\vspace{5mm}
\caption{Sensitivity Analysis of $d_1$}
\label{sen_omega}
\end{center}
\end{figure}

	\subsection{Parameter $\boldsymbol{d_2}$}
		
\begin{figure}[hbt!]
\begin{center}
\includegraphics[width=2.2in, height=1.8in, angle=0]{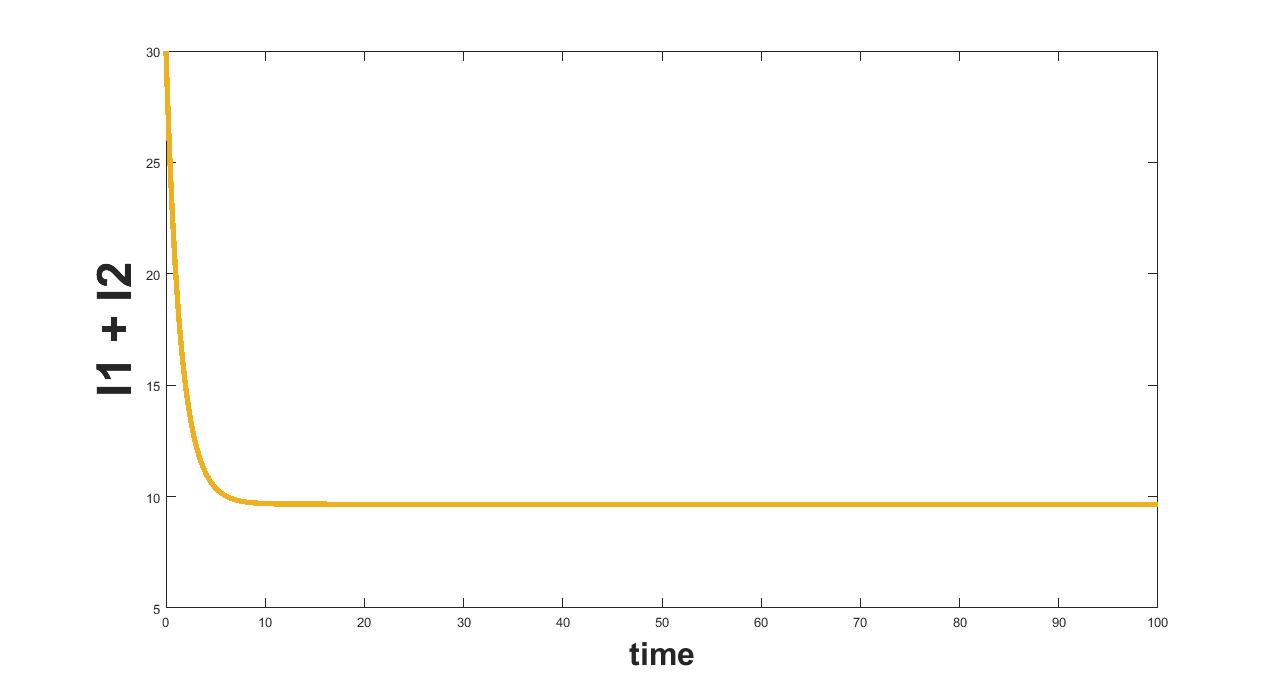}
\hspace{-.4cm}
\includegraphics[width=2.2in, height=1.8in, angle=0]{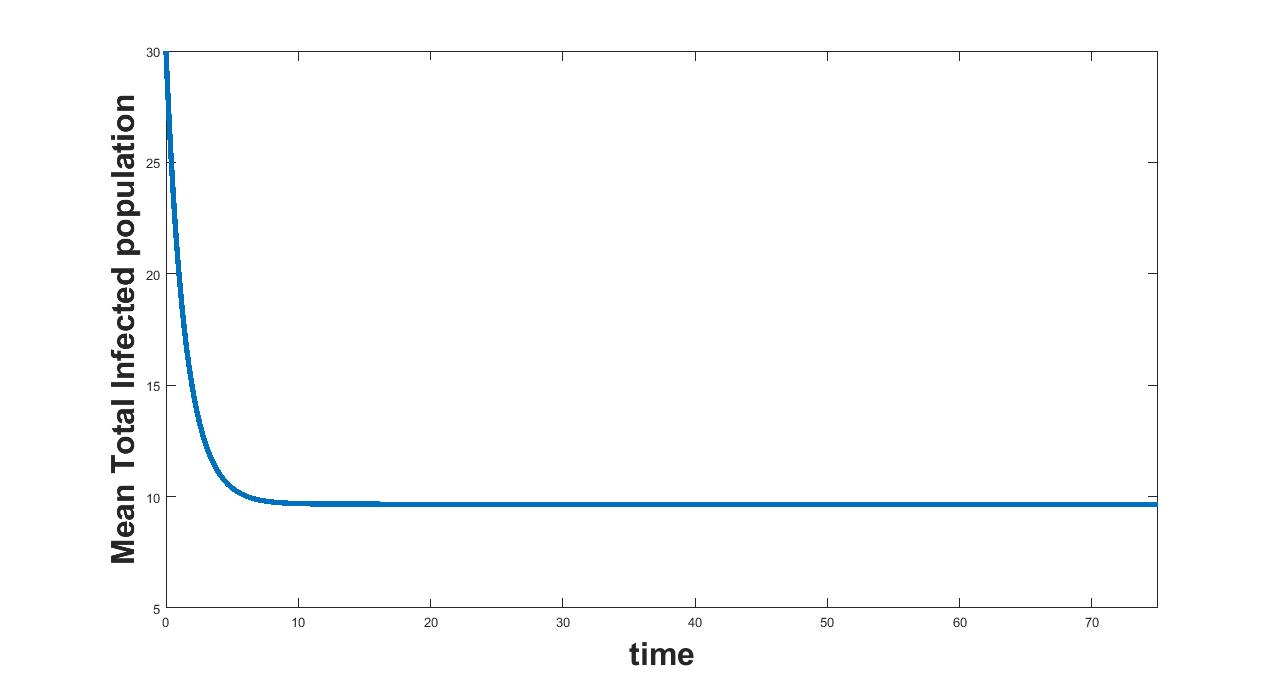}
\hspace{-.395cm}
\includegraphics[width=2.2in, height=1.8in, angle=0]{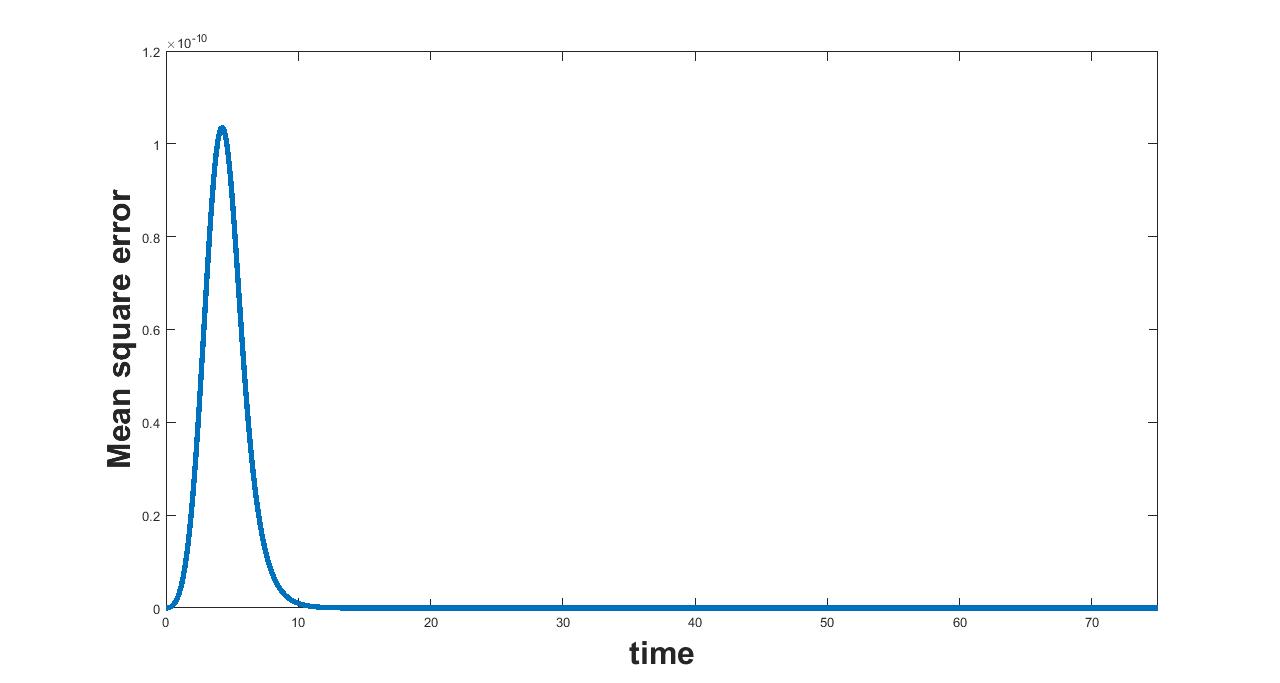}
\caption*{(a) Interval I : 0 to 0.0000913  }
\end{center}
\end{figure}

\vspace{-3mm}

\begin{figure}[hbt!]
\begin{center}
\includegraphics[width=2.2in, height=1.8in, angle=0]{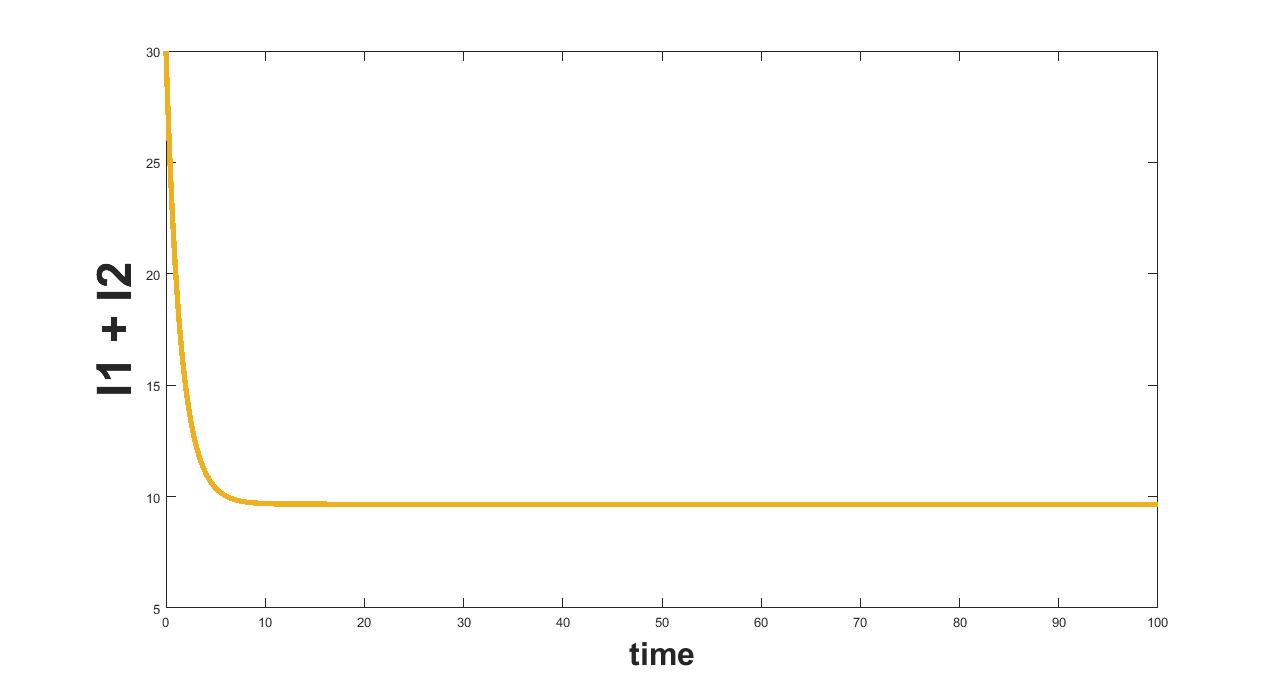}
\hspace{-.4cm}
\includegraphics[width=2.2in, height=1.8in, angle=0]{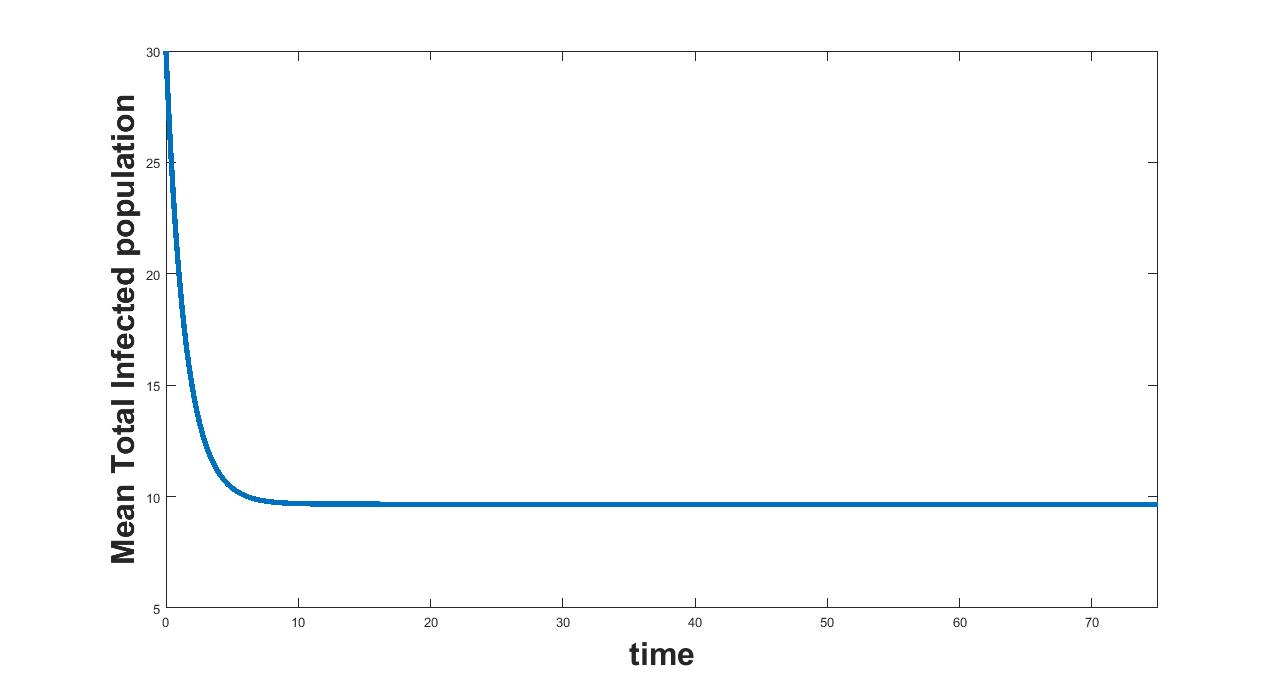}
\hspace{-.395cm}
\includegraphics[width=2.2in, height=1.8in, angle=0]{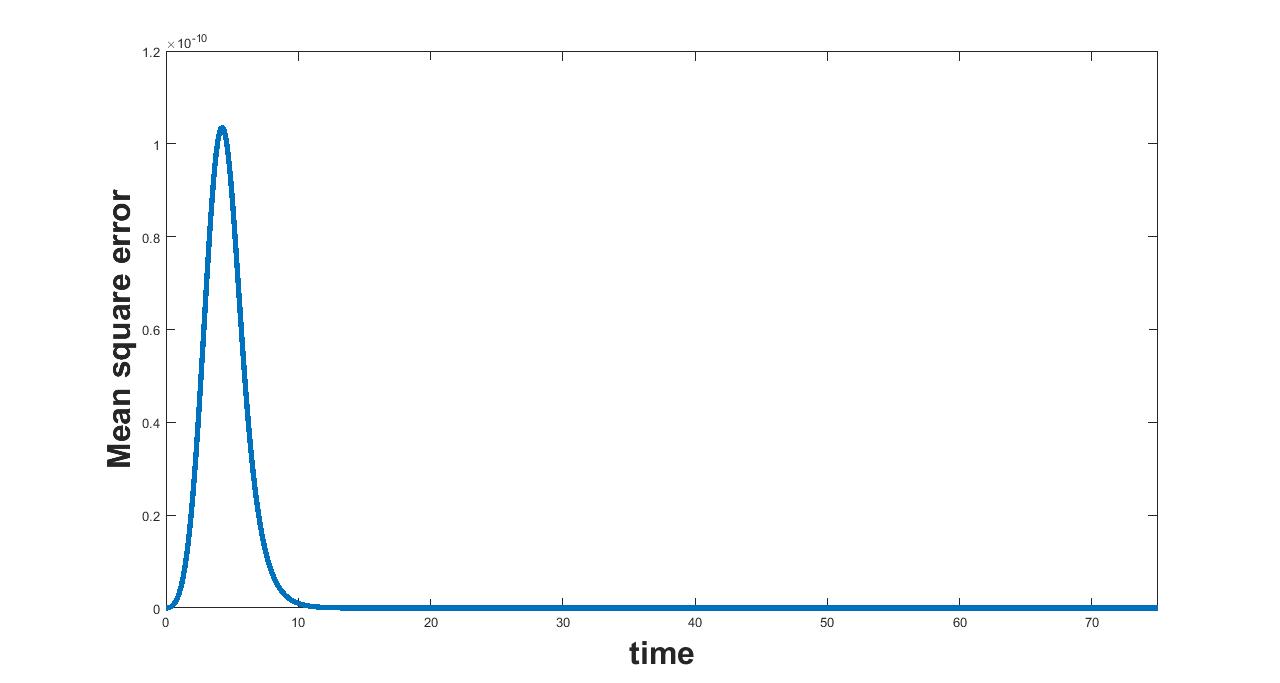}
\caption*{(b) Interval II : 0.0000913  to 2}

\vspace{5mm}
\caption{Sensitivity Analysis of $d_2$}
\label{sen_omega}
\end{center}
\end{figure}

\newpage		
		\subsection{Parameter $\boldsymbol{m}$}
		
\begin{figure}[hbt!]
\begin{center}
\includegraphics[width=2.2in, height=1.8in, angle=0]{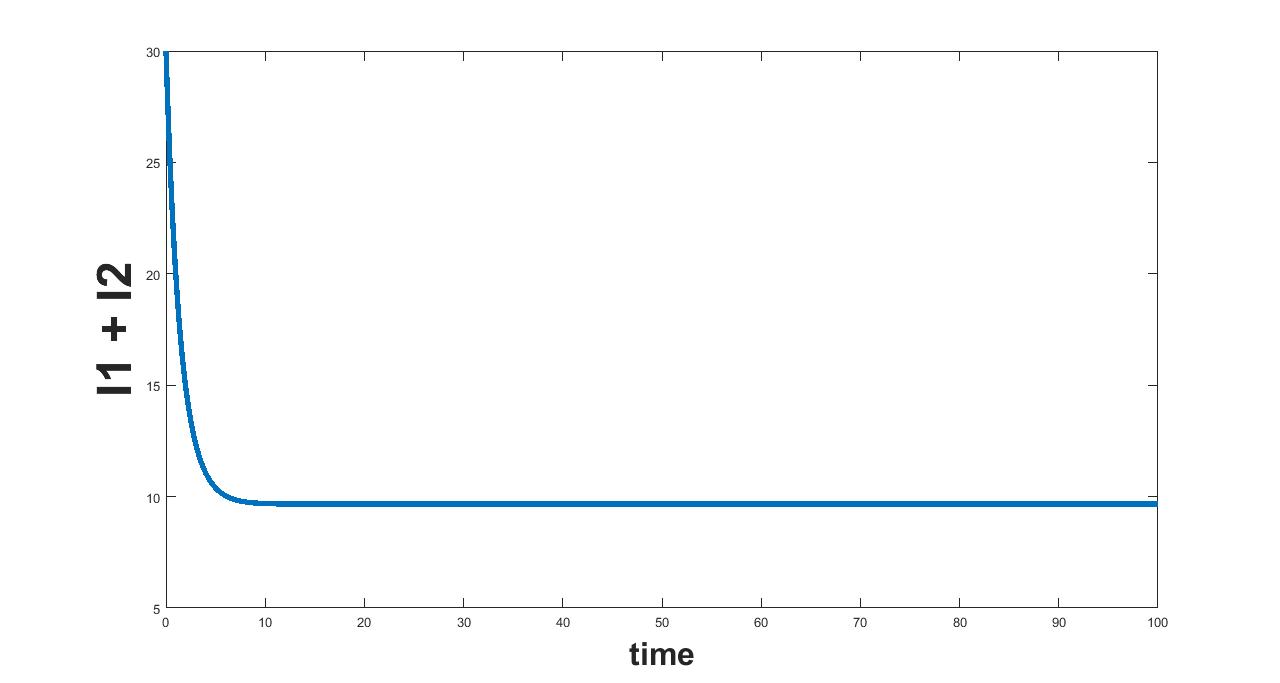}
\hspace{-.4cm}
\includegraphics[width=2.2in, height=1.8in, angle=0]{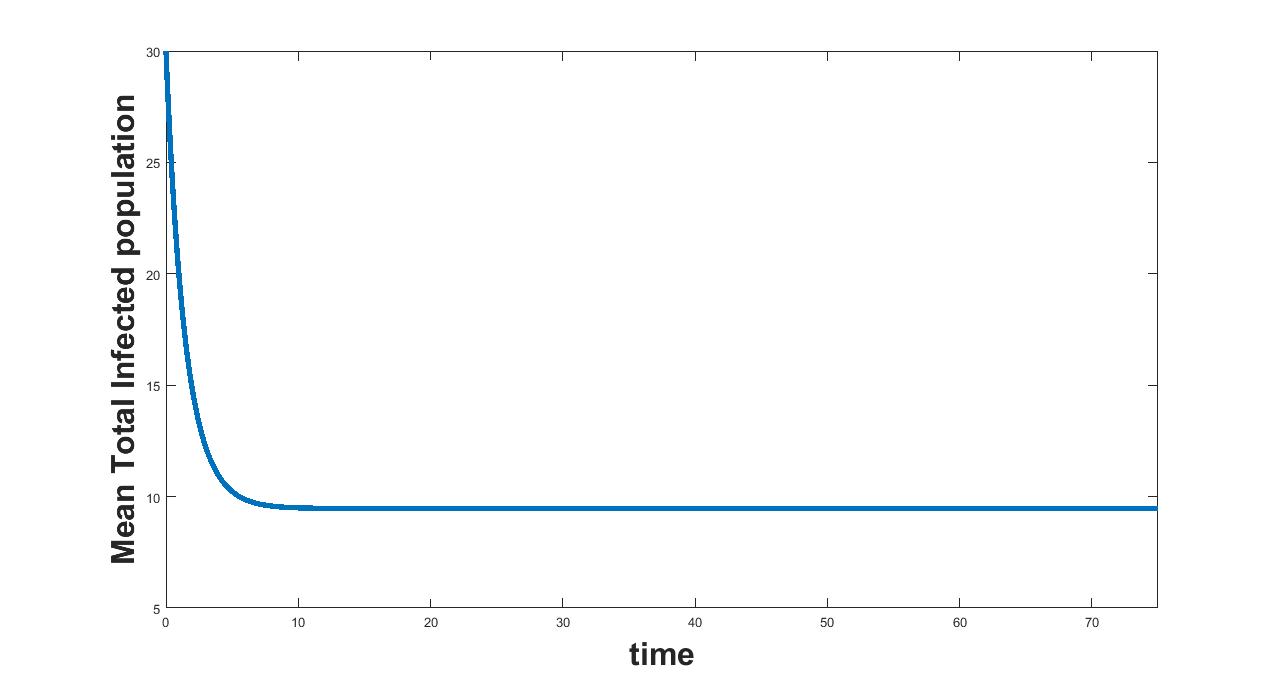}
\hspace{-.395cm}
\includegraphics[width=2.2in, height=1.8in, angle=0]{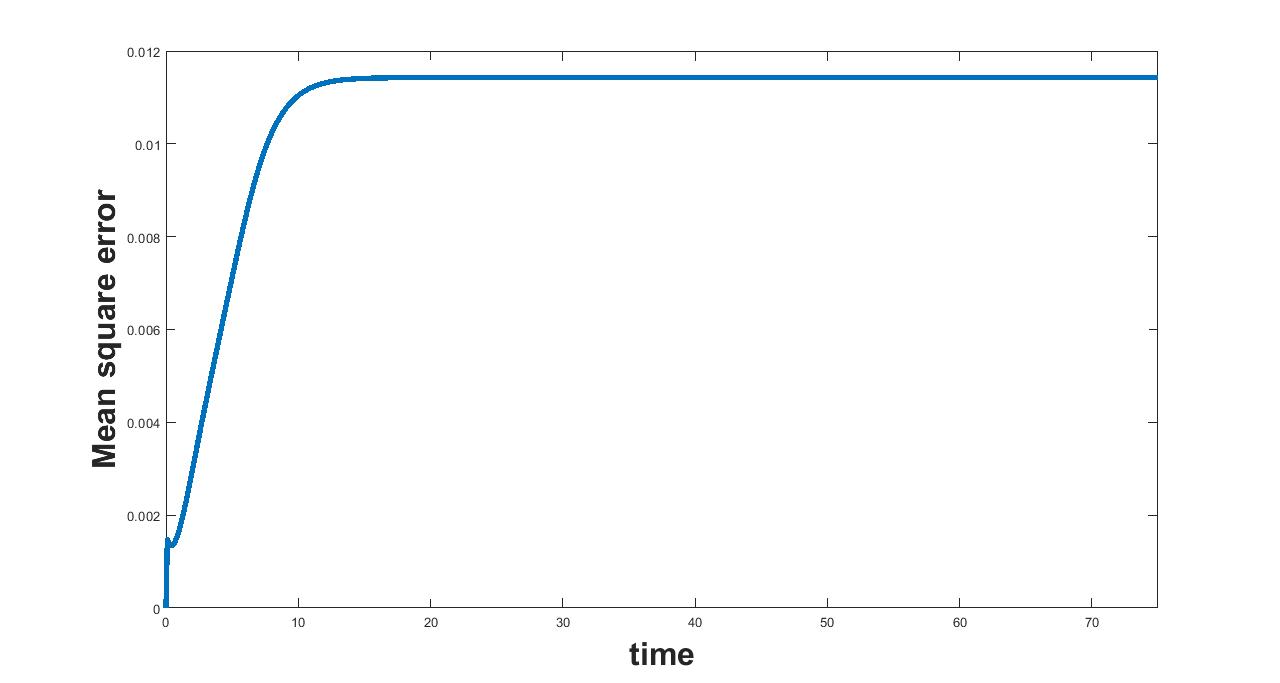}
\caption*{(a) Interval I : 0 to 0.00182}
\end{center}
\end{figure}

\vspace{-3mm}

\begin{figure}[hbt!]
\begin{center}
\includegraphics[width=2.2in, height=1.8in, angle=0]{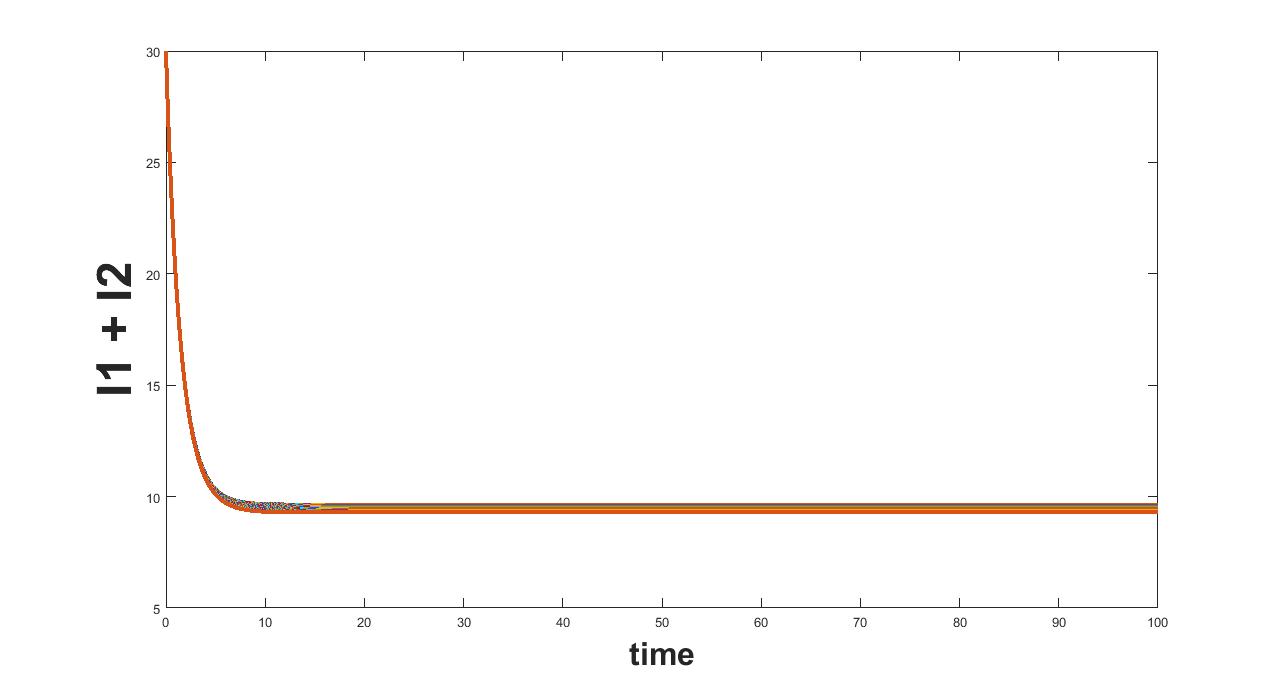}
\hspace{-.4cm}
\includegraphics[width=2.2in, height=1.8in, angle=0]{mmean}
\hspace{-.395cm}
\includegraphics[width=2.2in, height=1.8in, angle=0]{merror}
\caption*{(b) Interval II : 0.00182 to 1}

\vspace{5mm}

\caption{Sensitivity Analysis of $m$}
\label{sen_mu}
\end{center}
\end{figure}
		

\newpage
		\subsection{Parameter $\boldsymbol{u_{12}}$}
		
\begin{figure}[hbt!]
\begin{center}
\includegraphics[width=2.2in, height=1.8in, angle=0]{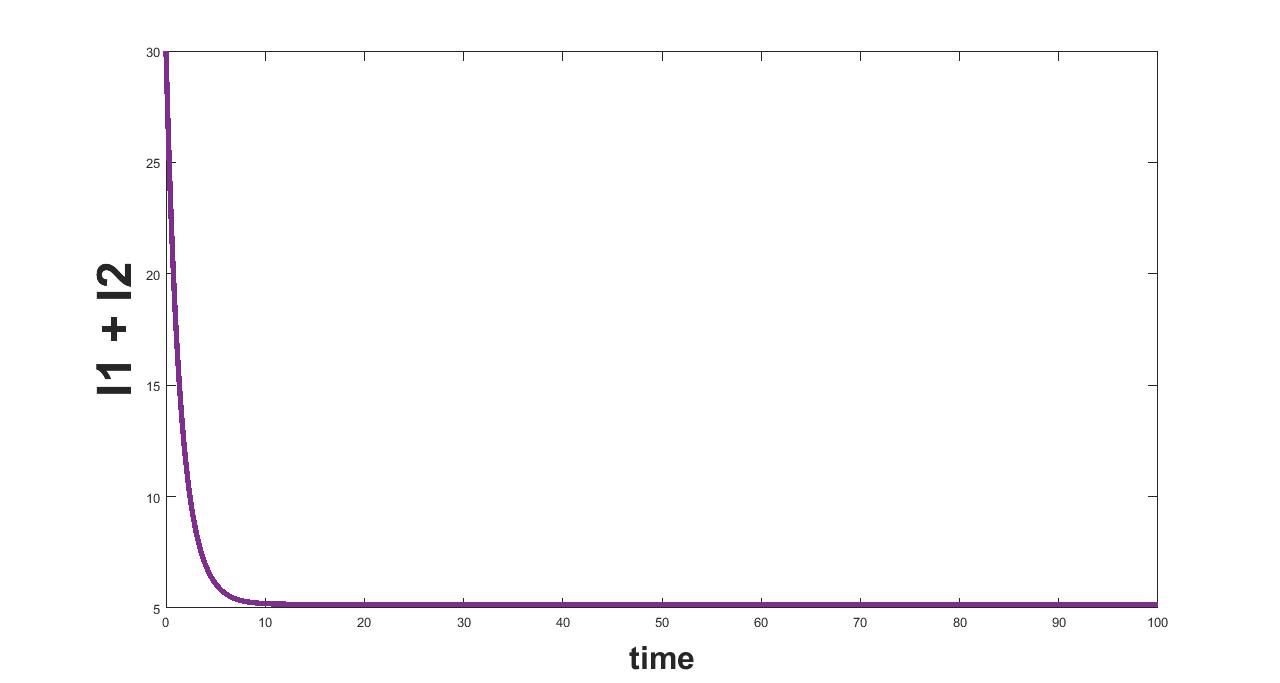}
\hspace{-.4cm}
\includegraphics[width=2.2in, height=1.8in, angle=0]{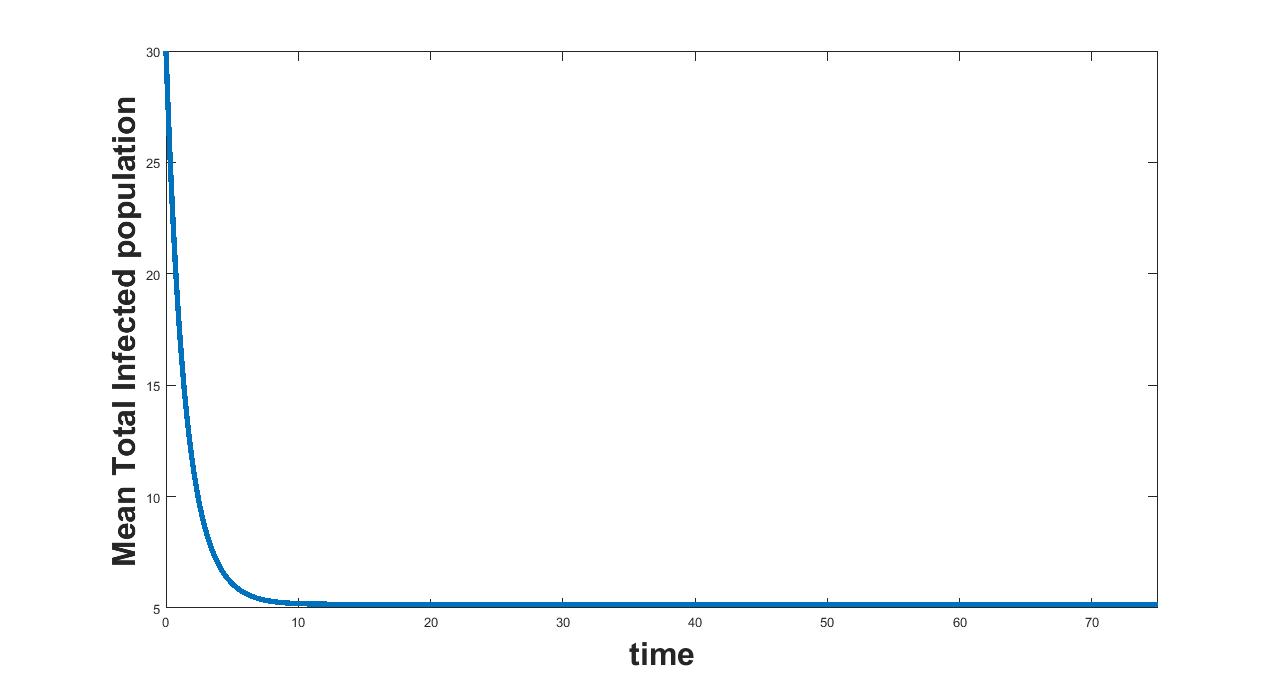}
\hspace{-.395cm}
\includegraphics[width=2.2in, height=1.8in, angle=0]{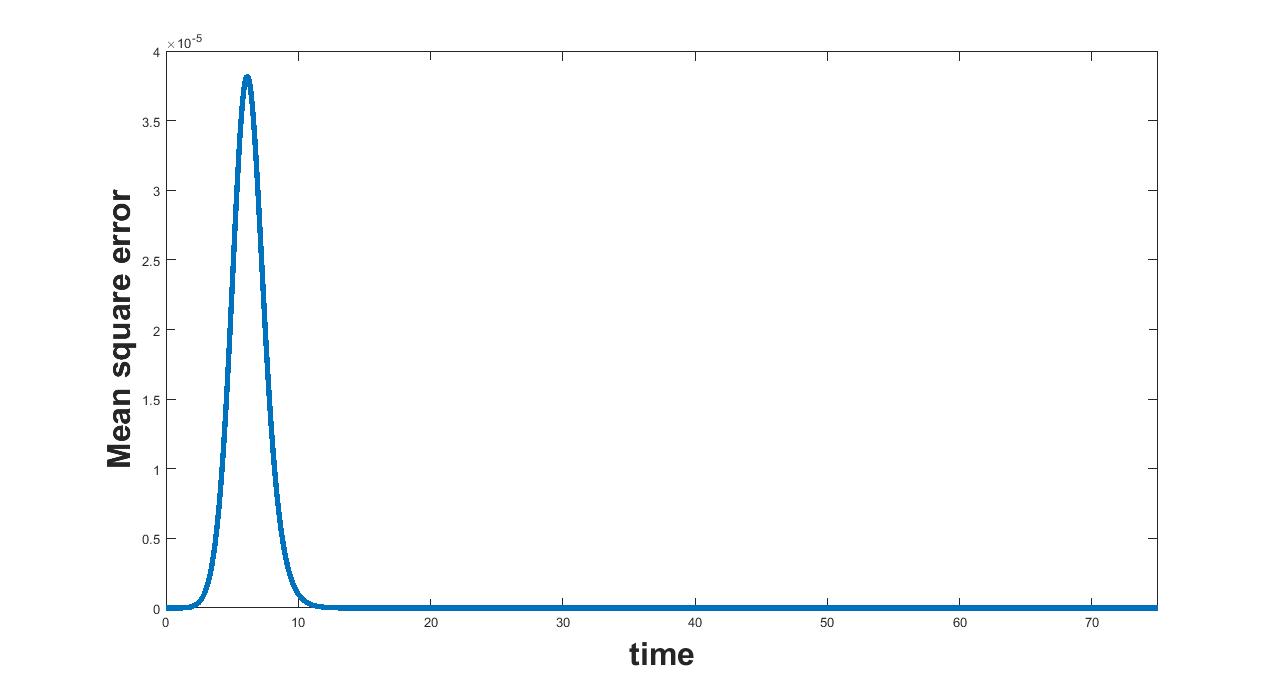}
\caption*{(a) Interval I : 0 to 0.5}
\end{center}
\end{figure}

\vspace{-3mm}

\begin{figure}[hbt!]
\begin{center}
\includegraphics[width=2.2in, height=1.8in, angle=0]{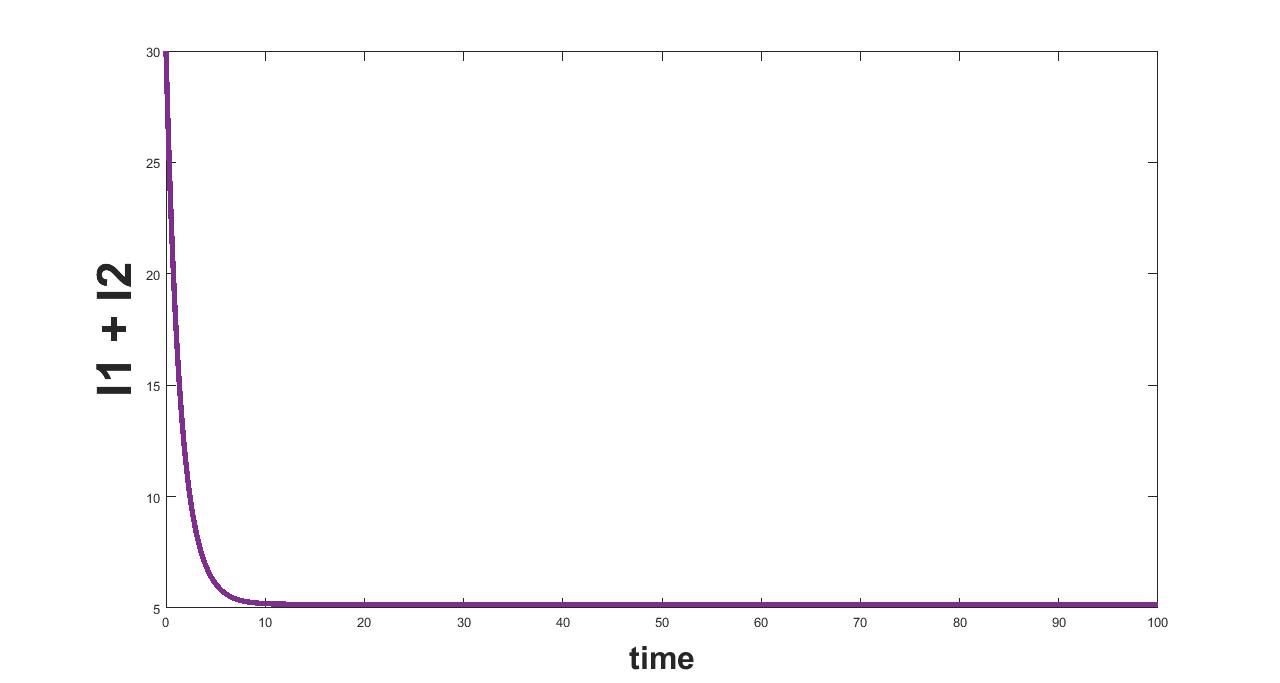}
\hspace{-.4cm}
\includegraphics[width=2.2in, height=1.8in, angle=0]{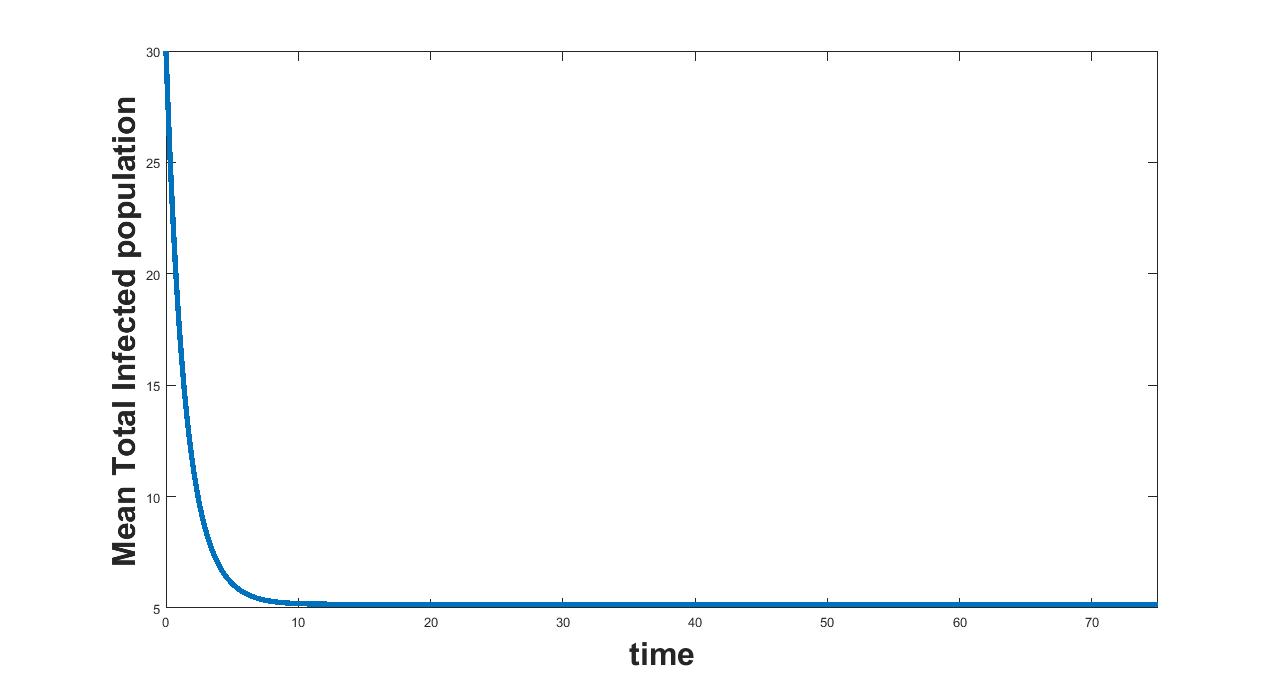}
\hspace{-.395cm}
\includegraphics[width=2.2in, height=1.8in, angle=0]{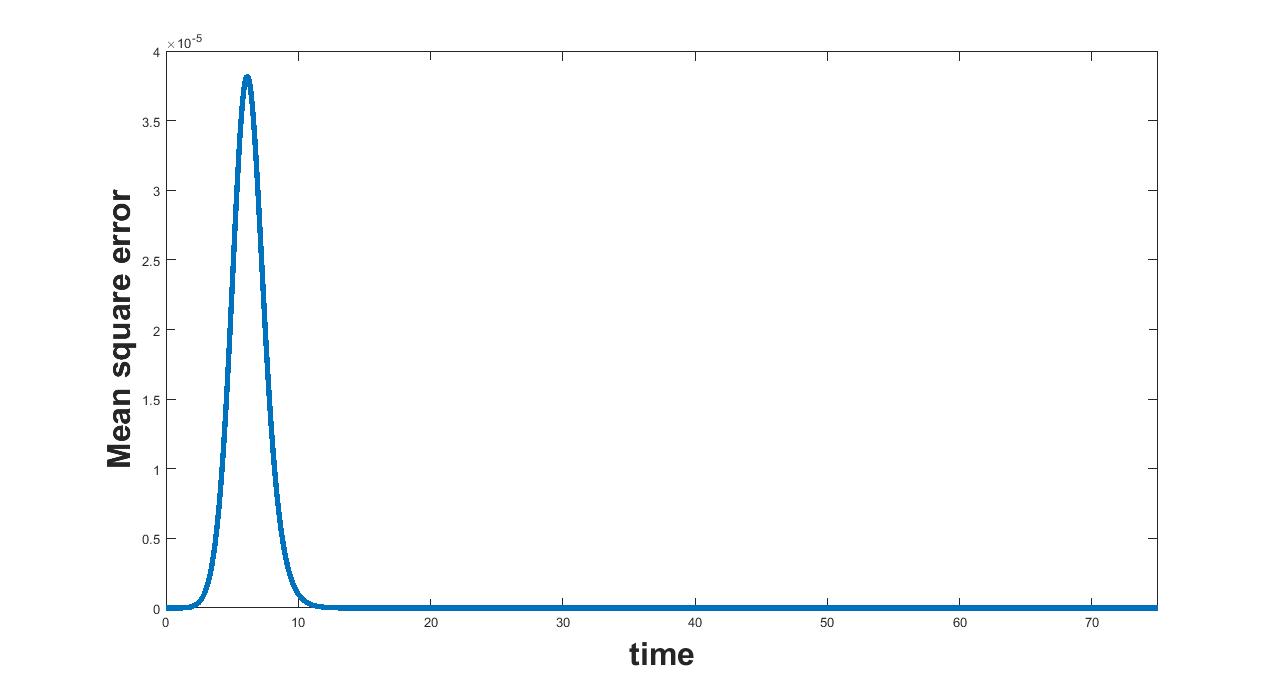}
\caption*{(a) Interval I : 0.5 to 2}
\vspace{5mm}
\caption{Sensitivity Analysis of $u_{12}$}
\end{center}
\end{figure}

\subsection{Parameter $\boldsymbol{\alpha}$}
		
\begin{figure}[hbt!]
\begin{center}
\includegraphics[width=2.2in, height=1.8in, angle=0]{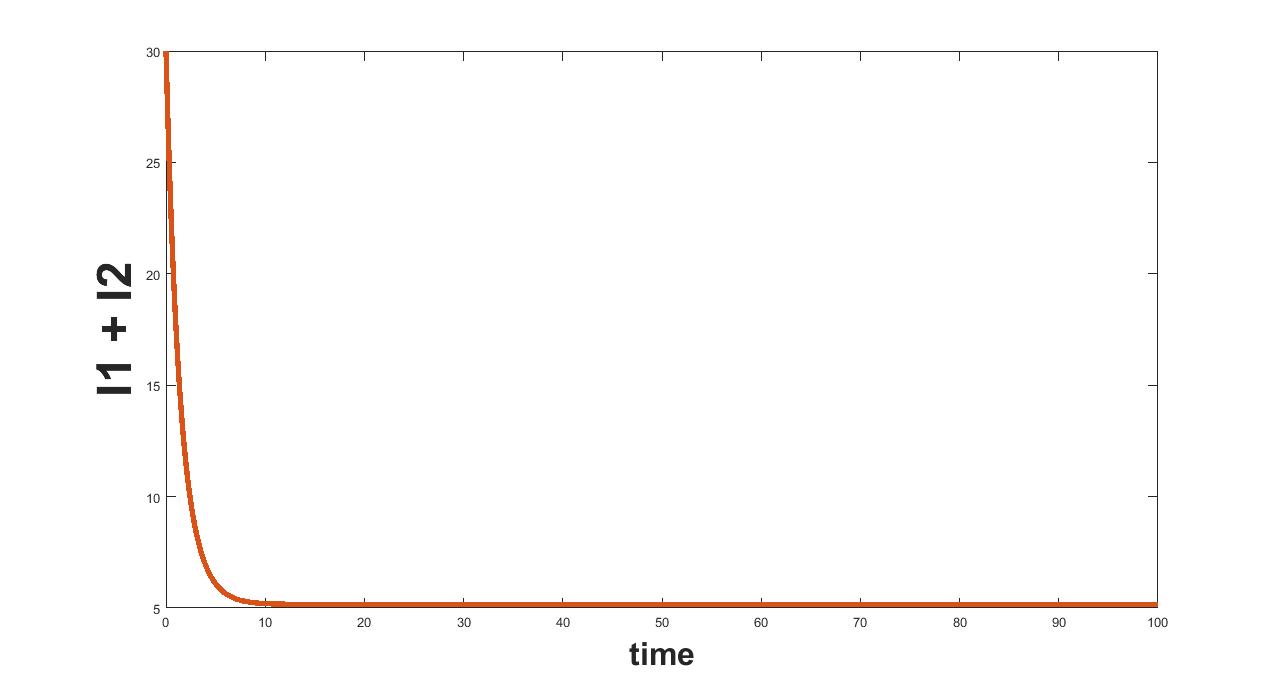}
\hspace{-.4cm}
\includegraphics[width=2.2in, height=1.8in, angle=0]{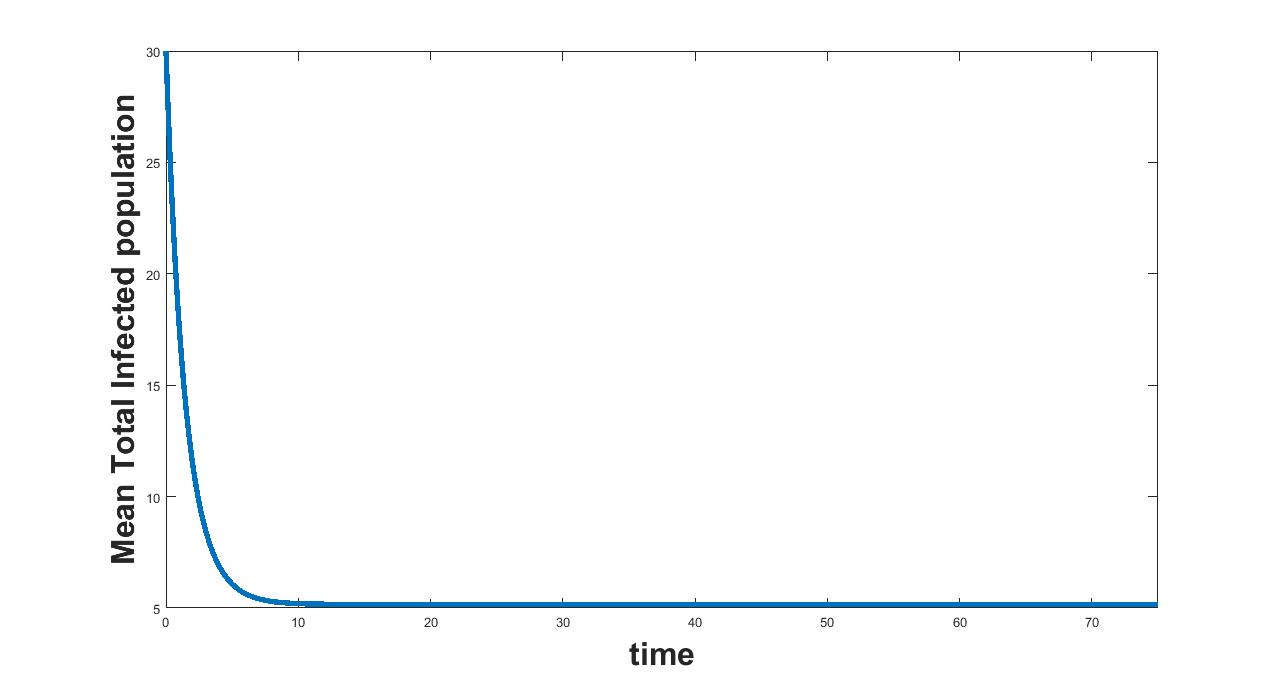}
\hspace{-.395cm}
\includegraphics[width=2.2in, height=1.8in, angle=0]{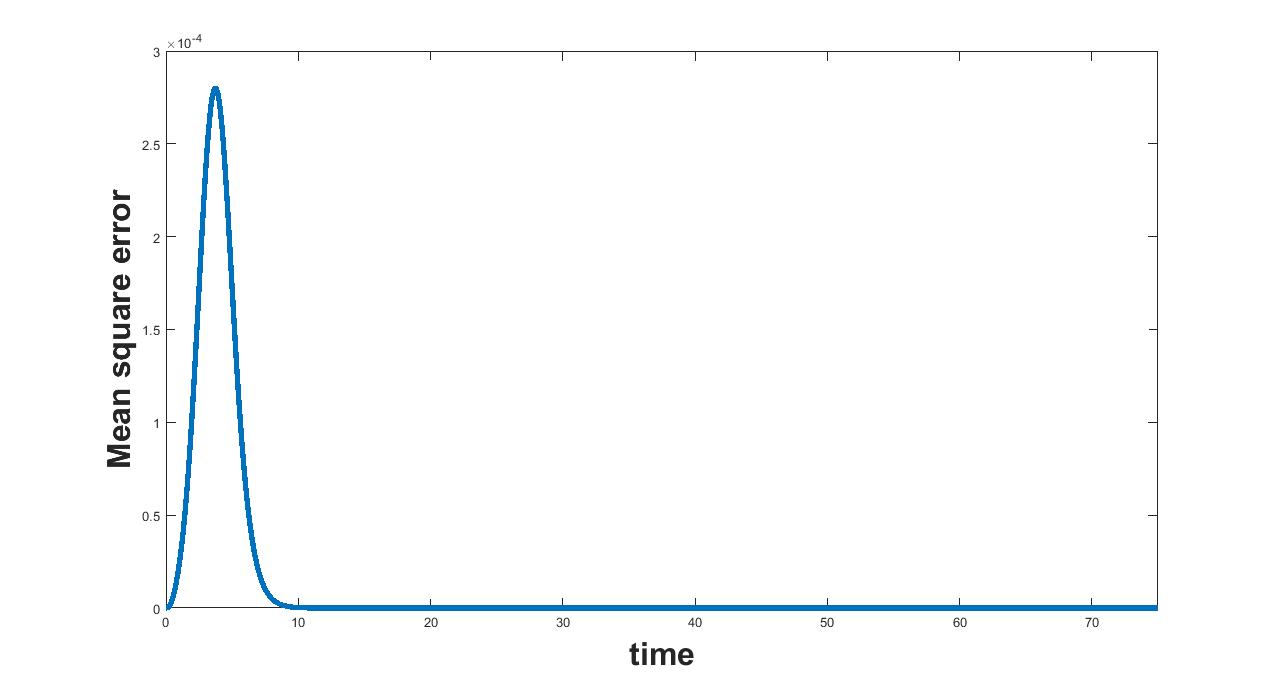}
\caption*{(a) Interval I : 0 to .5}
\end{center}
\end{figure}

\vspace{-3mm}

\begin{figure}[hbt!]
\begin{center}
\includegraphics[width=2.2in, height=1.8in, angle=0]{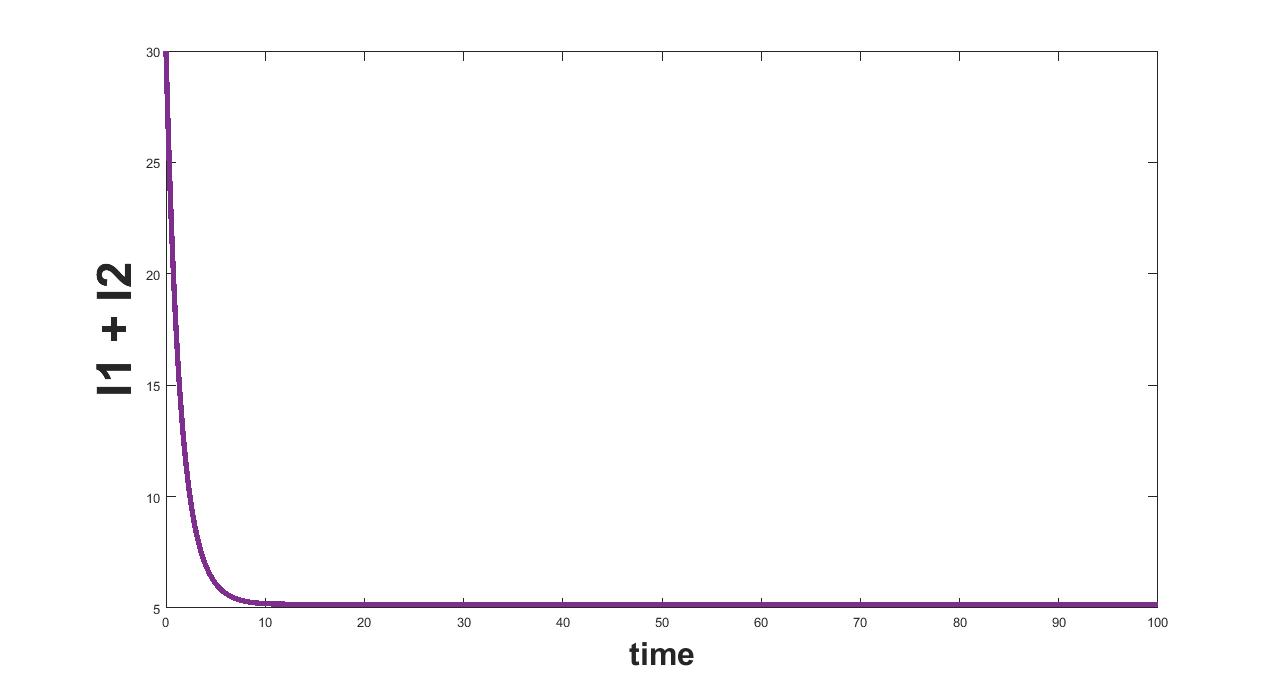}
\hspace{-.4cm}
\includegraphics[width=2.2in, height=1.8in, angle=0]{alpha_mean}
\hspace{-.395cm}
\includegraphics[width=2.2in, height=1.8in, angle=0]{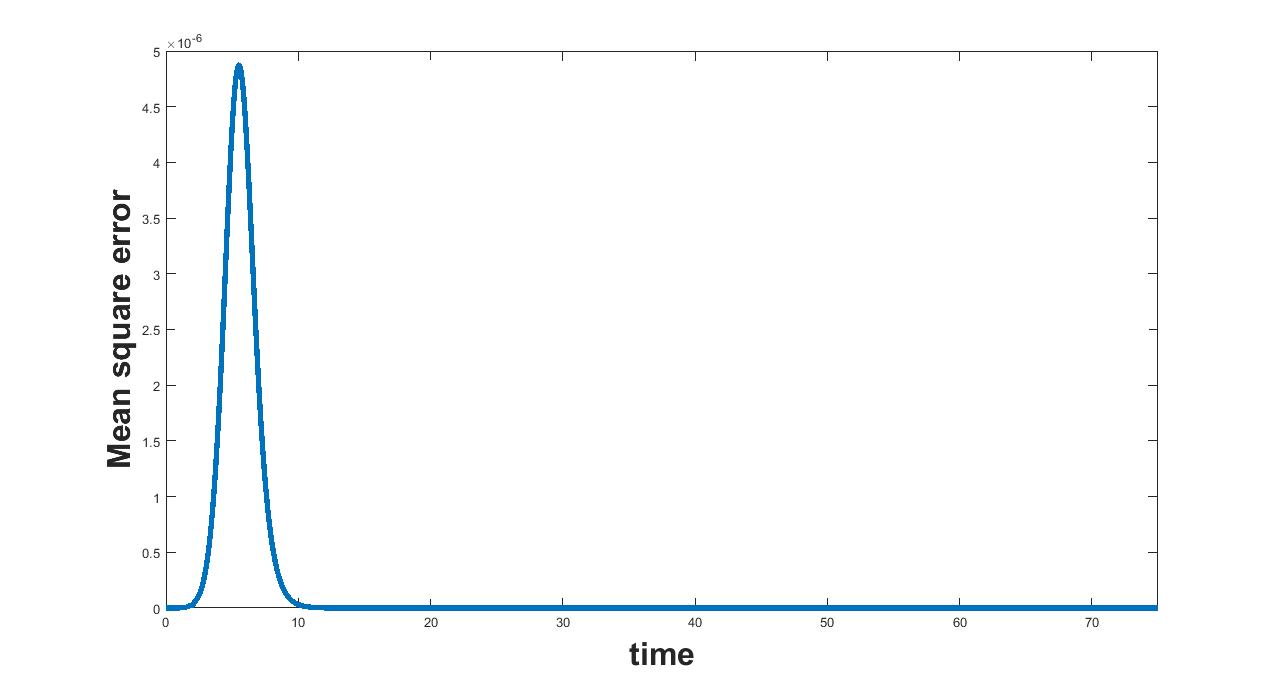}
\caption*{(b) Interval I : .5 to 2}
\vspace{5mm}
\caption{Sensitivity Analysis of $\alpha$}
\end{center}
\end{figure}
\newpage

	\subsection{Parameter $\boldsymbol{\delta_1}$}
		
\begin{figure}[hbt!]
\begin{center}
\includegraphics[width=2.2in, height=1.8in, angle=0]{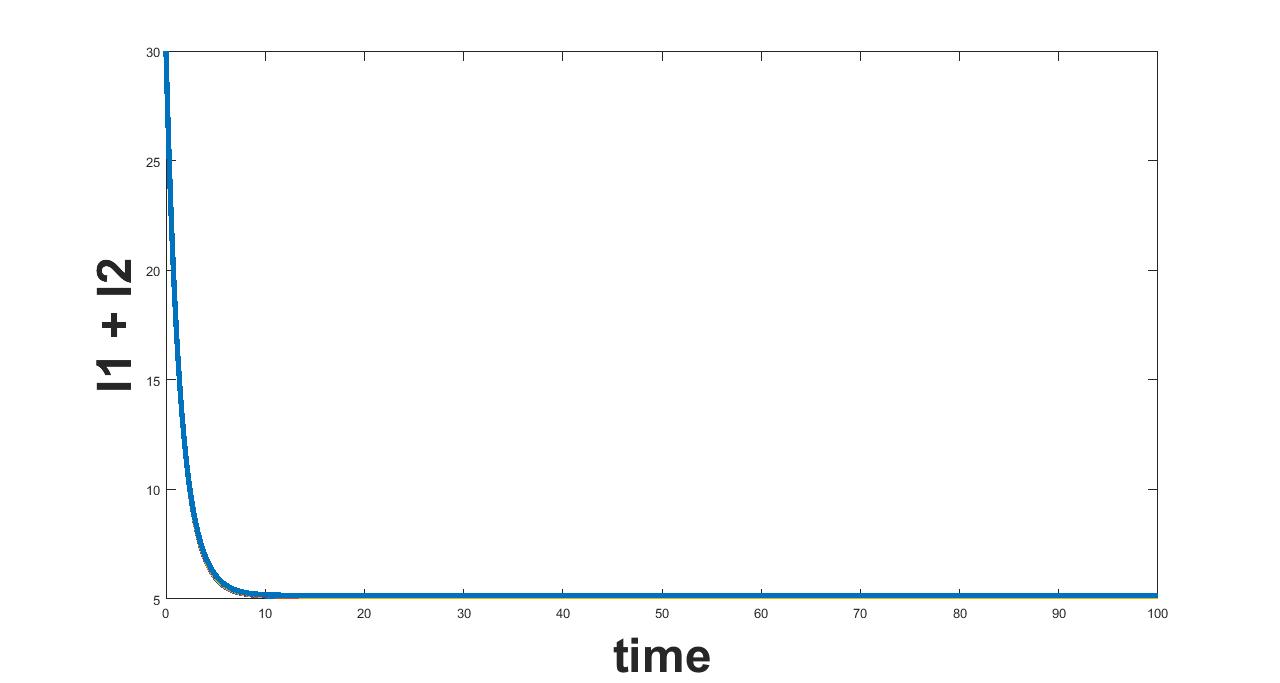}
\hspace{-.4cm}
\includegraphics[width=2.2in, height=1.8in, angle=0]{d1mean}
\hspace{-.395cm}
\includegraphics[width=2.2in, height=1.8in, angle=0]{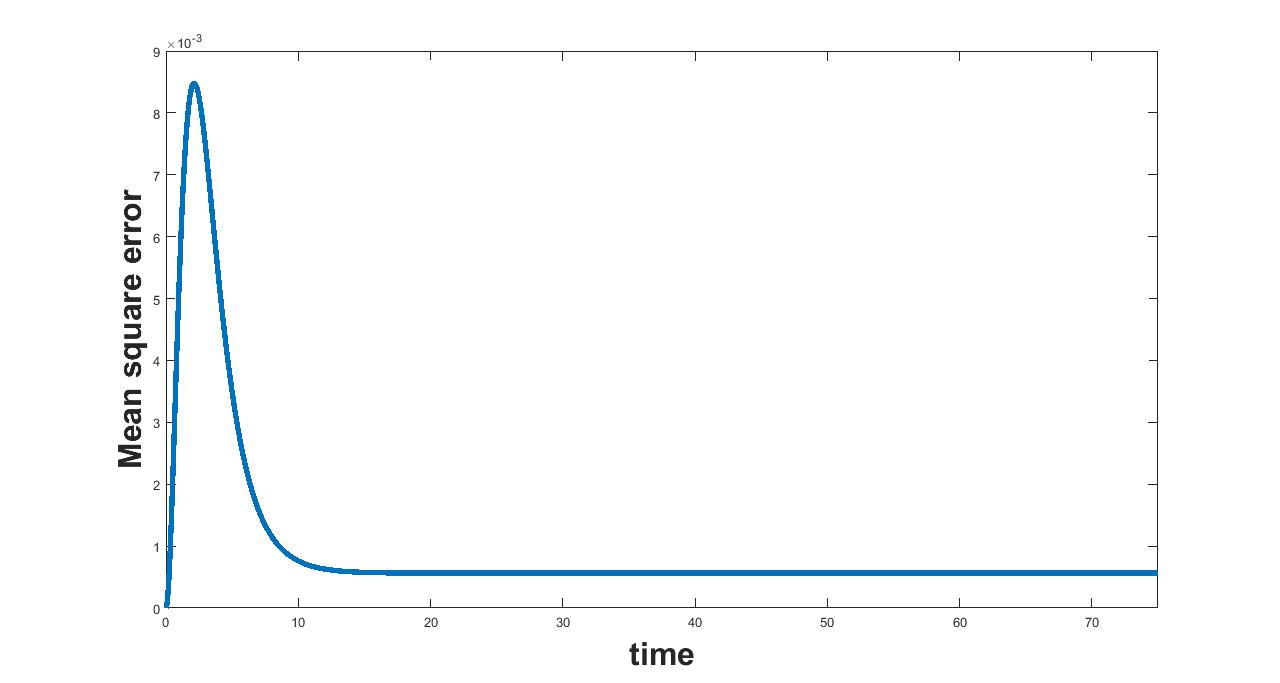}
\caption*{(a) Interval I : 0 to .0714 }
\end{center}
\end{figure}

\vspace{-3mm}

\begin{figure}[hbt!]
\begin{center}
\includegraphics[width=2.2in, height=1.8in, angle=0]{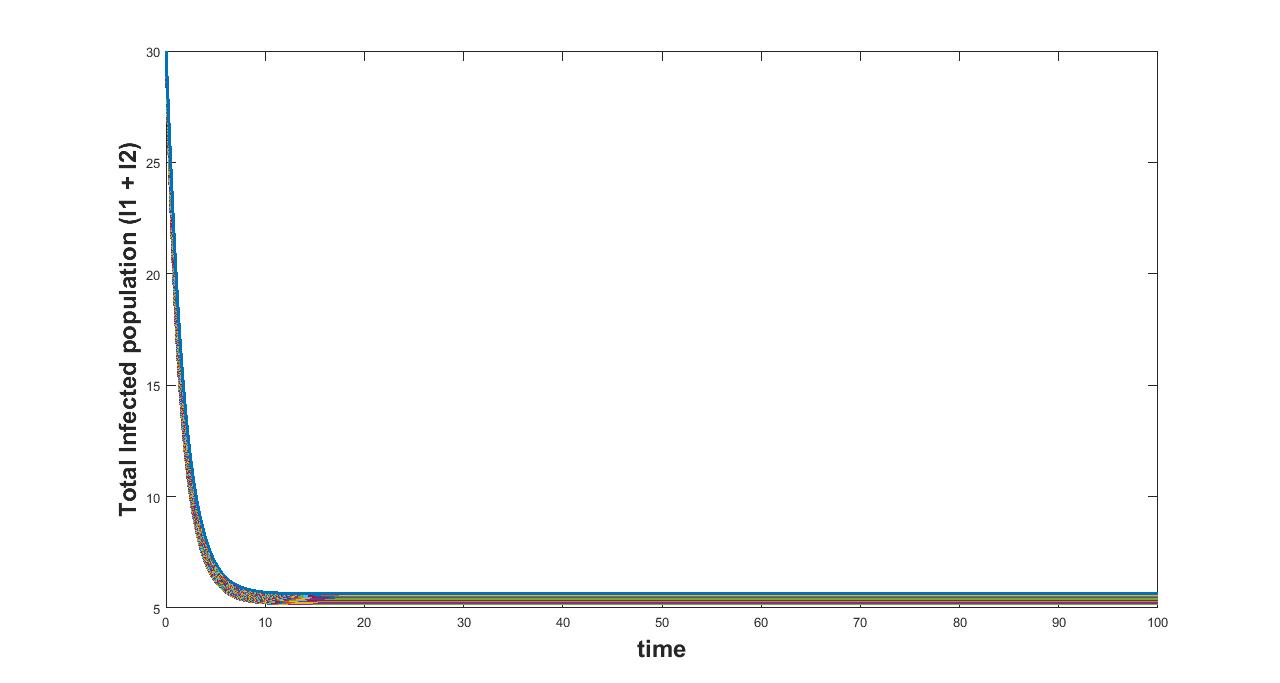}
\hspace{-.4cm}
\includegraphics[width=2.2in, height=1.8in, angle=0]{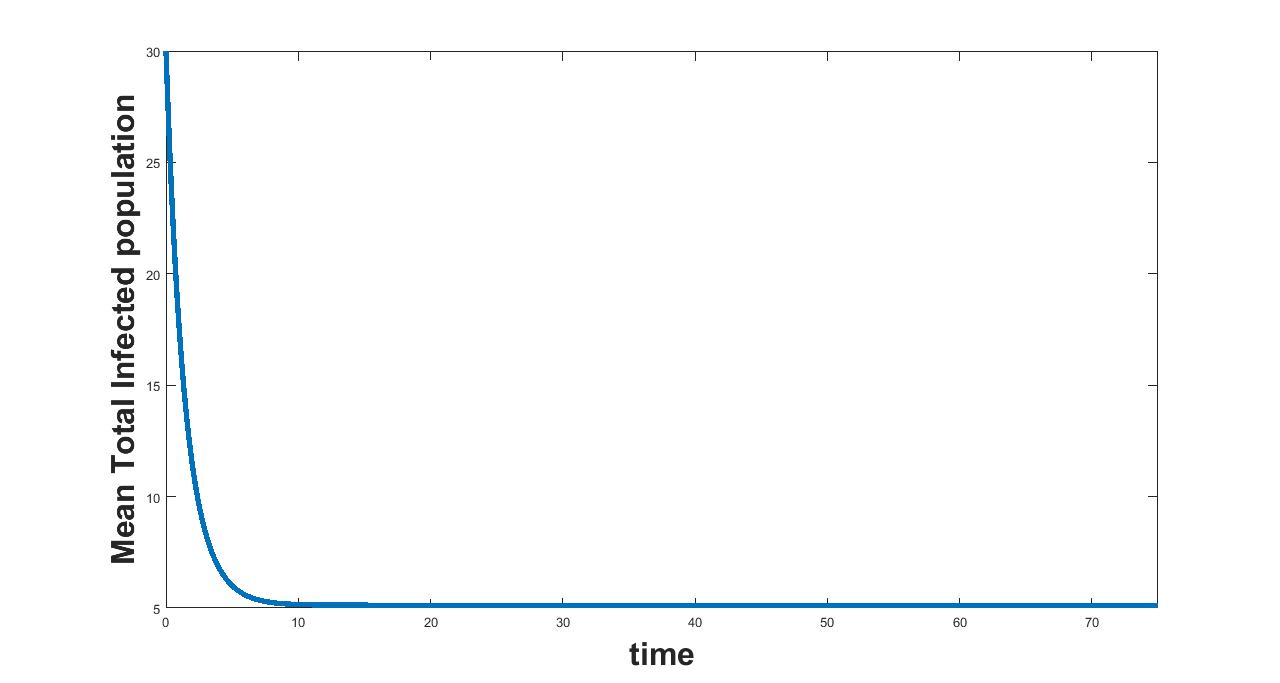}
\hspace{-.395cm}
\includegraphics[width=2.2in, height=1.8in, angle=0]{delta1erro}
\caption*{(b) Interval I : 0.0714 to 1}
\vspace{5mm}
\caption{Sensitivity Analysis of $\delta_1$}
\end{center}
\end{figure}

\newpage

	\subsection{Parameter $\boldsymbol{\delta_2}$}
		
\begin{figure}[hbt!]
\begin{center}
\includegraphics[width=2.2in, height=1.8in, angle=0]{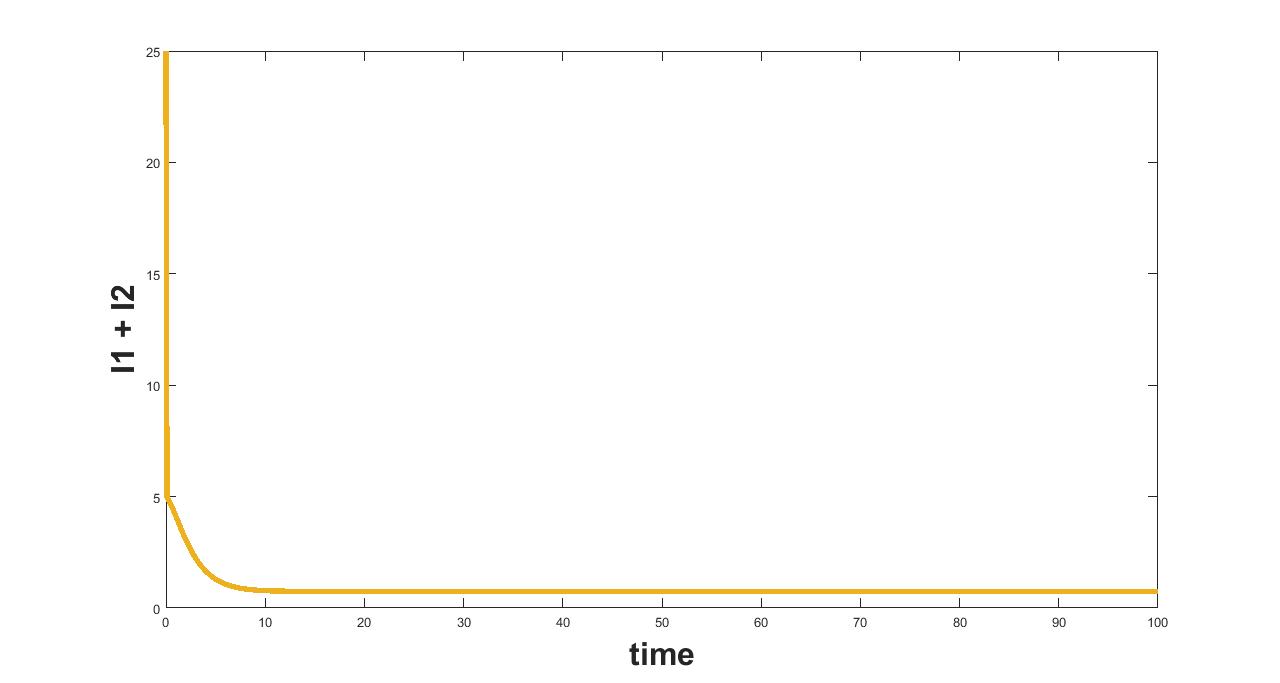}
\hspace{-.4cm}
\includegraphics[width=2.2in, height=1.8in, angle=0]{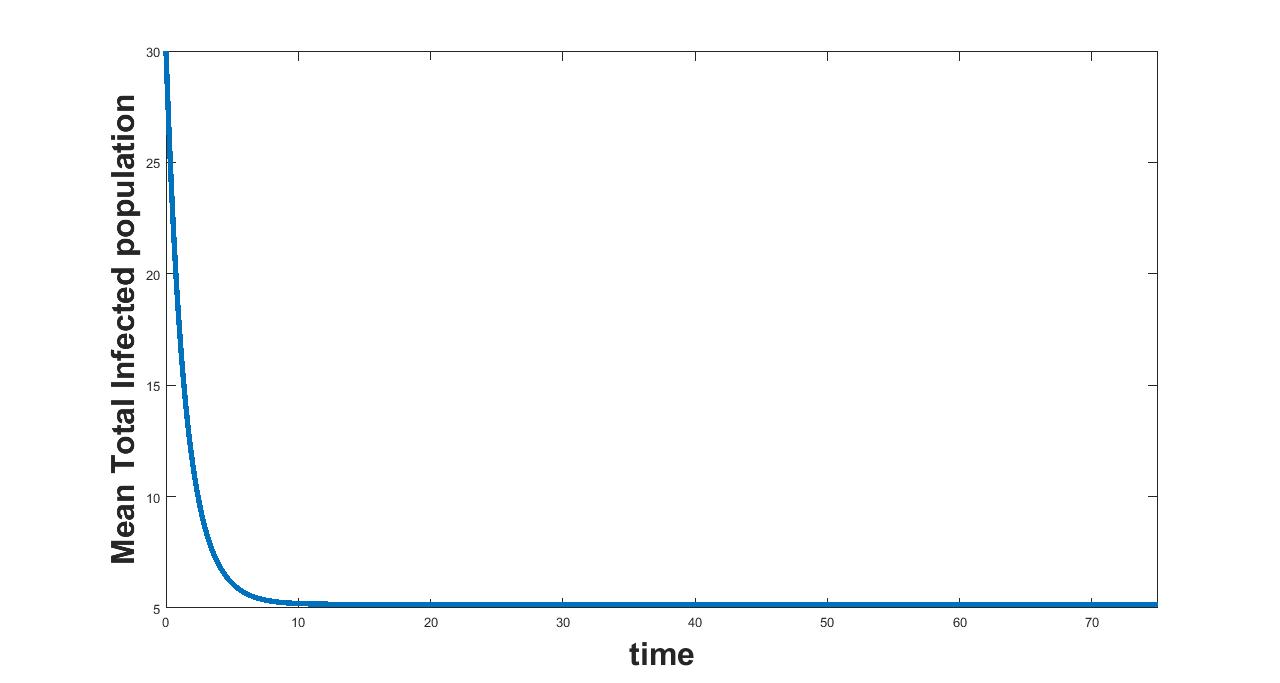}
\hspace{-.395cm}
\includegraphics[width=2.2in, height=1.8in, angle=0]{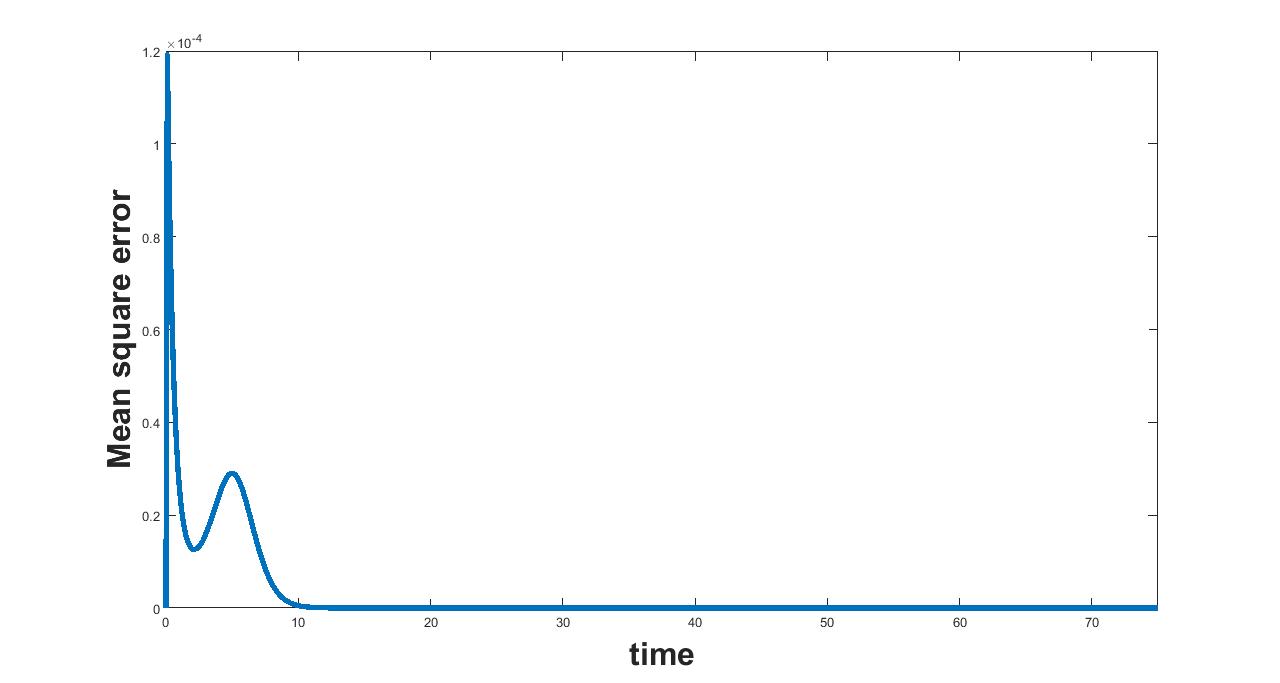}
\caption*{(a) Interval I : 0 to 0.0714}
\end{center}
\end{figure}

\vspace{-3mm}

\begin{figure}[hbt!]
\begin{center}
\includegraphics[width=2.2in, height=1.8in, angle=0]{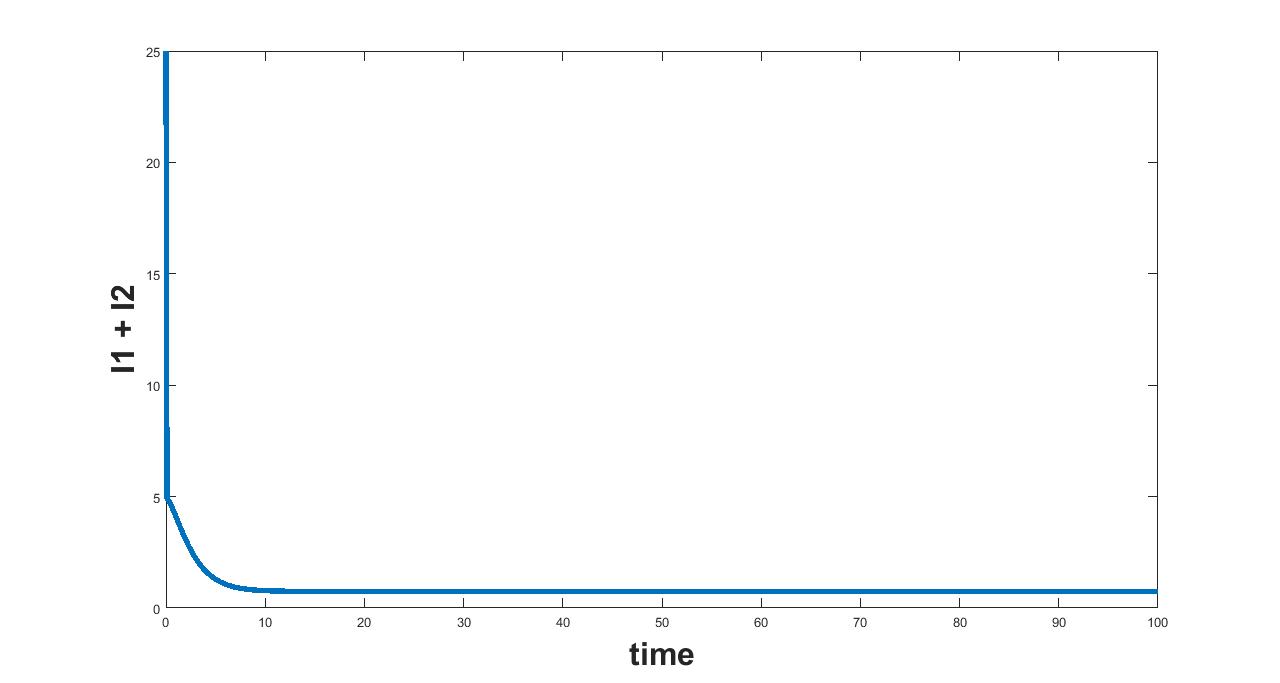}
\hspace{-.4cm}
\includegraphics[width=2.2in, height=1.8in, angle=0]{delta2mean}
\hspace{-.395cm}
\includegraphics[width=2.2in, height=1.8in, angle=0]{delta1erro}
\caption*{(b) Interval I : .0714 to 1}
\vspace{5mm}
\caption{Sensitivity Analysis of $\delta_2$}
\end{center}
\end{figure}

\end{document}